\documentclass[10pt]{article}

\usepackage{amsmath}
\usepackage{amsthm}
\usepackage{amssymb}
\usepackage{amscd}
\usepackage{amsfonts}  
\usepackage{makeidx}
\usepackage{hyperref}
\usepackage{makeidx}     
\usepackage{graphicx}    
\usepackage{multicol}    

\DeclareMathAlphabet\EuR{U}{eur}{m}{n}
\SetMathAlphabet\EuR{bold}{U}{eur}{b}{n}

\makeindex             

\newcounter{commentcounter}
\newcommand{\comment}[1]                      
{\stepcounter{commentcounter}
{\bf Comment \arabic{commentcounter}}: {\ttfamily #1} }

\newcommand{\indexnotation}[1]{\label{#1}}

\newcommand{\squarematrix}[4]                
{                                            
\left( \begin{array}{cc} #1 & #2 \\ #3 &
#4
\end{array} \right)
}

\newcommand{\calall}{{\mathcal A}\!{\mathcal L}\!{\mathcal L}}

\newcommand{\calfin}{{\mathcal F}\!{\mathcal I}\!{\mathcal N}}
\newcommand{\calcom}{{\mathcal C}\!{\mathcal O}\!{\mathcal M}}
\newcommand{\calcomop}{{\mathcal C}\!{\mathcal O}\!{\mathcal M}\!{\mathcal O}\!{\mathcal P}}
\newcommand{\calmfin}{{\mathcal M}\!{\mathcal F}\!{\mathcal I}\!{\mathcal N}}

\newcommand{\calvcyc}{{\mathcal V}\!{\mathcal C}\!{\mathcal Y}\!{\mathcal C}}
\newcommand{\caltr}{{\mathcal T}\!{\mathcal R}}

\newcommand{\calf}{{\cal F}}
\newcommand{\calg}{{\cal G}}
\newcommand{\calh}{{\cal H}}

\newcommand{\calp}{{\cal P}}

\newcommand{\calt}{{\cal T}}

\newcommand{\bbC}{{\mathbb C}}

\newcommand{\bbH}{{\mathbb H}}

\newcommand{\bbQ}{{\mathbb Q}} 
\newcommand{\bbR}{{\mathbb R}}

\newcommand{\bbZ}{{\mathbb Z}}

\newcommand{\bfK}{{\mathbf K}} 
\newcommand{\bfL}{{\mathbf L}}

\newcommand{\bfR}{{\mathbf R}}

\newcommand{\curs}{\EuR}

\newcommand{\Or}{\curs{Or}}

\newcommand{\asmb}{\operatorname{asmb}}

\newcommand{\cd}{\operatorname{cd}}

\newcommand{\colim}{\operatorname{colim}}

\newcommand{\Ext}{\operatorname{Ext}}

\newcommand{\Hei}{\operatorname{Hei}}
\newcommand{\id}{\operatorname{id}}
\newcommand{\im}{\operatorname{im}}

\newcommand{\map}{\operatorname{map}}

\newcommand{\Out}{\operatorname{Out}}
\newcommand{\pr}{\operatorname{pr}}
\newcommand{\pt}{\{\operatorname{pt.}\}}

\newcommand{\res}{\operatorname{res}}

\newcommand{\topo}{\operatorname{top}}

\newcommand{\tors}{\operatorname{tors}}
\newcommand{\vcd}{\operatorname{vcd}}

\newcommand{\Wh}{\operatorname{Wh}}

\newtheorem{theorem}{Theorem}[section]
\newtheorem{lemma}[theorem]{Lemma}

\newtheorem{definition}[theorem]{Definition}

\newtheorem{example}[theorem]{Example}
\newtheorem{remark}[theorem]{Remark}

\newtheorem{conjecture}[theorem]{Conjecture}
\newtheorem{notation}[theorem]{Notation}
\newtheorem{problem}[theorem]{Problem}
{\catcode`@=11\global\let\c@equation=\c@theorem}

\newcommand{\entry}[2]{#1, ~ \pageref{#2}}   
                                             
\newcommand{\EGF}[2]{E_{#2}(#1)}               

\newcommand{\JGF}[2]{J_{#2}(#1)}               
\newcommand{\OrGF}[2]{\Or_{#2}(#1)}               

\newcommand{\comsquare}[8]                   
{\begin{CD}
#1 @>#2>> #3\\
@V{#4}VV @VV{#5}V\\
#6 @>>#7> #8
\end{CD}
}

\begin{document}

\typeout{----------------------------  classi.tex  ----------------------------}

\title{Survey on Classifying Spaces for Families of Subgroups}
\author{
Wolfgang L\"uck\thanks{\noindent email:
lueck@math.uni-muenster.de\protect\\
www: ~http://www.math.uni-muenster.de/u/lueck/\protect\\
FAX: 49 251 8338370\protect}\\
Fachbereich Mathematik\\ Universit\"at M\"unster\\
Einsteinstr.~62\\ 48149 M\"unster\\Germany}
\maketitle


\typeout{-----------------------  Abstract  ------------------------}

\begin{abstract}
We define for a topological group $G$ and a family of subgroups
$\calf$ two versions for the classifying space for the family $\calf$,
the $G$-$CW$-version $\EGF{G}{\calf}$ and the numerable $G$-space
version $\JGF{G}{\calf}$. They agree if $G$ is discrete, or if $G$ is a
Lie group and each element in $\calf$ compact, or if $\calf$ is the family of compact
subgroups. We 
discuss special geometric models for these spaces for the family of
compact open groups in special cases such as
almost connected groups $G$ and word hyperbolic groups $G$. We deal with
the question whether there are finite models, models of finite type,
finite dimensional models. We also discuss the relevance of these
spaces for the Baum-Connes Conjecture about the topological $K$-theory
of the reduced group $C^*$-algebra, for the Farrell-Jones Conjecture
about the algebraic $K$- and $L$-theory of group rings,
for Completion Theorems and for classifying spaces for equivariant vector
bundles and for other situations.

\smallskip
\noindent
Key words: Family of subgroups, classifying spaces,\\
Mathematics Subject Classification 2000: 55R35, 57S99, 20F65, 18G99.
\end{abstract}


\typeout{--------------------   Section 0: Introduction --------------------------}
\setcounter{section}{-1}
\section{Introduction}
\label{sec: Introduction}

We define for a topological group $G$ and a family of subgroups
$\calf$ two versions for the classifying space for the family $\calf$,
the $G$-$CW$-version $\EGF{G}{\calf}$ and the numerable $G$-space
version $\JGF{G}{\calf}$. They agree, if $G$ is discrete, or if $G$ is a
Lie group and each element in $\calf$ is compact, or if each element
in $\calf$ is open, or if $\calf$ is the family of compact subgroups, but
not in general.

One motivation for the study of these classifying spaces comes from the fact that
they appear in the Baum-Connes Conjecture about the topological $K$-theory
of the reduced group $C^*$-algebra and in the Farrell-Jones Conjecture
about the algebraic $K$- and $L$-theory of group rings and that
they play a role in the formulations and constructions concerning
Completion Theorems and classifying spaces for equivariant vector
bundles and other situations. Because of the Baum-Connes Conjecture and the 
Farrell-Jones Conjecture the computation of the relevant $K$- and $L$-groups can be reduced
to the computation of certain equivariant homology groups applied to
these classifying spaces for the family of finite subgroups or the family of virtually
cyclic subgroups. Therefore it is important  to have nice geometric models for these spaces 
$\EGF{G}{\calf}$ and  $\JGF{G}{\calf}$ and in particular
for the orbit space $G\backslash \EGF{G}{\calfin}$.

The space $\EGF{G}{\calf}$ has for the family of
compact open subgroups or of finite subgroups nice geometric models 
for instance in the cases, where $G$ is an almost connected group $G$,
where $G$ is a discrete subgroup of a connected Lie group, where $G$ is a
 word hyperbolic group, arithmetic group, mapping class group, one-relator group and so on.
Models are given by symmetric spaces, Teichm\"uller spaces, outer space, Rips complexes, buildings,
trees and so on. On the other hand one can construct for any $CW$-complex $X$ a discrete group $G$ such that
$X$ and $G\backslash\EGF{G}{\calfin}$ are homotopy equivalent.

We deal with the question whether there are finite models, models of finite type,
finite dimensional models. In some sense the algebra of a discrete group $G$ 
is reflected in the geometry of
the spaces $\EGF{G}{\calfin}$. For torsionfree discrete groups $\EGF{G}{\calfin}$ 
is the same as $EG$. For discrete groups with torsion the space $\EGF{G}{\calfin}$ seems to carry 
relevant information which is not present in $EG$. For instance
for a discrete group with torsion $EG$ can never have a finite dimensional model,
whereas this is possible for $\EGF{G}{\calfin}$ and the minimal dimension is related 
to the notion of virtual cohomological dimension.

The space $\JGF{G}{\calcom}$ associated to the family of compact subgroups is
sometimes also called the classifying space for proper group actions. 
We will abbreviate it as $\underline{J}G$.
Analogously we often write $\underline{E}G$ instead of $\EGF{G}{\calcom}$.
Sometimes the abbreviation $\underline{E}G$ is used in the literature, especially
in connection with the Baum-Connes Conjecture, for the $G$-space
denoted in this article by $\underline{J}G = \JGF{G}{\calcom}$. This
does not really matter since we will show that the up to $G$-homotopy
unique $G$-map $\underline{E}G \to \underline{J}G$ is a $G$-homotopy equivalence.

A reader, who is only interested in discrete groups, can
skip Sections \ref{sec: Numerable G-Space-Version}  and 
\ref{sec: Comparision of the Two Versions} completely. 

Group means always locally compact Hausdorff topological group.
Examples are discrete groups and Lie groups but we will also consider
other groups. Space always means Hausdorff space. 
Subgroups are always assumed to be closed. Notice that isotropy groups
of $G$-spaces are automatically closed. A map is always understood to be continuous.

The author is grateful to Britta Nucinkis, Ian Leary and Guido Mislin for useful comments.

\tableofcontents


\typeout{--------------------   Section 1: G-CW-Complex-Version --------------------------}

\section{$G$-$CW$-Complex-Version}
\label{sec: G-CW-Complex-Version }

In this section we explain the $G$-$CW$-complex version of the
classifying space for a family $\calf$ of subgroups of a group $G$.


\subsection{Basics about $G$-$CW$-Complexes}
\label{subsec: Basics about G-CW-complexes} 

\begin{definition}[$G$-$CW$-complex]
\label{def: G-CW-complex}
A \emph{$G$-$CW$-complex}%
\index{G-CW-complex@$G$-$CW$-complex}
$X$ is a $G$-space together with a $G$-invariant filtration
$$\emptyset = X_{-1} \subseteq X_0 \subset X_1 \subseteq \ldots \subseteq
X_n \subseteq \ldots \subseteq\bigcup_{n \ge 0} X_n = X$$
such that $X$ carries the colimit topology%
\index{colimit topology}
with respect to this filtration
(i.e.\ a set $C \subseteq X$ is closed if and only if $C \cap X_n$ is closed
in $X_n$ for all $n \ge 0$) and $X_n$ is obtained
from $X_{n-1}$ for each $n \ge 0$ by attaching
equivariant $n$-dimensional cells, i.e.\ there exists
a $G$-pushout
$$\comsquare
{\coprod_{i \in I_n} G/H_i \times S^{n-1}}
{\coprod_{i \in I_n} q_i^n}
{X_{n-1}}
{}
{}
{\coprod_{i \in I_n} G/H_i \times D^{n}}
{\coprod_{i \in I_n} Q_i^n}
{X_n}
$$
\end{definition}

The space $X_n$ is called the \emph{$n$-skeleton}%
\index{skeleton}
of $X$. Notice that only the filtration by skeletons belongs to the
$G$-$CW$-structure but not the $G$-pushouts, only their existence is
required. An \emph{equivariant open $n$-dimensional cell}%
\index{cell!equivariant open $n$-dimensional cell}
is a $G$-component of $X_n-X_{n-1}$, i.e.\ the preimage of a path component
of $G\backslash (X_n-X_{n-1})$. The closure of an
equivariant open $n$-dimensional cell is called an
\emph{equivariant closed $n$-dimensional cell}%
\index{cell!equivariant closed $n$-dimensional cell}.
If one has chosen the $G$-pushouts in Definition \ref{def: G-CW-complex}, then the
equivariant open $n$-dimensional cells are the $G$-subspaces $Q_i(G/H_i \times (D^n-S^{n-1}))$
and the equivariant closed $n$-dimensional cells
are the $G$-subspaces $Q_i(G/H_i \times D^n)$.

\begin{remark}[Proper $G$-$CW$-complexes] \label{rem: Proper G-CW-complexes}
\em A $G$-space $X$ is called \emph{proper}%
\index{proper $G$-space}
if for each pair of points $x$ and $y$ in $X$ there are open neighborhoods
$V_x$ of $x$ and $W_y$ of $y$ in $X$ such that the closure of the subset
$\{g \in G \mid gV_x \cap W_y \not= \emptyset\}$ of $G$ is compact.
A $G$-$CW$-complex $X$ is proper
if and only if all its isotropy groups are compact
\cite[Theorem 1.23]{Lueck(1989)}. 
In particular a free $G$-$CW$-complex is always proper. However, not every free
$G$-space is proper. 
\em
\end{remark}

\begin{remark}[$G$-$CW$-complexes with open isotropy groups]
 \label{rem: G-CW-complexes with open isotropy groups} \em
Let $X$ be a $G$-space with  $G$-invariant filtration
$$\emptyset = X_{-1} \subseteq X_0 \subseteq X_1 \subseteq \ldots \subseteq
X_n \subseteq \ldots \subseteq\bigcup_{n \ge 0} X_n = X.$$ 
Then the following assertions are equivalent. i.) Every isotropy group of
$X$ is open and the filtration above yields a $G$-$CW$-structure on $X$.
ii.) The filtration above yields a (non-equivariant) $CW$-structure on $X$
such that each open cell $e \subseteq X$ and each $g\in G$ with $ge \cap e \not= \emptyset$
left multiplication with $g$ induces the identity on $e$.

In particular we conclude for a discrete group $G$ that 
a $G$-$CW$-complex $X$ is the same as a
$CW$-complex $X$ with $G$-action such that for
each open cell $e \subseteq X$ and each $g\in G$ with $ge \cap e \not= \emptyset$
left multiplication with $g$ induces the identity on $e$.
\em
\end{remark}

\begin{example}[Lie groups acting properly and smoothly on manifolds] 
\label{exa: Lie groups acting properly and smoothly on manifolds}
\em 
If $G$ is a Lie group and $M$ is a (smooth) proper $G$-manifold,
then an equivariant smooth triangulation
induces a $G$-$CW$-structure on $M$.
For the proof and for \emph{equivariant smooth triangulations}%
\index{equivariant smooth triangulation}
we refer to \cite[Theorem I and II]{Illman(2000)}. 
\end{example}

\begin{example}[Simplicial actions] \label{exa: simplicial actions}
\em Let $X$ be a simplicial complex
on which the group $G$ acts by simplicial automorphisms. Then all isotropy groups are
closed and open. Moreover, $G$ acts also on the
barycentric subdivision $X^{\prime}$ by simplicial automorphisms. The filtration of the 
barycentric subdivision
$X^{\prime}$ by the simplicial $n$-skeleton yields 
the structure of a $G$-$CW$-complex what 
is not necessarily true for $X$.
\em
\end{example}

A $G$-space
is called \emph{cocompact}%
\index{cocompact}
if $G\backslash X$ is compact. A $G$-$CW$-complex $X$ is \emph{finite}%
\index{G-CW-complex@$G$-$CW$-complex!finite}
if $X$ has only finitely many equivariant cells. A $G$-$CW$-complex
is finite if and only if it is cocompact. A $G$-$CW$-complex $X$ is
\emph{of finite type}%
\index{G-CW-complex@$G$-$CW$-complex!of finite type}
if each $n$-skeleton is finite. It is called
\emph{of dimension $\le n$}%
\index{G-CW-complex@$G$-$CW$-complex!of dimension $\le n$}
if $X = X_n$ and \emph{finite dimensional}%
\index{G-CW-complex@$G$-$CW$-complex!finite dimensional}
if it is of dimension $\le n$ for some integer $n$.
A free $G$-$CW$-complex $X$ is the same as
a $G$-principal bundle $X \rightarrow Y$ over a $CW$-complex $Y$
(see Remark \ref{rem: G-principal bundles over CW-complexes}).

\begin{theorem}[Whitehead Theorem for Families] 
\label{the: Whitehead theorem for families}
\index{Theorem!Whitehead Theorem for Families}
Let $f \colon Y \to Z$ be a $G$-map of $G$-spaces. Let $\calf$ be a
set of (closed) subgroups of $G$ which is closed under conjugation.
Then the following assertions are equivalent:

\begin{enumerate}

\item \label{the: Whitehead theorem for families: bijection}
For any $G$-$CW$-complex $X$, whose isotropy groups belong to $\calf$,
the map induced by $f$
$$f_*\colon  [X,Y]^G%
\indexnotation{[X,Y]^G}
\to [X,Z]^G,\hspace{5mm} [g] \mapsto [g \circ f]$$
between the set of $G$-homotopy classes of $G$-maps is bijective;

\item \label{the: Whitehead theorem for families: weak calf-homotopy equivalence}
For any $H \in \calf$ the map $f^H\colon Y^H \to Z^H$ is a weak homotopy equivalence%
\index{weak homotopy equivalence}
i.e.\ the map
$\pi_n(f^H,y)\colon \pi_n(Y^H,y) \to \pi_n(Z^H,f^H(y))$ is bijective
for any base point $y \in Y^H$ and $n \in \bbZ, n \ge 0$.

\end{enumerate}
\end{theorem}
\begin{proof}
\ref{the: Whitehead theorem for families: bijection} $\Rightarrow $
\ref{the: Whitehead theorem for families: weak calf-homotopy equivalence}
Evaluation at $1H$ induces for any $CW$-complex $A$ (equipped with the trivial $G$-action)
a bijection
$[G/H \times A, Y]^G \xrightarrow{\cong} [A,Y^H]$.
Hence for any $CW$-complex $A$ the map $f^H$ induces a bijection
$$(f^H)_*\colon [A,Y^H] \to [A,Z^H], \hspace{5mm} [g] \to [g \circ f^H].$$ 
This is equivalent to $f^H$ being a weak homotopy equivalence by the classical
non-equivariant Whitehead Theorem \cite[Theorem 7.17 in Chapter IV.7 on page 182]{Whitehead(1978)}.
\\[2mm]
\ref{the: Whitehead theorem for families: weak calf-homotopy equivalence} 
$\Rightarrow $
\ref{the: Whitehead theorem for families: bijection} 
We only give the proof in the case, where $Z$ is $G/G$ since this is the main important
case for us and the basic idea becomes already clear. The general case is treated for instance
in \cite[Proposition II.2.6 on  page 107]{Dieck(1987)}. We have to show
for any $G$-$CW$-complex $X$ such that two $G$-maps $f_0,f_1 \colon X \to Y$ are
$G$-homotopic provided that for any isotropy group $H$ of $X$ the $H$-fixed point set
$Y^H$ is \emph{weakly contractible}%
\index{weakly contractible}
i.e.\ $\pi_n(Y^H,y)$ consists of one element for all base points $y \in Y^H$.
Since $X$ is $\colim_{n \to \infty} X_n$ it suffices to construct inductively
over $n$ $G$-homotopies $h[n]\colon X_n \times [0,1]
\to Z$ such that $h[n]_i = f_i$ holds for $i = 0,1$ and $h[n]|_{X_{n-1} \times [0,1]} =
h[n-1]$. The induction beginning
$n = -1$ is trivial because of $X_{-1} = \emptyset$, the induction step from $n-1$ to $n
\ge 0$ done as follows. Fix a $G$-pushout
$$\comsquare
{\coprod_{i \in I_n} G/H_i \times S^{n-1}}
{\coprod_{i \in I_n} q_i^n}
{X_{n-1}}
{}
{}
{\coprod_{i \in I_n} G/H_i \times D^{n}}
{\coprod_{i \in I_n} Q_i^n}
{X_n}
$$
One easily checks that the desired $G$-homotopy $h[n]$ exists if and only if
we can find for each $i \in I$ an extension of the $G$-map
\begin{multline*} 
f_0 \circ Q_i^n \cup f_1 \circ Q_i^n \cup h[n-1] \circ (q^n_i \times \id_{[0,1]}) \colon 
\\
G/H_i \times D^n \times \{0\} \cup G/H_i \times D^n \times \{1\}
\cup G/H_i \times S^{n-1} \times [0,1]
 ~ \to ~ Y
\end{multline*}
to a $G$-map $G/H_i \times D^n \times [0,1] \to  Y$. This is the same problem as extending
the (non-equivariant) map 
$D^n \times \{0\} \cup D^n \times \{1\} \cup S^{n-1} \times [0,1] ~ \to ~ Y$, which is 
given by restricting the $G$-map above to $1H_i$, to a (non-equivariant) map
$D^n \times [0,1] \to Y^{H_i}$. Such an extension exists since $Y^{H_i}$ is weakly contractible.
This finishes the proof of Theorem \ref{the: Whitehead theorem for families}. 
\end{proof}

A $G$-map $f\colon X \rightarrow Y$ of $G$-$CW$-complexes is a
$G$-homotopy equivalence if and only if for any subgroup $H \subseteq G$ which occurs
as isotropy group of $X$ or $Y$ the induced map $f^H\colon X^H \rightarrow Y^H$
is a weak homotopy equivalence. This follows from the 
Whitehead Theorem for Families \ref{the: Whitehead theorem for families} above.

A $G$-map of $G$-$CW$-complexes
$f\colon X \rightarrow Y$ is \emph{cellular}%
\index{cellular!map}
if $f(X_n) \subseteq Y_n$ holds for all $n \ge 0$.
There is an equivariant version
of the \emph{Cellular Approximation Theorem},%
\index{Theorem!Equivariant Cellular Approximation Theorem}
namely, every $G$-map
of $G$-$CW$-complexes is $G$-homotopic to a cellular one and each $G$-homotopy
between cellular $G$-maps can be replaced by a cellular $G$-homotopy
\cite[Theorem II.2.1 on page 104]{Dieck(1987)}.


\subsection{The $G$-$CW$-Version for the Classifying Space for a Family}
\label{subsec: The G-CW-Version for the Classifying Space for a Family}

\begin{definition}[Family of subgroups]
 \label{def: family of subgroups} \em
A \emph{family $\calf$ of subgroups}%
\index{family of subgroups} 
of $G$ is a set of (closed) subgroups of $G$ which is closed under conjugation 
and finite intersections.
\end{definition}

Examples for $\calf$ are 
\\ 
\begin{tabular}{lll}
$\caltr$%
\indexnotation{caltr} 
& = & \{trivial subgroup\};
\\
$\calfin$%
\indexnotation{calfin} 
& = & \{finite subgroups\};
\\
$\calvcyc$%
\indexnotation{calvcyc} 
& = & \{virtually cyclic subgroups\};
\\
$\calcom$%
\indexnotation{calcom}
& = & \{compact subgroups\};
\\
$\calcomop$%
\indexnotation{calcomop}
& = & \{compact open subgroups\};
\\
$\calall$%
\indexnotation{calall}
& = & \{all subgroups\}.
\end{tabular}

\begin{definition}[Classifying $G$-$CW$-complex for a family of subgroups]
\label{def: Classifying G-CW- for a family of subgroups}
Let $\calf$ be a family of subgroups of $G$. A model
$\EGF{G}{\calf}$%
\indexnotation{EGF(G)(calf)}
for the \emph{classifying $G$-$CW$-complex for the family $\calf$ of subgroups}%
\index{classifying $G$-$CW$-complex for a family of subgroups}
is a $G$-$CW$-complex $\EGF{G}{\calf}$ which has the following properties:
i.) All isotropy groups of $\EGF{G}{\calf}$ belong to $\calf$. ii.) 
For any $G$-$CW$-complex $Y$, whose isotropy groups belong to $\calf$, 
there is up to $G$-homotopy precisely one $G$-map $Y \to X$. 

We abbreviate $\underline{E}G := \EGF{G}{\calcom}$%
\indexnotation{underline{E}G}
and call it the \emph{universal $G$-$CW$-complex for proper $G$-actions}.%
\index{universal G-CW-complex for proper G-actions@universal $G$-$CW$-complex for proper $G$-actions}
\end{definition}

In other words, $\EGF{G}{\calf}$ is a terminal object in the $G$-homotopy category of
$G$-$CW$-complexes, whose isotropy groups belong to $\calf$. 
In particular two models for $\EGF{G}{\calf}$ are $G$-homotopy equivalent
and for two families $\calf_0 \subseteq \calf_1$ there is
up to $G$-homotopy precisely one $G$-map $\EGF{G}{\calf_0} \to \EGF{G}{\calf_1}$.

\begin{theorem} [Homotopy characterization of $\EGF{G}{\calf}$]
\label{the: G-homotopy characterization of EGF(G)(calf)}
\index{Theorem!Homotopy characterization of $\EGF{G}{\calf}$}
Let $\calf$ be a family of subgroups.
\begin{enumerate}

\item \label{the: G-homotopy characterization of EGF(G)(calf): existence}
There exists a model for $\EGF{G}{\calf}$ for any family $\calf$;

\item \label{the: G-homotopy characterization of EGF(G)(calf): characterization}
A $G$-$CW$-complex $X$ is a model for $\EGF{G}{\calf}$ if and only if 
all its isotropy groups belong to $\calf$ and for each $H \in\calf$ the $H$-fixed point
set $X^H$ is weakly contractible.
\end{enumerate}

\end{theorem}
\begin{proof} 
\ref{the: G-homotopy characterization of EGF(G)(calf): existence}
A model can be obtained by attaching equivariant cells $G/H \times D^n$ for all 
$H \in \calf$ to make the $H$-fixed point sets weakly contractible. See for instance
\cite[Proposition 2.3 on page 35]{Lueck(1989)}.
\\[2mm]
\ref{the: G-homotopy characterization of EGF(G)(calf): characterization}
This follows from the  Whitehead Theorem for Families 
\ref{the: Whitehead theorem for families} 
applied to $f \colon X \to G/G$. 
\end{proof}

A model for $\EGF{G}{\calall}$ is $G/G$. In Section 
\ref{sec: Special Models} we will give many 
interesting geometric models for classifying spaces $\EGF{G}{\calf}$, in particular for the
case, where $G$ is discrete and  $\calf = \calfin$ or, more generally, where $G$ is a
(locally compact topological Hausdorff) group and $\calf = \calcom$. 
In some sense $\underline{E}G = \EGF{G}{\calcom}$ is the most interesting case.


\typeout{--------------------   Section 2: Numerable G-Space-Version --------------------------}

\section{Numerable $G$-Space-Version}
\label{sec: Numerable G-Space-Version}

In this section we explain the numerable $G$-space version of the
classifying space for a family $\calf$ of subgroups of group $G$.

\begin{definition}[$\calf$-numerable $G$-space]
\label{def: calf-numerable G-space}
\em A \emph{$\calf$-numerable $G$-space}%
\index{F-numerable G-space@$\calf$-numerable $G$-space}
is a $G$-space, for which there exists an open covering $\{U_i \mid i
\in I\}$ by $G$-subspaces such that there is for each $i \in I$ a
$G$-map $U_i \to G/G_i$ for some $G_i \in \calf$ and there is a
locally finite partition of unity $\{e_i \mid i \in I\}$ subordinate
to $\{U_i \mid i \in I\}$ by
$G$-invariant functions $e_i \colon X \to [0,1]$.
\end{definition}

Notice that we do not demand that the isotropy groups of a
$\calf$-numerable $G$-space belong to $\calf$. If $f \colon X \to Y$
is a $G$-map and $Y$ is $\calf$-numerable, then $X$ is also $\calf$-numerable.

\begin{lemma} \label{lem: G-CW-complexes are numerable}
Let $\calf$ be a family. Then a $G$-$CW$-complex is $\calf$-numerable if
all its isotropy groups belong to $\calf$.
\end{lemma}
\begin{proof} This follows from
the Slice Theorem for $G$-$CW$-complexes \cite[Theorem 1.37]{Lueck(1989)}
and the fact that $G\backslash X$ is a $CW$-complex 
and hence paracompact \cite{Miyazaki(1952)}. 
\end{proof}

\begin{definition}[Classifying numerable $G$-space for a family of subgroups]
\label{def: Classifying numerable G-space for a family of subgroups}
Let $\calf$ be a family of subgroups of $G$. 
A model $\JGF{G}{\calf}$%
\indexnotation{JGF(G)(calf)}
for the \emph{classifying numerable $G$-space for the family $\calf$ of subgroups}%
\index{classifying numerable $G$-space for a family of subgroups}
is a $G$-space which has  the following properties: i.)
$\JGF{G}{\calf}$ is $\calf$-numerable. ii.) For any $\calf$-numerable
$G$-space $X$ there is up to $G$-homotopy precisely one $G$-map $X \to \JGF{G}{\calf}$.

We abbreviate $\underline{J}G := \JGF{G}{\calcom}$%
\indexnotation{underline{J}G}
and call it the \emph{universal numerable $G$-space for proper $G$-actions}%
\index{universal numerable G-space for proper G-actions@universal numerable $G$-space for proper $G$-actions},
or briefly \emph{the universal space for proper $G$-actions}.%
\index{universal G-space for proper G-actions@universal $G$-space for proper $G$-actions}
\end{definition}

In other words, $\JGF{G}{\calf}$ is a terminal object in the $G$-homotopy category of
$\calf$-numerable $G$-spaces.
In particular two models for $\JGF{G}{\calf}$ are $G$-homotopy equivalent,
and for two families $\calf_0 \subseteq \calf_1$ there is
up to $G$-homotopy precisely one $G$-map $\JGF{G}{\calf_0} \to \JGF{G}{\calf_1}$.

\begin{remark}[Proper $G$-spaces] \label{rem: proper G-spaces} \em
A $\calcom$-numerable $G$-space $X$ is proper. Not every proper
$G$-space is $\calcom$-numerable. But a $G$-$CW$-complex $X$ is proper if
and only if it is $\calcom$-numerable (see Lemma \ref{lem: G-CW-complexes are numerable}).
\end{remark}

\begin{theorem} [Homotopy characterization of $\JGF{G}{\calf}$]
\label{the: G-homotopy characterization of JGF(G)(calf)}
\index{Theorem!Homotopy characterization of $\JGF{G}{\calf}$}
Let $\calf$ be a family of subgroups.
\begin{enumerate}

\item \label{the: G-homotopy characterization of JGF(G)(calf): existence}
For any family $\calf$ there exists a model for $\JGF{G}{\calf}$
whose isotropy groups belong to $\calf$;

\item \label{the: G-homotopy characterization of JGF(G)(calf): characterization}
Let $X$ be a $\calf$-numerable $G$-space. Equip $X \times X$ with the diagonal action and
let $\pr_i \colon X \times X \to X$ be the projection onto the $i$-th factor for $i = 1,2$.
Then $X$ is a model for $\JGF{G}{\calf}$ if and only if for each $H \in \calf$ there 
is $x \in X$ with $H \subseteq G_x$ and $\pr_1$ and $\pr_2$ are $G$-homotopic.

\item \label{the: G-homotopy characterization of JGF(G)(calf): necessary condition}
For $H \in \calf$ the $H$-fixed point set $\JGF{G}{\calf}^H$ is
contractible.
\end{enumerate}

\end{theorem}
\begin{proof} 
\ref{the: G-homotopy characterization of JGF(G)(calf): existence} 
A model for $\JGF{G}{\calf}$ is constructed in
\cite[Theorem I.6.6. on page 47]{Dieck(1987)} and
\cite[Appendix 1]{Baum-Connes-Higson(1994)}, namely, as the infinite join 
$\ast_{n =  1}^{\infty} Z$ for $Z = \coprod_{H \in \calf} G/H$. There
$G$ is assumed to be compact but the proof goes through for locally
compact topological Hausdorff groups. The isotropy groups are finite
intersections of the isotropy groups appearing in $Z$ and hence belong
to $\calf$.
\\[2mm]
\ref{the: G-homotopy characterization of JGF(G)(calf): characterization} 
Let $X$ be a model for the classifying space $\JGF{G}{\calf}$ for $\calf$.
Then $X \times X$ with the diagonal $G$-action is a $\calf$-numerable $G$-space.
Hence $\pr_1$ and $\pr_2$ are $G$-homotopic by the universal property.
Since for any $H \in  \calf$ the $G$-space $G/H$ is $\calf$-numerable,
there must exist a $G$-map $G/H \to X$ by the universal property of 
$\JGF{G}{\calf}$. If $x$ is the image under this map of $1H$, then
$H \subseteq G_x$.

Suppose that $X$ is a $G$-space such that for each $H \in \calf$ there 
is $x \in X$ with $H \subseteq G_x$ and $\pr_1$ and $\pr_2$ are $G$-homotopic. 
We want to show that then $X$ is a model for $\JGF{G}{\calf}$.
Let $f_0,f_1 \colon Y \to X$ be two $G$-maps. Since $\pr_i \circ (f_0 \times f_1) = f_i$
holds for $i = 0,1$, $f_0$ and $f_1$ are $G$-homotopic. It remains 
to show for any $\calf$-numerable
$G$-space $Y$ that there exists a $G$-map $Y \to X$.
Because of the universal property
of $\JGF{G}{\calf}$ it suffices to do this in the case, where
$Y = \ast_{n =  1}^{\infty} L$ for $L = \coprod_{H \in \calf} G/H$.
By assumption there is a $G$-map $L \to X$. Analogous to the construction in
\cite[Appendix 2]{Baum-Connes-Higson(1994)} one uses a $G$-homotopy
from $\pr_1$ to $\pr_2$ to construct a $G$-map
$\ast_{n =  1}^{\infty} L \to X$. 
\\[2mm]
\ref{the: G-homotopy characterization of JGF(G)(calf): necessary condition}
Restricting to $1H$ yields a bijection 
$$[G/H \times \JGF{G}{\calf}^H, \JGF{G}{\calf}]^G \xrightarrow{\cong} 
[\JGF{G}{\calf}^H, \JGF{G}{\calf}^H],$$
where we consider $X^H$ as a $G$-space with trivial $G$ action. Since
$G/H \times X^H$ is a $\calf$-numerable $G$-space, 
$[\JGF{G}{\calf}^H, \JGF{G}{\calf}^H]$ consists of one element. Hence 
$\JGF{G}{\calf}^H$ is contractible. 
\end{proof}

\begin{remark}\label{rem: Homotopy characterization of JGF{G}{calf} reversed} 
\em We do not know whether the converse of Theorem
\ref{the: G-homotopy characterization of JGF(G)(calf)}
\ref{the: G-homotopy characterization of JGF(G)(calf): necessary condition}
is true, i.e.\ whether a $\calf$-numerable
$G$-space $X$ is a model for $\JGF{G}{\calf}$ if $X^H$ is
contractible for each $H \in \calf$.
\em
\end{remark}

\begin{example}[Numerable $G$-principal bundles]
\label{exa: Numerable G-principal bundles} \em
A \emph{numerable (locally trivial) $G$-principal bundle}%
\index{numerable $G$-principal bundle}
$p\colon E \to B$ consists by definition of a $\caltr$-numerable $G$-space $E$,
a space $B$ with trivial action and a surjective $G$-map $p: E \to B$ such that the induced
map $G\backslash E \to B$ is a homeomorphism. A numerable $G$-principal bundle
$p \colon EG \to BG$ is \emph{universal}%
\index{numerable $G$-principal bundle!universal}
if and only if each numerable $G$-bundle admits a $G$-bundle map 
to $p$ and two such $G$-bundle maps are $G$-bundle homotopic. 
A numerable $G$-principal bundle is universal if and only if $E$ is
contractible. This follows from \cite[7.5 and  7.7]{Dold(1963)}.
More information about numerable $G$-principal bundles can be found for instance in
\cite[Section 9 in Chapter 4]{Husemoeller(1966)} 
\cite[Chapter I Section 8]{Dieck(1987)}.

If $p \colon E \to B$ is a universal numerable $G$-principal bundle, then
$E$ is a model for $\JGF{G}{\caltr}$. Conversely, $\JGF{G}{\caltr} \to G\backslash\JGF{G}{\caltr}$
is a model for the universal numerable $G$-principal bundle. We conclude that
a $\caltr$-numerable $G$-space $X$ is a model for $\JGF{G}{\caltr}$ if and only if
$X$ is contractible (compare Remark 
\ref{rem: Homotopy characterization of JGF{G}{calf} reversed}).
\em
\end{example}

\begin{remark}[$G$-Principal bundles over $CW$-complexes]
\label{rem: G-principal bundles over CW-complexes}
\em Let $p\colon E \to B$ be a (locally trivial) $G$-principal bundle over a $CW$-complex. 
Since any $CW$-complex is paracompact  \cite{Miyazaki(1952)}, 
it is automatically a numerable $G$-principal bundle.
The $CW$-complex structure on $B$ pulls back  to $G$-$CW$-structure on
$E$ \cite[1.25 on page 18]{Lueck(1989)}. Conversely,
if $E$ is a free $G$-CW-complex, then $E \to G\backslash E$ is a 
numerable $G$-principal bundle over a $CW$-complex by 
Lemma \ref{lem: G-CW-complexes are numerable}

The classifying bundle map from $p$ above to $\JGF{G}{\caltr} \to G\backslash\JGF{G}{\caltr}$ 
lifts to a $G$-bundle map from $p$ to $\EGF{G}{\caltr} \to G\backslash\EGF{G}{\caltr}$ 
and two such $G$-bundle maps from $p$ to  $\EGF{G}{\caltr} \to G\backslash\EGF{G}{\caltr}$ 
are $G$-bundle homotopic. Hence for $G$-principal bundles over $CW$-complexes
one can use  $\EGF{G}{\caltr} \to G\backslash\EGF{G}{\caltr}$ as the universal object.
\em
\end{remark}

We will compare the spaces $\EGF{G}{\calf}$ and $\JGF{G}{\calf}$ in Section
\ref{sec: Comparision of the Two Versions}. 
In Section \ref{sec: Special Models} we will give many 
interesting geometric models for $\EGF{G}{\calf}$ and $\JGF{G}{\calf}$
in particular in the case $\calf = \calcom$. In some sense
$\underline{J}G = \JGF{G}{\calcom}$ is the most interesting case.


\typeout{----- Section 3: Comparison of the Two Versions --------------------------}

\section{Comparison of the Two Versions}
\label{sec: Comparision of the Two Versions}

In this section we compare the two classifying spaces $\EGF{G}{\calf}$
and $\JGF{G}{\calf}$ and show that the two classifying spaces $\underline{E}G$
and $\underline{J}G$ agree up to $G$-homotopy equivalence.

Since $\EGF{G}{\calf}$ is a $\calf$-numerable space by Lemma 
\ref{lem: G-CW-complexes are numerable},
there is up to $G$-homotopy precisely one $G$-map
\begin{eqnarray}
u \colon  \EGF{G}{\calf} & \to & \JGF{G}{\calf}.
\label{map E_FG to J_FG}
\end{eqnarray}

\begin{lemma}
 \label{lem: Comparison of E_FG and J_FG}
The following assertions are equivalent for a family $\calf$ of subgroups of $G$:

\begin{enumerate}

\item \label{lem: Comparison of E_FG and J_FG: u G-homotopy equivalence} 
The map $u \colon  \EGF{G}{\calf}  \to  \JGF{G}{\calf}$ defined in 
\ref{map E_FG to J_FG}  is a $G$-homotopy equivalence;

\item \label{lem: Comparison of E_FG and J_FG: existence of a G-homotopy equivalence} 
The $G$-spaces $\EGF{G}{\calf}$ and $\JGF{G}{\calf}$ are $G$-homotopy equivalent;

\item \label{lem: Comparison of E_FG and J_FG: (J_FG) G-homotopc to calf-G-CW-complex}
The $G$-space $\JGF{G}{\calf}$ is $G$-homotopy equivalent to a $G$-$CW$-complex, whose
isotropy groups belong to $\calf$;

\item \label{lem: Comparison of E_FG and J_FG: map J?F(G) to G-CW-complex X with iso in calf}
There exists a $G$-map $\JGF{G}{\calf} \to Y$ to a $G$-$CW$-complex $Y$, 
whose isotropy groups belong to $\calf$;

\end{enumerate}
\end{lemma}
\begin{proof}
This follows from the
universal properties of $\EGF{G}{\calf}$ and $\JGF{G}{\calf}$.
\end{proof}

\begin{lemma} \label{lem: E_FG = J_FG if all H in calf are open}
Suppose either that every element $H \in\calf$ is an open (and closed) subgroup of $G$ or that
$G$ is a Lie group and $\calf \subseteq \calcom$.
Then the map $u \colon  \EGF{G}{\calf}  \to  \JGF{G}{\calf}$ defined in 
\ref{map E_FG to J_FG}  is a $G$-homotopy equivalence.
\end{lemma}
\begin{proof} We have to inspect the construction in 
\cite[Lemma 6.13 in Chapter I on page 49]{Dieck(1987)}
and will use the same notation as in that paper. Let $Z$ be a $\calf$-numerable $G$-space.
Let $X = \coprod_{H \in \calf} G/H$. Then $\ast_{n = 1}^{\infty} X$  is a model for
$\JGF{G}{\calf}$ by \cite[Lemma 6.6 in Chapter I on page 47]{Dieck(1987)}. 
We inspect the construction of a $G$-map $f\colon Z \to \ast_{n = 1}^{\infty} X$.
One constructs a countable covering $\{U_n \mid n = 1,2, \ldots\}$ of $Z$ by $G$-invariant
open subsets of $Z$ together with a locally finite subordinate partition of unity
$\{v_n\mid n = 1,2, \ldots\}$ by $G$-invariant functions $v_n \colon Z \to [0,1]$ and
$G$-maps $\phi_n \colon U_n \to X$. Then one obtains a $G$-map
$$f\colon Z \to \ast_{n = 1}^{\infty} X, \hspace{5mm} 
z \mapsto (v_1(z)\phi_1(z), v_2(z)\phi_2(z), \ldots ),$$
where $v_n(z)\phi_n(z)$ means $0x$ for any $x \in X$ if $z \not\in U_n$.  Let 
$i_k \colon \ast_{n = 1}^{k} X \to \ast_{n = 1}^{\infty} X$ and
$j_k \colon \ast_{n = 1}^{k} X \to \ast_{n = 1}^{k+1} X$
be the obvious inclusions. Denote by 
$\alpha_k \colon \ast_{n = 1}^{k} X \to \colim_{k \to \infty} \ast_{n = 1}^{k} X$
the structure map and by
$i \colon \colim_{k \to \infty} \ast_{n = 1}^{k} X \to  \ast_{n = 1}^{\infty} X$
the map induced by the system $\{i_k \mid k = 1,2, \ldots\}$.
This $G$-map is a (continuous) bijective $G$-map but not necessarily a $G$-homeomorphism. 
Since the partition $\{v_n\mid n = 1,2, \ldots\}$ is locally finite, we can find for
each $z \in Z$ an open $G$-invariant neighborhood $W_z$ of $z$ in $Z$ and a positive integer $k_z$
such that $v_n$ vanishes on $W_z$  for $n > k_z$. Define a map
$$f^{\prime}_z \colon W_z \to \ast_{n=1}^{k_z} X, \hspace{5mm} 
z \mapsto (v_1(z)\phi_1(z), v_2(z)\phi_2(z), \ldots , v_{k_z}(z)\phi_{k_z}(z)).$$
Then $\alpha_{k_z} \circ f^{\prime}_z \colon W_z \to \colim_{k \to \infty} \ast_{n = 1}^{k} X$ 
is a well-defined $G$-map whose composition with 
$i \colon \colim_{k \to \infty} \ast_{n = 1}^{k} X \to  \ast_{n = 1}^{\infty} X$ is
$f|_{W_z}$. Hence the system of the maps $\alpha_{k_z} \circ f^{\prime}_z$ defines a $G$-map
$$f^{\prime} \colon Z \to  \colim_{k \to \infty} \ast_{n = 1}^{k} X$$ 
such that $i \circ f^{\prime} = f$ holds. 

Let
$$\Delta_{n-1} ~ = ~ \{(t_1, t_2 \ldots t_n) \mid t_i \in [0,1], \sum_{i=1}^n t_i = 1\} ~ 
\subseteq ~ \prod_{n=1}^k [0,1]$$
be the standard $(n-1)$-simplex. Let
$$p \colon \left(\prod_{n=1}^k X\right) \times \Delta_n \to \ast_{n = 1}^{k} X, 
\hspace{5mm}
 (x_1, \ldots, x_n),(t_1, \ldots ,t_n) \mapsto (t_1x_1, \ldots ,t_nx_n)$$
be the obvious projection. It is a surjective continuous map but in general not an identification.
Let $\overline{\ast}_{n = 1}^{k} X$ be the topological space
whose underlying set is the same as for $\ast_{n = 1}^{k} X$ but 
whose topology is the quotient topology with respect to
$p$. The identity induces a (continuous) map $\overline{\ast}_{n = 1}^{k} X \to \ast_{n = 1}^{k} X$
which is not a homeomorphism in general.
Choose for $n \ge 1$ a (continuous) function $\phi_n \colon [0,1] \to [0,1]$ which 
satisfies $\phi_n^{-1}(0) = [0,4^{-n}]$.  Define
\begin{multline*}
u_k \colon \ast_{n = 1}^{k} X \to \overline{\ast}_{n = 1}^{k} X, \hspace{2mm}
\\
(t_nx_n \mid n = 1, \ldots, k) ~ \mapsto  ~
\left(\left.\frac{\phi_n(t_n)}{\sum_{n = 1}^k \phi_n(t_n)}x_n \right| n = 1, \ldots, k\right).
\end{multline*}
It is not hard to check that this $G$-map is continuous.
If $\overline{j}_k \colon \overline{\ast}_{n = 1}^{k} X \to \overline{\ast}_{n = 1}^{k+1} X$
is the obvious inclusion, we have $u_{k+1} \circ j_k = \overline{j}_k \circ u_k$ for all $k \ge 1$.
Hence the system of the maps $u_k$ induces a $G$-map
$$u \colon \colim_{k \to \infty} \ast_{n = 1}^{k} X \to \colim_{k \to \infty} \overline{\ast}_{n = 1}^{k} X.$$

Next we want to show that each $G$-space $\overline{\ast}_{n = 1}^{k} X$ has the $G$-homotopy type
of a $G$-$CW$-complex, whose isotropy groups belong to $\calf$. We first show that 
$\overline{\ast}_{n = 1}^{k} X$ is a $\left(\prod_{n = 1}^k G\right)$-$CW$-complex. It suffices to
treat the case $k = 2$, the general case follows by induction over $k$.
We can rewrite $X \overline{\ast} X$ as a $G \times G$-pushout
$$\comsquare{X \times X}{i_1}{CX \times X}{i_2}{}{X \times CX}{}{X \overline{\ast} X}$$
where $CX$ is the cone over $X$ and $i_1$ and $i_2$ are the obvious
inclusions. 
The product of two finite dimensional $G$-$CW$-complexes
is in a canonical way a finite dimensional $(G \times G)$-$CW$-complex, and, if $(B,A)$ is a $G$-$CW$-pair, $C$ a
$G$-$CW$-complex  and $f \colon B \to C$ is a cellular $G$-map, then
$A \cup_f C$ inherits a $G$-$CW$-structure in a canonical way. Thus $X \overline{\ast} X$ inherits
a $(G \times G)$-$CW$-complex structure.

The problem is now to decide whether the $\left(\prod_{n = 1}^k G\right)$-$CW$-complex
$\overline{\ast}_{n = 1}^{k} X$ regarded as a $G$-space by the diagonal action has the
$G$-homotopy type of a $G$-$CW$-complex.
If each $H \in \calf$ is open, then each isotropy group of the $G$-space
$\ast_{n = 1}^{k} X$ is open and we conclude from 
Remark \ref{rem: G-CW-complexes with open isotropy groups} 
that $\overline{\ast}_{n = 1}^{k} X$ with the diagonal $G$-action  is a $G$-$CW$-complex 
Suppose that $G$ is a Lie group and each $H \in \calf$ is compact. 
Example  \ref{exa: Lie groups acting properly and smoothly on manifolds} implies that for any
compact subgroup $K \subseteq \prod_{n=1}^k G$ the space
$\left( \prod_{n=1}^k G\right)/K$ regarded as $G$-space by the diagonal 
action has the $G$-homotopy type of a $G$-$CW$-complex. We conclude from
\cite[Lemma 7.4 on page 121]{Lueck(1989)}
that $\overline{\ast}_{n = 1}^{k} X$ with the diagonal $G$-action 
has the $G$-homotopy type of a $G$-$CW$-complex.
The isotropy groups $\overline{\ast}_{n = 1}^{k} X$ belong to $\calf$ 
since $\calf$ is closed under finite intersections and conjugation.
It is not hard to check that each
$G$-map $\overline{j}_k$ is a $G$-cofibration. 
Hence $\colim_{k \to \infty} \overline{\ast}_{n = 1}^{k} X$
has the $G$-homotopy type of a $G$-$CW$-complex, whose isotropy groups belong to $\calf$. 

Thus we have shown for every $\calf$-numerable $G$-space $Z$ that it admits a $G$-map to a 
$G$-$CW$-complex whose isotropy groups belong to $\calf$.
Now Lemma \ref{lem: E_FG = J_FG if all H in calf are open}
follows from Lemma \ref{lem: Comparison of E_FG and J_FG}. 
\end{proof}

\begin{definition}[Totally disconnected group]
 \label{def: totally disconnected group}
A (locally compact topological Hausdorff) group $G$ is called
\emph{totally disconnected}%
\index{group!totally disconnected}
if it satisfies one of the following equivalent conditions:
\begin{itemize}

\item[(T)] $G$ is totally disconnected as a topological space, i.e.
each component consists of one point;

\item[(D)]
The covering dimension of the topological  space $G$ is zero;

\item[(FS)]
Any element of $G$ has a fundamental system of compact open neighborhoods.

\end{itemize}
\end{definition}

We have to explain why these three conditions are equivalent. The implication
(T) $\Rightarrow$ (D) $\Rightarrow$ (FS) is shown in
\cite[Theorem 7.7 on page 62]{Hewitt-Ross(1979)}. It remains to prove
(FS) $\Rightarrow $ (T). Let $U$ be a subset of $G$ containing two distinct points 
$g$ and $h$. Let $V$ be a compact open neighborhood of $x$ which does not contain $y$. 
Then $U$ is the disjoint union of the open non-empty  sets $V \cap U$ and $V^c \cap U$
and hence disconnected. 

\begin{lemma} \label{lem: totally disconnected groups and calcom}
Let $G$ be a totally disconnected group and $\calf$ a family satisfying
$\calcomop \subseteq \calf \subseteq \calcom$.
Then the following square commutes up to $G$-homotopy
and consists of $G$-homotopy equivalences
$$\comsquare{\EGF{G}{\calcomop}}{u}{\JGF{G}{\calcomop}}
{}{}
{\EGF{G}{\calf}}{u}{\JGF{G}{\calf}}$$
where all maps come from the universal properties.
\end{lemma}
\begin{proof} We first show that any compact subgroup $H \subseteq G$ is contained in a
compact open subgroup. From 
\cite[Theorem 7.7 on page 62]{Hewitt-Ross(1979)} we get a compact open subgroup $K \subseteq G$.
Since $H$ is compact, we can find finitely many elements $h_1$, $h_2$, $\ldots$ , $h_s$ in $H$ 
such that $H \subseteq \bigcup_{i=1}^s h_iK$. Put $L := \bigcap_{h \in H} hKh^{-1}$.
Then $hLh^{-1} = L$ for all $h \in H$. Since $L = \bigcap_{i=1}^s h_iKh_i^{-1}$ holds,
$L$ is compact open. Hence $LH$ is a compact open subgroup containing $H$.

This implies that $\JGF{G}{\calf}$ is $\calcomop$-numerable. 
Obviously $\JGF{G}{\calcomop}$ is $\calf$-numerable. 
We conclude from the universal properties
that $\JGF{G}{\calcomop} \to \JGF{G}{\calf}$ is a $G$-homotopy equivalence.

The map $u \colon \EGF{G}{\calcomop} \to \JGF{G}{\calcomop}$ is a $G$-homotopy equivalence
because of  Lemma \ref{lem: E_FG = J_FG if all H in calf are open}.

This and Theorem \ref{the: G-homotopy characterization of JGF(G)(calf)} 
\ref{the: G-homotopy characterization of JGF(G)(calf): necessary condition}
imply that $\EGF{G}{\calcomop}^H$ is contractible for all $H \in \calf$.
Hence $\EGF{G}{\calcomop} \to \EGF{G}{\calf}$ is a $G$-homotopy equivalence
by Theorem \ref{the: G-homotopy characterization of EGF(G)(calf)}
\ref{the: G-homotopy characterization of EGF(G)(calf): characterization}. 
\end{proof}

\begin{definition}[Almost connected group] \label{def: almost connected} 
Given a group $G$, let $G^0$%
\indexnotation{G^0} be the normal subgroup given by the component of the identity%
\index{component of the identity}
and $\overline{G} = G/G^0$%
\indexnotation{overline(G)}
be the \emph{component group}.%
\index{component group}
We call $G$ \emph{almost connected}%
\index{group!almost connected}
if its component group $\overline{G}$ is compact.
\end{definition}
A Lie group $G$ is almost connected if and only if it has finitely many path components. 
In particular a discrete group is almost connected if it is finite.

\begin{theorem}[Comparison of $\EGF{G}{\calf}$ and $\JGF{G}{\calf}$]
\label{the: Examples where E_FG =J_FG}
\index{Theorem!Comparison of $\EGF{G}{\calf}$ and $\JGF{G}{\calf}$}
The map $u \colon  \EGF{G}{\calf}  \to  \JGF{G}{\calf}$ defined in 
\ref{map E_FG to J_FG}  is a $G$-homotopy equivalence if one of the following
conditions is satisfied:

\begin{enumerate}

\item \label{the: Examples where E_FG =J_FG: H in calf open}
Each element in $\calf$ is an open subgroup of $G$;

\item \label{the: Examples where E_FG =J_FG: G discrete}
The group $G$ is discrete;

\item \label{the: Examples where E_FG =J_FG: G Lie and H in calf compact}
The group $G$ is a Lie group and every element $H \in \calf$ is compact;

\item \label{the: Examples where E_FG =J_FG: G totally disconnected}
The group $G$ is totally disconnected and $\calf = \calcom$ or $\calf = \calcomop$;

\item \label{the: Examples where E_FG =J_FG: G almost connected}
The group $G$ is almost connected and each element in $\calf$ is compact.

\end{enumerate}

\end{theorem}
\begin{proof} Assertions  \ref{the: Examples where E_FG =J_FG: H in calf open}, 
\ref{the: Examples where E_FG =J_FG: G discrete},
\ref{the: Examples where E_FG =J_FG: G Lie and H in calf compact}
and \ref{the: Examples where E_FG =J_FG: G totally disconnected}  have already been proved in
Lemma \ref{lem: E_FG = J_FG if all H in calf are open}
and Lemma \ref{lem: totally disconnected groups and calcom}.
Assertion \ref{the: Examples where E_FG =J_FG: G almost connected} follows from
Lemma  \ref{lem: Comparison of E_FG and J_FG} 
and Theorem
\ref{the: almost connected groups}. 
\end{proof}

The following example shows that  the map 
$u \colon  \EGF{G}{\calf}  \to  \JGF{G}{\calf}$ defined in 
\ref{map E_FG to J_FG}  is in general not a $G$-homotopy equivalence.
\begin{example}[Totally disconnected groups and $\caltr$]
 \label{exa: E_TRG not equl to J_TRG}
\em
Let $G$ be totally disconnected. We claim that
$u \colon  \EGF{G}{\caltr}  \to  \JGF{G}{\caltr}$ defined in 
\ref{map E_FG to J_FG}  is a $G$-homotopy equivalence if and only if $G$ is discrete.
In view of Theorem \ref{the: G-homotopy characterization of JGF(G)(calf)}
\ref{the: G-homotopy characterization of JGF(G)(calf): necessary condition}
and Lemma \ref{lem: E_FG = J_FG if all H in calf are open} this is equivalent to the statement
that $\EGF{G}{\caltr}$ is contractible
if and only if $G$ is discrete. If $G$ is discrete, we already know that
$\EGF{G}{\caltr}$ is contractible. Suppose now that
$\EGF{G}{\caltr}$ is contractible. We obtain a numerable
$G$-principal bundle $G \to \EGF{G}{\caltr} \to G\backslash \EGF{G}{\caltr}$
by Remark \ref{rem: G-principal bundles over CW-complexes}.
This implies that it is a fibration by a result of Hurewicz
\cite[Theorem on p. 33]{Whitehead(1978)}.
Since $\EGF{G}{\caltr}$ is contractible, $G$ and the loop space
$\Omega(G\backslash \EGF{G}{\caltr})$ are homotopy equivalent
\cite[6.9$^*$ on p. 137, 6.10$^*$ on p. 138,
Corollary 7.27 on p. 40]{Whitehead(1978)}.
Since $G\backslash \EGF{G}{\caltr}$ is a $CW$-complex,
$\Omega(G\backslash \EGF{G}{\caltr})$ has the homotopy type of a $CW$-complex
\cite{Milnor(1959)}. Hence there exists a homotopy equivalence
$f\colon G \to X$ be from $G$ to a $CW$-complex $X$.
Then the induced map $\pi_0(G) \to \pi_0(X)$ between the set of path components
is bijective.
Hence the  preimage of each path component of $X$ is a path component of $G$ and therefore 
a point since $G$ is totally disconnected.
Since $X$ is locally path-connected
each path component of $X$ is open in $X$.
We conclude that $G$ is the disjoint union
of the preimages of the path components of $X$ and each of these preimages
is open in $G$ and consists of one point.  Hence $G$ is discrete.
\end{example}

 \begin{remark}[Compactly generated spaces] \label{rem: compactly generated spaces} \em
 In the following theorem we will work in the category of compactly
 generated spaces. This convenient category is explained in detail in
 \cite{Steenrod(1967)} and \cite[I.4]{Whitehead(1978)}. A reader may
 ignore this technical point in the following theorem without harm, but
 we nevertheless give a short explanation. 

 A Hausdorff
 space $X$ is called \emph{compactly generated} if a subset $A
 \subseteq X$ is closed if and only if $A \cap K$ is closed for every
 compact subset $K \subseteq X$. Given a topological space $X$, let
 $k(X)$ be the compactly generated topological space with the same
 underlying set as $X$ and the topology for which a subset $A \subseteq
 X$ is closed if and only if for every compact subset $K \subseteq X$ the
 intersection $A \cap K$ is closed in the given topology on $X$. The
 identity induces a continuous map $k(X) \to X$ which is a
 homeomorphism if and only if $X$ is compactly generated. 
 The spaces $X$ and $k(X)$ have the same compact
 subsets. Locally compact Hausdorff spaces and every Hausdorff space
 which satisfies the first axiom of countability are compactly
 generated. In particular metrizable spaces are compactly generated.

  Working in the category of compactly generated
 spaces means that one only considers compactly generated spaces and
 whenever a topological construction such as the cartesian product 
 or the mapping space leads out of this category, one retopologizes the
 result as described above to get a
 compactly generated space. The advantage is for example that in the
 category of compactly generated spaces the exponential
 map $\map(X \times Y,Z) \to \map(X,\map(Y,Z))$ is always a
 homeomorphism, for an identification $p \colon X \to Y$ the map $p
 \times \id_Z \colon X \times Z \to Y \times Z$ is always an
 identification and for a filtration by closed subspaces
 $X_1 \subset X_2 \subseteq \ldots \subseteq X$ such that $X$ is the
 colimit $\colim_{n \to \infty} X_n$, we always get 
 $X \times Y = \colim_{n \to \infty} (X_n \times Y)$. In particular the
 product of a $G$-$CW$-complex $X$ with a $H$-$CW$-complex $Y$ is in a
 canonical way a $G \times H$-$CW$-complex $X \times Y$. Since
 we are assuming that $G$ is a locally compact Hausdorff group, any
 $G$-$CW$-complex $X$ is compactly generated. 
 \em
 \end{remark}

The following result has grown out of discussions with Ralf Meyer.

 \begin{theorem}[Equality of $\underline{E}G$ and $\underline{J}G$]
 \label{the : Equality of underline{E}G and underline{J}}
 Let $G$ be a locally compact topological Hausdorff group. Then  the
 canonical $G$-map $\underline{E}G \to \underline{J}G$ is a
 $G$-homotopy equivalence.
 \end{theorem}
 \begin{proof}
 In the sequel of the proof we work in the category of compactly generated spaces
 (see Remark~\ref{rem: compactly generated spaces}). Notice that the
 model mentioned in Theorem~\ref{the: G-homotopy characterization of JGF(G)(calf)}
\ref{the: G-homotopy characterization of JGF(G)(calf): existence}  is
metrizable and hence compactly generated (see \cite[Appendix 1]{Baum-Connes-Higson(1994)}).
Because of Lemma~\ref{lem: Comparison of E_FG and J_FG} it suffices
to construct a $G$-$CW$-complex $Z$ with compact isotropy groups
together with a $G$-map $\underline{J}G \to Z$.

Let $G^0$ be the component of the identity which is a normal closed subgroup.
Let $p \colon G \to G/G^0$ be the projection. The groups $G^0$ and $G/G^0$ are locally
compact Hausdorff groups and $G/G^0$ is totally disconnected. We conclude from
Lemma~\ref{lem: totally disconnected groups and calcom} that there is a $G$-map
$\underline{J}(G/G^0) \to \EGF{G/G^0}{\calcomop}$. Since $\underline{J}G$ is
$\calcomop$-numerable, the $G/G^0$-space $G^0\backslash(\underline{J}G)$ is
$\calcom$-numerable and hence there exists a $G/G^0$-map 
$G^0\backslash(\underline{J}G) \to \underline{J}(G/G^0) $. Thus we get a $G$-map
$u \colon \underline{J}G \to \res_p  \EGF{G/G^0}{\calcomop}$, where the $G$-$CW$-complex 
$\res_p  \EGF{G/G^0}{\calcomop}$ is obtained from the  $G/G^0$-$CW$-complex
$\EGF{G/G^0}{\calcomop}$ by letting $g \in G$ act by $p(g)$.
We obtain a $G$-map
$\id  \times f \colon \underline{J}G \to \underline{J}G \times \res_p
\EGF{G/G^0}{\calcomop}$. Hence it suffices to
construct a $G$-$CW$-complex $Z$ with compact isotropy groups together with a $G$-map
$f \colon \underline{J}G \times \res_p  \EGF{G/G^0}{\calcomop} \to Z$.
For this purpose we construct a sequence of $G$-$CW$-complexes
$Z_{-1} \subseteq Z_0 \subseteq Z_1 \subseteq \ldots$ such that 
$Z_n$ is a $G$-$CW$-subcomplex of $Z_{n+1}$ and each $Z_n$ has compact
isotropy groups, and $G$-homotopy equivalences 
$f_n \colon \res_p  \EGF{G/G^0}{\calcomop}_n \times \underline{J}G\to Z_n$. 
with $f_{n+1}|_{\res_p  \EGF{G/G^0}{\calcomop}_n} = f_n$, where $\EGF{G/G^0}{\calcomop}_n$ is the $n$-skeleton of
$\EGF{G/G^0}{\calcomop}$. The canonical $G$-map 
$$ \colim_{n \to \infty} \left(\underline{J}G \times \res_p
\EGF{G/G^0}{\calcomop}_n\right)  
~ \to ~
\underline{J}G \times \res_p\EGF{G/G^0}{\calcomop}$$
is a $G$-homeomorphism. The $G$-space $Z = \colim_{n \to \infty} Z_n$
is a $G$-$CW$-complex with compact isotropy groups. 
Hence we can define the desired $G$-map by $f = \colim_{n \to \infty} f_n$ after we have
constructed the $G$-maps $f_n$. This will be done by induction over $n$. The induction
beginning $n = -1$ is given by $\id \colon \emptyset \to \emptyset$. The induction step
from $n$ to $(n+1)$ is done as follows. Choose a $G/G^0$-pushout
$$\comsquare{\coprod_{i \in I} (G/G^0)/H_i \times S^n}{}{\EGF{G/G^0}{\calcomop}_n}
{}{}
{\coprod_{i \in I} (G/G^0)/H_i \times D^{n+1}}{}{\EGF{G/G^0}{\calcomop}_{n+1}}
$$
where each $H_i$ is a compact open subgroup of $G/G^0$.  We obtain a $G$-pushout
$$\comsquare{\coprod_{i \in I} \res_p \left((G/G^0)/H_i \times S^n\right) \times
  \underline{J}G}
{}
{\res_p\EGF{G/G^0}{\calcomop}_n \times \underline{J}G}
{}{}
{\coprod_{i \in I} \res_p\left((G/G^0)/H_i \times D^{n+1}\right)\times
  \underline{J}G}{}
{\res_p\EGF{G/G^0}{\calcomop}_{n+1}  \times \underline{J}G}
$$
In the sequel let $K_i \subseteq G$ be the open almost connected subgroup $p^{-1}(H_i)$. 
The $G$-spaces $\res_p (G/G^0)/H_i$ and $G/K_i$ agree. We have the $G$-homeomorphism
$$G \times_{K_i} \res_G^{K_i} \underline{J}G  \xrightarrow{\cong}  
G/K_i \times \underline{J}G, \quad (g,x) \mapsto (gK_i,gx).$$
Thus we obtain a $G$-pushout
\begin{eqnarray}
& \comsquare{\left(\coprod_{i \in I} G \times_{K_i} (\res_G^{K_i} \underline{J}G)\right)\times S^n}
{w}
{\res_p\EGF{G/G^0}{\calcomop}_n \times \underline{J}G}
{\id \times i}{}
{\left(\coprod_{i \in I} G \times_{K_i} (\res_G^{K_i}
    \underline{J}G)\right) \times D^{n+1}}
{}
{\res_p\EGF{G/G^0}{\calcomop}_{n+1}  \times \underline{J}G} 
&
\label{pushout for X_{n+1}}
\end{eqnarray}
where $i \colon S^n \to D^{n+1}$ is the obvious inclusion.

Let $X$ be a $\calcom$-numerable $K_i$-space. Then 
the $G$-space $G \times_{K_i} \underline{J}K_i$ is a
$\calcom$-numerable and hence admits a $G$-map to $\underline{J}G$.
Its restriction to $\underline{J}K_i = K_i \times_{K_i} \underline{J}K_i$
defines a $K_i$-map $f \colon X \to \res_G^{K_i}\underline{J}G$. 
If $f_1$ and $f_2$ are $K_i$-maps $X \to \res_G^{K_i}\underline{J}G$,
we obtain $G$-maps $\overline{f_k} \colon G \times_{K_i} X \to
\underline{J}G$  by sending $(g,x) \to gf_k(x)$ for $k =
0,1$. By the universal property of $\underline{J}G$ these two $G$-maps
are $G$-homotopic. Hence $f_0$ and $f_1$ are $K_i$-homotopic. 
Since $K_i \subseteq G$ is open, $\res_G^{K_i}\underline{J}G$ is a $\calcom$-numerable
$K_i$-space. Hence the $K_i$-space $\res_G^{K_i}\underline{J}G$ is a model for $\underline{J}K_i$.
Since $K_i$ is almost connected, there is a $K_i$-homotopy equivalence
$ \underline{E}K_i \to \res_G^{K_i}\underline{J}G$ by 
Theorem~\ref{the: Examples where E_FG =J_FG}
\ref{the: Examples where E_FG =J_FG: G almost connected}. Hence we obtain
a $G$-homotopy equivalence
$$u_i \colon G\times_{K_i} \underline{E}K_i 
~ \to ~
G \times_{K_i} (\res_G^{K_i} \underline{J}G)
$$
with a $K_i$-$CW$-complex with compact isotropy groups as source. 

In the sequel we abbreviate
\begin{eqnarray*}
X_n & := & \res_p\EGF{G/G^0}{\calcomop}_n \times \underline{J}G;
\\
Y & := & \coprod_{i \in I} G \times_{K_i} (\res_G^{K_i}
\underline{J}G);
\\
Y' & = & \coprod_{i \in I} G\times_{K_i} \underline{E}K_i.
\end{eqnarray*}
Choose a $G$-homotopy equivalence $v \colon Y' \to Y$.
By the equivariant cellular Approximatiom Theorem we can
find a $G$-homotopy $h \colon Y' \times S^n \times [0,1]\to Z_n$ 
such that $h_0 = f_n \circ w \circ  (v \times \id_{S^n})$ and the
$G$-map $h_1 \colon Y \times S^n \to Z_n$ is
cellular. Consider the following commutative diagram of $G$-spaces
$$
\begin{CD}
Y \times D^{n+1} @< \id_Y \times i << Y \times S^n @> w >> X_n
\\
@V \id VV @V \id VV  @V f_n VV
\\
Y \times D^{n+1} @< \id_Y \times i << Y \times S^n @> f_n \circ w >> Z_n
\\
@A v \times \id_{D^{n+1}} AA @A v \times \id_{S^n} AA  @A \id AA
\\
Y' \times D^{n+1} @< \id_{Y'} \times i << Y' \times S^n @> f_n \circ w
\circ (v \times \id_{S^n})>> Z_n
\\
@V j_0 VV @V j_0 VV @ V \id VV
\\
Y' \times D^{n+1} \times [0,1] @< \id_{Y'} \times i \times \id_{[0,1]} <<
Y' \times S^n \times [0,1]@>  h >> Z_n
\\
@A j_1 AA @A j_1AA @A \id AA
\\
Y' \times D^{n+1} @< \id_{Y'} \times i << Y' \times S^n @> h_1 >> Z_n
\end{CD}
$$
where $j_0$ and $j_1$ always denotes the obvious inclusions.
The $G$-pushout of the top row is $X_{n+1}$ by
\eqref{pushout for X_{n+1}}. Let $Z_{n+1}$ be the $G$-pushout of the
bottom row. This is a $G$-$CW$-complex with compact isotropy groups
containing $Z_n$ as $G$-$CW$-subcomplex.  Let $W_2$ and $W_3$ and $W_4$ be the $G$-pushout of the
second, third and fourth row. The diagram above induces a sequence of
$G$-maps
$$X_{n+1} \xrightarrow{u_1} W_2 \xleftarrow{u_2} W_3 \xrightarrow{u_3}
W_4 \xleftarrow{u_4} Z_{n+1}$$
The left horizontal arrow in each row is a $G$-cofibration as $i$ is a
cofibration. Each of the vertical arrows is a $G$-homotopy
equivalence. This implies that each of the maps $u_1$, $u_2$, $u_3$
and $u_4$ are $G$-homotopy equivalences. Notice that we can consider
$Z_n$ as a subspace of $W_2$, $W_3$, $W_4$ such that the inclusion
$Z_n \to W_k$ is a $G$-cofibration. Each of the maps $u_2$, $u_3$ and
$u_4$ induces the identity on $Z_n$, whereas $u_1$ induces $f_n$ on
$X_n$. By a cofibration argument one can find $G$-homotopy inverses
$u_2^{-1}$ and $u_4^{-1}$ of $u_2$ and $u_4$ which induce the identity
on $Z_n$. Now define the desired $G$-homotopy equivalence $f_{n+1} 
\colon X_{n+1} \to Z_{n+1}$ to be the composition $u_4^{-1} \circ u_3
\circ u_2^{-1} \circ u_1$.  This finishes the proof of
Theorem~\ref{the : Equality of underline{E}G and underline{J}}.
\end{proof}


\typeout{--------------------   Section 4: Special Models --------------------------}

\section{Special Models}
\label{sec: Special Models}

In this section we present some interesting geometric models for the
space $\EGF{G}{\calf}$ and $\JGF{G}{\calf}$ focussing on
$\underline{E}G$ and $\underline{J}G$. In particular we are interested in cases,
where these models satisfy finiteness conditions such as being
finite, finite dimensional or of finite type.

One extreme case is, where we take $\calf$ to be the family $\calall$ of
all subgroups. Then a model for both $\EGF{G}{\calall}$ and $\JGF{G}{\calall}$
is $G/G$. The other extreme case is the family $\caltr$ consisting of the trivial subgroup.
This case has already been treated in Example \ref{exa: Numerable G-principal bundles},
Remark \ref{rem: G-principal bundles over CW-complexes} and Example \ref{exa: E_TRG not equl to J_TRG}.


\subsection{Operator Theoretic Model}
\label{subsec: Operator Theoretic Model}

Let $G$ be a locally compact Hausdorff topological group. 
Let $C_0(G)$ be the Banach space of complex valued functions of $G$ vanishing at infinity
with the supremum-norm. The group $G$ acts isometrically on $C_0(G)$ by 
$(g\cdot f)(x) := f(g^{-1}x)$ for $f \in C_0(G)$ and $g,x \in G$.
Let $PC_0(G)$%
\indexnotation{PC_0(G)} 
be the subspace of $C_0(G)$ consisting of functions $f$ such that
$f$ is not identically zero and has non-negative real numbers as values. 

The next theorem is due to Abels \cite[Theorem 2.4]{Abels(1978)}.

\begin{theorem}[Operator theoretic model]
\label{the: Operator theoretic model}
\index{Theorem!Operator theoretic model}
The $G$-space $PC_0(G)$ is a model for $\underline{J}G$.
\end{theorem}

\begin{remark}\em
Let $G$ be discrete.  
Another model for $\underline{J}G$ is the space
$$X_G = \{f \colon G \to [0,1] \mid f \text{ has finite support, }
\sum_{g \in G} f(g) = 1\}$$
with the topology coming from the supremum norm
\cite[page 248]{Baum-Connes-Higson(1994)}.
Let $P_{\infty}(G)$ be the geometric realization of the simplicial
set whose $k$-simplices consist of $(k+1)$-tupels
$(g_0,g_1, \ldots , g_k)$ of elements $g_i$ in $G$. This also 
a model for $\underline{E}G$ \cite[Example 2.6]{Abels(1978)}. 
The spaces $X_G$ and $P_{\infty}(G)$ have the same underlying sets
but in general they have different topologies. The identity map
induces a (continuous) $G$-map $P_{\infty}(G) \to X_G$ which is a $G$-homotopy equivalence, but
in general not a $G$-homeomorphism (see also \cite[A.2]{Valette(2002)}). \em
\end{remark}


\subsection{Almost Connected  Groups}
\label{subsec: Almost Connected Groups}

The next result is due to Abels \cite[Corollary 4.14]{Abels(1978)}.
\begin{theorem}[Almost connected groups] \label{the: almost connected groups}
\index{Theorem!Almost connected groups}
Let $G$ be a (locally compact Hausdorff) topological group.
Suppose that $G$ is almost connected, i.e.\ the group $G/G^0$ is compact for
$G^0$ the component of the identity element. Then $G$ contains a
maximal compact subgroup $K$ which is unique up to conjugation.
The $G$-space $G/K$ is a model for $\underline{J}G$.
\end{theorem}

The next result follows from
Example \ref{exa: Lie groups acting properly and smoothly on manifolds}, 
Theorem \ref{the: Examples where E_FG =J_FG}
\ref{the: Examples where E_FG =J_FG: G Lie and H in calf compact} and 
Theorem \ref{the: almost connected groups}.

\begin{theorem}[Discrete subgroups of almost connected Lie groups] 
\label{the: Discrete subgroups of almost connected Lie groups}
\index{Theorem!Discrete subgroups of almost connected Lie groups}
Let $L$ be a Lie group with finitely many path components.
Then $L$ contains a maximal compact subgroup $K$ which is unique up to conjugation.
The $L$-space $L/K$ is a model for $\underline{E}L$.

If $G \subseteq L$ is a discrete subgroup of $L$, then 
$L/K$ with the obvious left $G$-action is a finite dimensional
$G$-$CW$-model for $\underline{E}G$.
\end{theorem}


\subsection{Actions on Simply Connected Non-Positively Curved Manifolds}
\label{subsec: Actions on Simply Connected Non-Positively Curved Manifolds}

The next theorem is due to Abels \cite[Theorem 4.15]{Abels(1978)}.
 
\begin{theorem}[Actions on simply connected non-positively curved manifolds]
\label{the: actions on simply connected non-posively curved manifolds}
\index{Theorem!Actions on simply connected non-positively curved manifolds}
Let $G$ be a (locally compact Hausdorff) topological group. Suppose
that $G$ acts properly and isometrically on the
simply-connected complete Riemannian manifold $M$ with non-positive sectional curvature.
Then $M$ is a model for $\underline{J}G$.
\end{theorem}


\subsection{Actions on CAT(0)-spaces}
\label{subsec: Actions on CAT(0)-spaces}

\begin{theorem}[Actions on CAT(0)-spaces]
\label{the: Actions on CAT(0)-spaces}
\index{Theorem!Actions on CAT(0)-spaces}
Let $G$ be a (locally compact Hausdorff) topological group. Let $X$ be a proper $G$-$CW$-complex.
Suppose that $X$ has the structure of a complete {\rm CAT(0)}-space for which $G$ acts by isometries.
Then $X$ is a model for $\underline{E}G$.
\end{theorem}
\begin{proof}
By \cite[Corollary II.2.8 on page 179]{Bridson-Haefliger(1999)} the $K$-fixed point set
of $X$ is a non-empty convex subset of $X$ and hence contractible for any compact subgroup $K \subset G$. 
\end{proof}

This result contains as special case 
Theorem~\ref{the: actions on simply connected non-posively curved  manifolds} 
and partially Theorem~\ref{the: Actions on trees} since
simply-connected complete Riemannian manifolds with non-positive sectional curvature
and trees are CAT(0)-spaces.


\subsection{Actions on Trees and Graphs of Groups}
\label{subsec: Actions on Trees and Graphs of Groups}

A \emph{tree}%
\index{tree}
is a $1$-dimensional $CW$-complex which is contractible.

\begin{theorem}[Actions on trees]
\label{the: Actions on trees}
\index{Theorem!Actions on trees}
Suppose that $G$ acts continuously on a tree $T$
such that for each element
$g \in G$ and each open cell $e$ with $g \cdot e \cap e \not= \emptyset$ we have
$gx = x$ for any $x \in e$. Assume that the isotropy group
of each $x \in T$ is compact. 

Then $G$ can be written as an extension
$1 \to K \to G \to \overline{G} \to 1$ of a compact group containing
$G^0$ and a totally disconnected group $\overline{G}$ such that $K$
acts trivially and $T$ is a $1$-dimensional model for 
$$\EGF{G}{\calcom} = \JGF{G}{\calcom} =\EGF{G}{\calcomop} = \JGF{G}{\calcomop}.$$
\end{theorem}
\begin{proof} 
We conclude from Remark \ref{rem: G-CW-complexes with open isotropy groups} 
that $T$ is a $G$-$CW$-complex and all isotropy groups are compact
open. Let $K$ be the intersection of all the isotropy groups of points
of $T$. This is a normal compact  subgroup of $G$ which contains the component of the
identity $G^0$. Put $\overline{G} = G/K$. This is a totally
disconnected group. Let $H \subseteq G$ be compact. If $e_0$ is a zero-cell in
$T$, then $H \cdot e_0$ is a compact discrete set and hence finite.
Let $T'$ be the union of all geodesics with extremities in $H \cdot
e$. This is a $H$-invariant subtree of $T$ of finite diameter. One shows now
inductively over the diameter of $T'$ that $T'$ has a vertex which is
fixed under the $H$-action (see \cite[page 20]{Serre(1980)} or 
\cite[Proposition 4.7 on page 17]{Dicks-Dunwoody(1989)}). Hence $T^H$ is non-empty. 
If $e$ and $f$ are vertices in $T^H$, the geodesic in $T$ from $e$ to
$f$ must be $H$-invariant. Hence $T^H$ is a connected $CW$-subcomplex
of the tree $T$ and hence is itself a tree. This shows that $T^H$ is
contractible. Hence $T$ is a model for $\EGF{G}{\calcom} = \EGF{\overline{G}}{\calcom}$. Now apply
Lemma~\ref{lem: totally disconnected groups and calcom}.
\end{proof}

Let $G$ be a locally compact Hausdorff group.
Suppose that $G$ acts continuously on a tree $T$
such that for each element
$g \in G$ and each open cell $e$ with $g \cdot e \cap e \not= \emptyset$ we have
$gx = x$ for any $x \in e$. If the $G$-action on a tree has possibly not
compact isotropy groups, one can nevertheless get nice models for $\EGF{G}{\calcomop}$
as follows. Let $V$ be the set of equivariant $0$-cells and $E$ be the set of equivariant $1$-cells of $T$.
Then we can choose a $G$-pushout
\begin{eqnarray}
&\comsquare{\coprod_{e \in E} G/H_e \times \{-1,1\}}{q}{T_0 = \coprod_{v \in V}G/K_v}
{}{}
{\coprod_{e \in E} G/H_e \times [-1,1]}{}{T}
&
\label{G-pushout for the tree T}
\end{eqnarray}
where the left vertical arrow is the obvious inclusion.
Fix $e \in E$ and $\sigma \in \{-1,1\}$. Choose elements
$v(e,\sigma) \in V$ and $g(e,\sigma) \in G$ such that
$q$ restricted to $G/H_e \times \{\sigma\}$ is the $G$-map $G/H_e \to G/K_{v(e,\sigma)}$
which sends $1H_e$ to $g(e,\sigma)K_{v(e,\sigma)}$. Then conjugation with
$g(e,\sigma)$ induces a group homomorphism $c_{g(e,\sigma)} \colon H_e \to K_{v(e,\pm 1)}$
and there is an up to equivariant homotopy unique $c_{g(e,\sigma)}$-equivariant cellular map
$f_{g(e,\sigma)} \colon \EGF{H_e}{\calcomop} \to \EGF{K_{e(g,\sigma)}}{\calcomop}$.
Define a $G$-map
$$Q \colon \coprod_{e \in E} G \times_{H_e} \EGF{H_e}{\calcomop} \times \{-1,1\}
~ \to ~ 
\coprod_{v \in V} G \times_{K_v} \EGF{K_v}{\calcomop}$$
by requiring that the restriction of
$Q$ to $G \times_{H_e} \EGF{H_e}{\calcomop} \times \{\sigma \}$ is the $G$-map
$$G \times_{H_e} \EGF{H_e}{\calcomop} ~ \to G \times_{K_{v(e,\sigma)}} \EGF{K_{g(e,\sigma)}}{\calcomop},
\hspace{5mm} (g,x) \mapsto (g,f_{g(e,\sigma)}(x)).$$
Let $T_{\calcomop}$ be the $G$-pushout
$$
\comsquare{\coprod_{e \in E} G \times_{H_e} \EGF{H_e}{\calcomop} \times \{-1,1\}}{Q}
{\coprod_{v \in V} G \times_{K_v} \EGF{K_v}{\calcomop}}
{}{}
{\coprod_{e \in E} G \times_{H_e} \EGF{H_e}{\calcomop} \times [-1,1]}{}{T_{\calcomop}}
$$
The $G$-space $T_{\calcomop}$ inherits a canonical $G$-$CW$-structure
with compact open isotropy groups.
Notice that for any open subgroup $L \subseteq G$ one can choose as model for $\EGF{L}{\calcomop}$ the 
restriction $\res_G^L \EGF{G}{\calcomop}$ of $\EGF{G}{\calcomop}$ to $L$
and that there is a $G$-homeomorphism 
$G \times_L \res_G^L \EGF{G}{\calcomop}\xrightarrow{\cong} G/L \times \EGF{G}{\calcomop}$ 
which sends $(g,x)$ to $(gL,gx)$. This implies that $T_{\calcomop}$ is
$G$-homotopy equivalent to $T \times \EGF{G}{\calcomop}$ with the diagonal $G$-action.
If $H \subseteq G$ is compact open, then $T^H$ is contractible. Hence
$(T \times \EGF{G}{\calcomop})^H$ is contractible for compact open subgroup $H \subseteq G$.
Theorem \ref{the: G-homotopy characterization of EGF(G)(calf)}
\ref{the: G-homotopy characterization of EGF(G)(calf): characterization}
shows
\begin{theorem}[Models based on actions on trees]
 \label{the: Models based on actions on trees}
\index{Theorem!Models based on actions on trees}
The $G$-$CW$-complex $T_{\calcomop}$ is a model for $\EGF{G}{\calcomop}$.
\end{theorem}
The point is that it may be possible to choose nice models for
the various spaces $\EGF{H_e}{\calcomop}$ and $\EGF{K_v}{\calcomop}$ and thus 
get a nice model for $\EGF{G}{\calcomop}$. 
If all isotropy groups of the $G$-action on $T$ are compact, we can
choose all spaces $\EGF{H_e}{\calcomop}$ and $\EGF{K_v}{\calcomop}$ to be
$\pt$ and we rediscover Theorem \ref{the: Actions on trees}.

Next we recall which discrete groups $G$ act on trees. 
Recall that an oriented graph $X$ is a $1$-dimensional $CW$-complex together with
an orientation for each $1$-cell. This can be codified by
specifying a triple $(V,E,s \colon E \times \{-1,1\} \to V)$ consisting of two sets
$V$ and $E$ and a map $s$. The associated oriented graph is the pushout
$$\comsquare{E \times \{-1,1\}}{s}{V}{}{}{E \times[0,1]}{}{X}$$
So $V$ is the set of vertices, $E$ the set of edges, and for a edge $e \in E$ 
its initial vertex is $s(e,-1)$ and its terminal vertex is $s(e,1)$.
A \emph{graph of groups}%
\index{graph of groups}
$\calg$ on a connected oriented graph $X$ consists of two sets of groups
$\{K_v \mid v \in V\}$ and $\{H_e\mid e \in E\}$ with $V$ and $E$ as index sets
together with injective group homomorphisms $\phi_{v,\sigma} \colon H_e \to
K_{s(e,\sigma)}$ for each $e \in E$. Let $X_0 \subseteq X$ be some maximal tree.
We can associate to these data \emph{the fundamental group}%
\index{graph of groups!fundamental group}
\index{fundamental group of a graph of groups}
$\pi = \pi(\calg,X,X_0)$%
\indexnotation{pi(calg,X,X_0)}
 as follows. Generators of $\pi$ are the elements in $K_v$ for each $v \in V$ and 
the set $\{t_e \mid e \in E\}$. The relations are the relations in each group $K_v$ for each
$v \in V$, the relation $t_e = 1$ for $e \in V$ if $e$  belongs to $X_0$,
and for each $e \in E$ and $h \in H_e$ we require
$t_e^{-1}\phi_{e,-1}(h)t_e = \phi_{e,+1}(h)$. It turns out that the obvious map
$K_v \to \pi$ is an injective group homomorphism for each $v \in V$ and we will identify
in the sequel $K_v$ with its image in $\pi$ \cite[Corollary 7.5 on page 33]{Dicks-Dunwoody(1989)},
\cite[Corollary 1 in 5.2 on page 45]{Serre(1980)}.
We can assign to these
data a tree $T = T(X,X_0,\calg)$ with $\pi$-action as follows. Define a 
$\pi$-map
$$q \colon \coprod_{e \in E} \pi/\im(\phi_{e,-1}) \times \{-1,1\}
~ \to ~ \coprod_{v \in V} \pi/K_v$$
by requiring that its restriction to $\pi/\im(\phi_{e,-1}) \times \{-1\}$
is the $\pi$-map given by the projection 
$\pi/\im(\phi_{e,-1}) \to \pi/K_{s(e,-1)}$ and
its restriction to $\pi/\im(\phi_{e,-1}) \times \{1\}$
is the $\pi$-map  $\pi/\im(\phi_{e,-1}) \to \pi/K_{s(e,1)}$
which sends $g\im(\phi_{e,-1})$ to $gt_e\im(\phi_{e,1})$.
Now define a $1$-dimensional $G$-$CW$-complex $T = T(\calg,X,X_0)$%
\indexnotation{T(calg,X,X_0)}
using this $\pi$-map $q$ and the $\pi$-pushout analogous to 
\eqref{G-pushout for the tree T}. It turns out that $T$ is
contractible \cite[Theorem 7.6  on page 33]{Dicks-Dunwoody(1989)},
\cite[Theorem 12 in 5.3 on page 52]{Serre(1980)}.

On the other hand, suppose that $T$ is a $1$-dimensional $G$-$CW$-complex.
Choose a $G$-pushout \eqref{G-pushout for the tree T}. Let $X$ be the
connected oriented graph $G\backslash T$. It has a set of vertices $V$ 
and as set of edges the set $E$.
The required map $s \colon E \times \{-1,1\} \to V$ sends $s(e,\sigma)$ to the vertex
for which $q(G/H_e \times \{\sigma\})$ meets and hence is equal to $G/K_{s(e,\sigma)}$.
Moreover, we get a graph of groups $\calg$ on $X$ as follows. Let $\{K_v \mid v \in V\}$
and $\{H_e \mid e \in E\}$ be the set of groups given by
\eqref{G-pushout for the tree T}.  Choose an element $g(e,\sigma) \in G$ such that
the $G$-map induced by $q$ from $G/H_e$ to $G/K_{s(e,\sigma)}$ sends $1H_e$ to
$g(e,\sigma)K_{s(e,\sigma)}$. Then conjugation with $g(e,\sigma)$ induces
a group homomorphism $\phi_{e,\sigma} \colon H_e \to K_{s(e,\sigma)}$.
After a choice of a maximal tree $X_0$ in $X$ one obtains an isomorphism
$G \cong \pi(\calg,X,X_0)$. (Up to isomorphism) we get a bijective
correspondence between pairs $(G,T)$ consisting of a group $G$ acting on
an oriented tree $T$ and a graph of groups on  connected oriented
graphs. For details we refer for instance to 
\cite[I.4 and I.7]{Dicks-Dunwoody(1989)}  and 
\cite[\S 5]{Serre(1980)}.

\begin{example}[The graph associated to amalgamated products]
 \label{exa: amalgamated products} \em
Consider the graph $D$ with one edge $e$ and two vertices $v_{-1}$ and $v_1$
and the map $s \colon \{e\} \times \{-1,1\} \to \{v_{-1},v_1\}$ which sends
$(e,\sigma)$ to $v_{\sigma}$. Of course this is just the graph consisting of a single
segment which is homeomorphic to $[-1,1]$. 
Let $\calg$ be a graph of groups on $D$. This is the same as
specifying a group $H_e$ and groups $K_{-1}$ and $K_1$ together with
injective group homomorphisms $\phi_{\sigma} \colon H_e \to K_{\sigma}$ for
$\sigma \in \{-1,1\}$. There is only one choice of a maximal subtree in $D$, namely
$D$ itself. Then the fundamental group $\pi$ of this graph of groups
is the amalgamated product of $K_{-1}$ and $K_1$ over $H_e$ with respect to
$\phi_{-1}$ and $\phi_1$, i.e.\ the pushout of groups
$$\comsquare{H_e}{\phi_{-1}}{K_{-1}}{\phi_1}{}{K_1}{}{\pi}$$
Choose $\phi_{\sigma}$-equivariant maps
$f_{\sigma} \colon \underline{E}H_e \to \underline{E}K_{\sigma}$.
They induce $\pi$-maps 
$$F_{\sigma} \colon \pi \times_{H_e} \underline{E}H_e ~ \to ~ 
 \pi \times_{K_{\sigma}} \underline{E}K_{\sigma} , \hspace{5mm} 
(g,x) \mapsto (g,f_{\sigma}(x)).$$
We get a model for $\underline{E}\pi$ as
the $\pi$-pushout
$$\comsquare{\pi \times_{H_e} \underline{E}H_e \times \{-1,1\}}{F_{-1} \coprod F_1}
{\pi \times_{K_{-1}} \underline{E}K_{-1} \coprod \pi \times_{K_1} \underline{E}K_1}{}{}
{\pi \times_{H_e} \underline{E}H_e \times [-1,1]}{}{\underline{E}\pi}
$$
\em
\end{example}

\begin{example}[The graph associated to an $HNN$-extension] 
\label{exa: HNN-extension} \em
Consider the graph $S$ with one edge $e$ and one vertex $v$.
There is only one choice for the  map $s \colon \{e\} \times \{-1,1\} \to \{v\}$. 
Of course this graph is homeomorphic to $S^1$. Let $\calg$ be a graph of
groups on $S$. It consists of two groups $H_e$ and $K_v$ and two injective group
homomorphisms $\phi_{\sigma} \colon H_e \to K_v$ for $\sigma \in \{-1,1\}$.
There is only one choice of a maximal subtree, namely $\{v\}$. The fundamental group
$\pi$ of $\calg$ is the so called $HNN$-extension associated to the data
$\phi_{\sigma} \colon H_e \to K_v$ for $\sigma \in \{-1,1\}$, i.e.
the group generated by the elements of $K_v$ and a letter $t_v$ whose relations are those
of $K_v$ and the relations $t_v^{-1}\phi_{-1}(h)t_v = \phi_1(h)$ for all $h \in H_e$.
Recall that the natural map $K_v \to \pi$ is injective and we will identify
$K_v$ with its image in $\pi$.
Choose  $\phi_{\sigma}$-equivariant maps $f_{\sigma} \colon \underline{E}H_e \to \underline{E}K_v$.
Let $F_{\sigma} \colon \pi \times_{\phi_{-1}} \underline{E}H_e \to \pi \times\underline{E}K_v$ be
the $\pi$-map which sends $(g,x)$ to $gf_{-1}(x)$ for $\sigma = -1$ and 
to $gt_ef_1(x)$ for $\sigma = 1$. Then a model for $\underline{E}\pi$ is given by
the $\pi$-pushout
$$\comsquare{\pi \times_{\phi_{-1}} \underline{E}H_e \times \{-1,1\}}
{F_{-1} \coprod F_1}{\pi \times_{K_v} \underline{E}K_v}  
{}{}
{\pi \times_{\phi_{-1}} \underline{E}H_e \times [-1,1]}{}{\underline{E}\pi}
$$
Notice that this looks like a telescope construction which is infinite to both sides.
Consider the special case, where $H_e = K_v$, $\phi_{-1} = \id$ and $\phi_1$ is an
automorphism. Then $\pi$ is the semidirect product $K_v \rtimes_{\phi_1} \bbZ$. Choose
a $\phi_1$-equivariant map $f_1 \colon \underline{E}K_v \to \underline{E}K_v$. Then a model for
$\underline{E}\pi$ is given by the to both side infinite mapping telescope
of $f_1$ with the $K_v \rtimes_{\phi_1} \bbZ$ action,
for which $\bbZ$ acts by shifting to the right and $k \in K_v$ acts on the part 
belonging to $n \in \bbZ$ by multiplication with $\phi_1^n(k)$.
If we additionally assume that $\phi_1 =  \id$, then $\pi = K_v \times \bbZ$
and we get $\underline{E}K_v \times \bbR$ as model for $\underline{E}\pi$.
\em
\end{example}  

\begin{remark} \label{rem: graphs of groups and EG} \em
All these constructions yield also models for $EG = \EGF{G}{\caltr}$ if one replaces
everywhere the spaces $\underline{E}H_e$ and $\underline{E}K_v$ by the spaces
$EH_e$ and $EK_v$.
\em
\end{remark}


\subsection{Affine Buildings}
\label{subsec: Affine Buildings}

Let $\Sigma$ be an \emph{affine building},%
\index{building!affine}
sometimes also called  \emph{Euclidean building}.%
\index{building!Euclidean}
This is a simplicial complex together with a system of subcomplexes called 
\emph{apartments}
\index{apartments}
satisfying the following axioms:

\begin{enumerate}

\item  Each apartment is isomorphic to an affine Coxeter complex;

\item  Any two simplices of $\Sigma$ are contained in some common apartment;
 
\item  If two apartments both contain two simplices $A$ and $B$ of $\Sigma$,
then there is an isomorphism of one apartment onto the other which fixes the two
simplices $A$ and $B$ pointwise. 

\end{enumerate}

The precise definition of an affine Coxeter complex, which is sometimes called
also Euclidean Coxeter complex, can be found in 
\cite[Section 2 in Chapter VI]{Brown(1998)}, where also more 
information about affine buildings is given. An affine building comes with metric
$d \colon \Sigma \times \Sigma \to [0,\infty)$ 
which is non-positively curved and complete. The building with this metric is a
CAT(0)-space. A simplicial automorphism
of $\Sigma$ is always an isometry with respect to $d$.
For two points $x,y$ in the affine building
there is a unique line segment $[x,y]$ joining $x$ and $y$.
It is the set of points $\{z \in \Sigma \mid d(x,y) = d(x,z) + d(z,y)\}$.
For $x,y \in \Sigma$ and $t \in [0,1]$ let $tx + (1-t)y$ be the point $z \in \Sigma$
uniquely determined by the property that $d(x,z) = td(x,y)$ and $d(z,y) = (1-t)d(x,y)$.
Then the map 
$$r \colon \Sigma \times \Sigma \times [0,1] \to \Sigma, 
\hspace{5mm} (x,y,t) \mapsto tx + (1-t)y$$ 
is continuous. This implies that
$\Sigma$ is contractible. All these facts are taken from 
\cite[Section 3 in Chapter VI]{Brown(1998)} and
\cite[Theorem 10A.4 on page 344]{Bridson-Haefliger(1999)}.

Suppose that the group $G$ acts on $\Sigma$ by isometries.
If $G$ maps a non-empty bounded subset $A$ of $\Sigma$ to itself,
then the $G$-action has a fixed point 
\cite[Theorem 1 in Section 4 in Chapter VI on page 157]{Brown(1998)}.
Moreover the $G$-fixed point set must be contractible since for two points $x,y \in \Sigma^G$
also the segment $[x,y]$ must lie in $\Sigma^G$ and hence
the map $r$ above induces a continuous map
$\Sigma^G \times \Sigma^G \times[0,1] \to \Sigma^G$. This implies
together with Theorem~\ref{the: G-homotopy characterization of EGF(G)(calf)}
\ref{the: G-homotopy characterization of EGF(G)(calf): characterization},
Example \ref{exa: simplicial actions}, 
Lemma \ref{lem: E_FG = J_FG if all H in calf are open} and
Lemma \ref{lem: totally disconnected groups and calcom}

\begin{theorem}[Affine buildings]
\label{the: affine buildings}
\index{Theorem!Affine buildings}
Let $G$ be a topological (locally compact Hausdorff group).
Suppose that $G$ acts on the affine building by simplicial automorphisms
such that each isotropy group is compact. Then each isotropy group is compact open,
$\Sigma$ is a model
for $\JGF{G}{\calcomop}$  and the barycentric subdivision $\Sigma^{\prime}$
is a model for both $\JGF{G}{\calcomop}$ and  $\EGF{G}{\calcomop}$.
If we additionally assume that $G$ is totally disconnected or is a Lie group, then
$\Sigma$ is a model for both $\underline{J}G$ and $\underline{E}G$.
\end{theorem}

\begin{example}[Bruhat-Tits building] \label{exa: Bruhat-Tis building} \em
An important example is the case of a reductive $p$-adic algebraic group $G$ and 
its associated affine Bruhat-Tits building $\beta(G)$
\cite{Tits(1974)},\cite{Tits(1979)}.
Then $\beta(G)$ is a model for $\underline{J}G$ and
$\beta(G)^{\prime}$ is a model for $\underline{E}G$
by Theorem \ref{the: affine buildings}.
\end{example}


\subsection{The Rips Complex of a Word-Hyperbolic Group}
\label{subsec: The Rips Complex of a Word-Hyperbolic Group}

Let $G$ be a finitely generated discrete group. Let $S$ be a finite set of
generators. We will always assume that $S$ is \emph{symmetric}, %
\index{symmetric set of generators} 
i.e.\ that the identity element $1 \in G$ does not belong to $S$ and 
$s \in S$ implies $s^{-1} \in S$. For $g_1,g_2 \in G$  let $d_S(g_1,g_2)$
be the minimal natural number $n$ such that $g_1^{-1}g_2$ can be
written as a word $s_1s_2\ldots s_n$. This defines a left
$G$-invariant metric on $G$, the so called \emph{word metric}.%
\index{word metric} 

A metric space $X = (X,d)$ is called \emph{$\delta$-hyperbolic}%
\index{delta-hyperbolic@$\delta$-hyperbolic!metric space}
for a given real number $\delta \ge 0$ if for any four points $x,y,z,t$ the
following inequality holds
\begin{eqnarray}
d(x,y) + d(z,t) ~ \le ~  \max\{d(x,z) + d(y,t),d(x,t) + d(y,z)\} + 2\delta.
\label{inequality defining delta-hyperbolic}
\end{eqnarray}

A group $G$ with a finite symmetric set $S$ of generators is called
\emph{$\delta$-hyperbolic}%
\index{delta-hyperbolic@$\delta$-hyperbolic!group} if the metric space $(G,d_S)$ is
$\delta$-hyperbolic. 

The \emph{Rips complex}%
\index{Rips complex}
$P_d(G,S)$%
\indexnotation{P_d(G,S)}
of a group $G$ with a symmetric finite set $S$ of generators for a
natural number $d$ is the geometric realization of the simplicial set whose set of
$k$-simplices consists of $(k+1)$-tuples $(g_0,g_1, \ldots g_k)$ of
pairwise distinct elements $g_i \in G$ satisfying $d_S(g_i,g_j) \le d$
for all $i,j \in \{0,1,\ldots ,k\}$. The obvious
$G$-action by simplicial automorphisms on $P_d(G,S)$ induces a
$G$-action by simplicial automorphisms on the barycentric
subdivision $P_d(G,S)^{\prime}$ (see Example
\ref{exa: simplicial actions}). The following result is 
proved in \cite{Meintrup(2000)}, \cite{Meintrup-Schick(2002)}.

\begin{theorem}[Rips complex] 
\label{the: Rips complex}
\index{Theorem!Rips complex}
Let $G$ be a (discrete) group with a finite symmetric set of generators.
Suppose that $(G,S)$ is $\delta$-hyperbolic for the real number
$\delta \ge 0$. Let $d$ be a natural number with $d \ge 16\delta +8$. 
Then the barycentric subdivision  of
the Rips complex $P_d(G,S)^{\prime}$ is a finite $G$-$CW$-model for
$\underline{E}G$.
\end{theorem}

A metric space is called \emph{hyperbolic}%
\index{hyperbolic!metric space}
\index{metric space!hyperbolic}
if it is
$\delta$-hyperbolic for some real number $\delta \ge 0$. 
A finitely generated group $G$ is called \emph{hyperbolic}%
\index{hyperbolic!group}
\index{group!hyperbolic}
if for one (and hence all) finite symmetric set $S$ of generators the
metric space $(G,d_S)$ is a hyperbolic metric space.  Since for metric spaces the property hyperbolic
is invariant under quasiisometry and for two symmetric finite sets
$S_1$ and $S_2$ of generators of $G$ the metric spaces
$(G,d_{S_1})$ and $(G,d_{S_2})$ are quasiisometric, the choice of $S$
does not matter. Theorem \ref{the: Rips complex} implies
that for a hyperbolic group there is a finite $G$-$CW$-model for
$\underline{E}G$.

The notion of a hyperbolic group is due to Gromov and has intensively
been studied (see for example \cite{Bridson-Haefliger(1999)}, \cite{Ghys-Harpe(1990)}, \cite{Gromov(1987)}).
The prototype is the fundamental group of a closed hyperbolic
manifold.


\subsection{Arithmetic Groups}
\label{subsec: Arithmetic Groups}
Arithmetic groups in a semisimple connected linear $\bbQ$-algebraic group 
possess finite models for
$\underline{E}G$. Namely, let $G(\bbR)$ be the $\bbR$-points of a semisimple
$\bbQ$-group $G(\bbQ)$  and let $K\subseteq G(\bbR)$ a maximal compact subgroup.
If $A \subseteq G(\bbQ)$ is an arithmetic group, then $G(\bbR)/K$ with the left 
$A$-action is a model for $\EGF{A}{\calfin}$ as already explained in 
Theorem \ref{the: Discrete subgroups of almost connected Lie groups}. 
The $A$-space $G(\bbR)/K$ is not necessarily cocompact.
The Borel-Serre completion  of  $G(\bbR)/K$
  (see \cite{Borel-Serre(1973)}, \cite{Serre(1979)}) is a finite $A$-$CW$-model 
for $\EGF{A}{\calfin}$ as pointed out in \cite[Remark 5.8]{Adem-Ruan(2001)},
where a private communication with Borel and Prasad is mentioned.


\subsection{Outer Automorphism Groups of Free groups}
\label{subsec: Outer Automorphism Groups of Free groups}

Let $F_n$%
\indexnotation{free group in n letters}
be the free group of rank $n$. Denote by 
$\Out(F_n)$%
\indexnotation{outer group of free group} 
the group of outer automorphisms
of $F_n$, i.e.\ the quotient of the group of all automorphisms of $F_n$ by the normal subgroup of inner
automorphisms. 
Culler and Vogtmann \cite{Culler-Vogtmann(1986)}, 
\cite{Vogtmann(2003)}
have constructed a space $X_n$ called \emph{outer space}%
\index{outer space}
on which $\Out(F_n)$
acts with finite isotropy groups. It is analogous to the 
Teich\-m\"ul\-ler space of a surface with the action of the mapping class group of the surface.
Fix a graph $R_n$ with one vertex $v$ and $n$-edges and identify $F_n$
with $\pi_1(R_n,v)$.
A \emph{marked metric graph}%
\index{marked metric graph}
 $(g,\Gamma)$ consists of a graph $\Gamma$ with all vertices
of valence at least three, a homotopy equivalence $g \colon R_n \to \Gamma$ called marking and
to every edge of $\Gamma$ there is assigned a positive length 
which makes $\Gamma$ into a metric space by the path metric. We call two marked metric graphs 
$(g,\Gamma)$ and $(g',\Gamma')$ equivalent of there is a homothety 
$h \colon \Gamma \to \Gamma'$ such that $g \circ h$ and $h'$ are homotopic.
Homothety means that there is a constant $\lambda > 0$ with $d(h(x),h(y)) = \lambda \cdot d(x,y)$ for all
$x,y$. Elements in outer space $X_n$ are equivalence classes of marked graphs.
The main result in \cite{Culler-Vogtmann(1986)} is that $X$ is contractible.
Actually, for each  finite subgroup $H \subseteq \Out(F_n)$ the $H$-fixed point set $X_n^H$ is contractible
\cite[Propostion 3.3 and Theorem 8.1]{Krstic-Vogtmann(1993)}, 
\cite[Theorem 5.1]{White(1993)}.

The space $X_n$ contains a \emph{spine}%
\index{outer space!spine of}
$K_n$ which is an $\Out(F_n)$-equivariant deformation retraction.
This space $K_n$ is a simplicial complex of dimension $(2n-3)$
on which the $\Out(F_n)$-action is by simplicial automorphisms
and cocompact. Actually the group of simplicial automorphisms of $K_n$ is $\Out(F_n)$
\cite{Bridson-Vogtmann(2001)}. 
Hence the barycentric subdivision $K_n^{\prime}$ is a finite
$(2n-3)$-dimensional model of $\underline{E}\Out(F_n)$.


\subsection{Mapping Class groups}
\label{subsec: Mapping Class groups}

Let $\Gamma^s_{g,r}$%
\indexnotation{Gamma^s_g,r}
\index{mapping class group}
 be the \emph{mapping class group} of an orientable compact surface
$F$ of genus $g$ with $s$ punctures and $r$ boundary components.
This is the group of isotopy classes of orientation preserving
selfdiffeomorphisms $F_g \to F_g$, which  preserve the punctures individually and restrict to the
identity on the boundary. We require that the isotopies
leave the boundary pointwise fixed. We will always assume
that $2g +s +r > 2$, or, equivalently, that the Euler characteristic
of the punctured surface $F$ is negative. It is well-known that
the associated  \emph{Teich\-m\"ul\-ler space}%
\index{Teichm\"uller space}
$\calt^s_{g,r}$%
\indexnotation{calt^s_g,r}
 is a contractible
space on which  $\Gamma^s_{g,r}$ acts properly. Actually 
$\calt^s_{g,r}$ is a model for $\EGF{\Gamma^s_{g,r}}{\calfin}$ by the results of Kerckhoff
 \cite{Kerckhoff(1983)}.

We could not find a clear reference in the literature for the to
experts known statement that there
exist a finite $\Gamma^s_{g,r}$-$CW$-model for 
 $\EGF{\Gamma^s_{g,r}}{\calfin}$. The work of Harer \cite{Harer(1986)}
on the existence of a spine and the construction of the spaces
$T_S(\epsilon)^H$ due to Ivanov \cite[Theorem 5.4.A]{Ivanov(2002)} seem to lead
to such models.


\subsection{Groups with Appropriate Maximal Finite Subgroups}
\label{subsec: Groups with Appropriate Maximal Finite Subgroups}

Let $G$ be a discrete group.
Let $\calmfin$ be the subset  of $\calfin$ consisting
of elements in $\calfin$ which are maximal in $\calfin$. 
Consider the following assertions concerning $G$:
\begin{itemize}

\item[(M)] Every non-trivial finite subgroup of $G$ is contained in a unique maximal finite subgroup;

\item[(NM)] $M \in \calmfin, M \not= \{1\} ~ \Rightarrow ~ N_GM = M$;

\end{itemize}

For such a group there is a  nice model for $\underline{E}G$ with as few non-free cells as possible.
Let $\{(M_i) \mid i \in I\}$ be the set of conjugacy classes of maximal
finite subgroups of $M_i \subseteq Q$. By attaching free $G$-cells we get
an inclusion of $G$-$CW$-complexes  
$j_1 \colon \coprod_{i \in I} G \times_{M_i} EM_i ~ \to EG$,
where $EG$ is the same as $\EGF{G}{\caltr}$, i.e.\ a contractible free $G$-$CW$-complex.
Define $\underline{E}G$ as the $G$-pushout
\begin{eqnarray}
& \comsquare{\coprod_{i \in I} G \times_{M_i} EM_i}{j_1}{EG}
{u_1}{f_1}{\coprod_{i \in I} G/M_i}{k_1}{\underline{E}G}, &EG
\label{pushout for underline E G}
\end{eqnarray}
where $u_1$ is the obvious $G$-map obtained by collapsing each $EM_i$ to a point. 

We have to explain why $\underline{E}G$ is a model for the classifying space for
proper actions of $G$. Obviously it is a $G$-$CW$-complex. Its  isotropy groups
are all finite. We have to show for $H \subseteq G$ finite
that $(\underline{E}G)^H$ contractible. We begin with the case $H \not= \{1\}$.
Because of conditions (M) and (NM) 
there is precisely one index $i_0 \in I$ such that $H$ is subconjugated to 
$M_{i_0}$ and is not subconjugated to $M_i$ for $i \not= i_0$ and we get
$$\left(\coprod_{i \in I} G/M_i\right)^H  ~ = ~ \left(G/M_{i_0}\right)^H ~ = ~ \pt.$$
Hence $\underline{E}G^H = \pt$.
It remains to treat $H = \{1\}$. Since $u_1$ is a non-equivariant homotopy equivalence
and $j_1$ is a cofibration, $f_1$ is a non-equivariant homotopy equivalence and hence
$\underline{E}G$ is contractible (after forgetting the group action).

Here are some examples of groups $Q$ which satisfy conditions (M) and (NM):
\begin{itemize}

\item Extensions $1 \to \bbZ^n \to G \to F \to 1$ for finite $F$ such that the conjugation
  action of $F$ on $\bbZ^n$ is free outside $0 \in \bbZ^n$. \\[1mm]
The conditions (M), (NM)  are satisfied by
\cite[Lemma 6.3]{Lueck-Stamm(2000)}. 

\item Fuchsian groups $F$ \\[1mm]
The conditions (M), (NM) are satisfied
(see for instance \cite[Lemma 4.5]{Lueck-Stamm(2000)}). 
In  \cite{Lueck-Stamm(2000)} the larger class of cocompact planar
groups (sometimes also called cocompact NEC-groups) is treated.

\item One-relator groups $G$\\[1mm]
Let $G$ be a one-relator group. Let $G = \langle (q_i)_{i \in I} \mid r \rangle$ 
be a presentation with one relation.  We only
have to consider the case, where $G$ contains torsion.
Let $F$ be the free group with basis $\{q_i \mid i \in I\}$. Then $r$ is an element in
$F$. There exists an element $s \in F$ and an integer $m \ge 2$
such that $r = s^m$, the cyclic  subgroup $C$
generated by the class $\overline{s} \in G$ represented 
by $s$ has order $m$, any finite subgroup of $G$ is subconjugated to $C$ and 
for any $g \in G$ the implication
$g^{-1}Cg \cap C \not= 1 \Rightarrow g \in C$ holds.
These claims follows from
\cite[Propositions 5.17, 5.18 and 5.19 in II.5 on pages 107 and 108]{Lyndon-Schupp(1977)}.
Hence $G$ satisfies (M) and (NM). 
\end{itemize}


\subsection{One-Relator Groups}
\label{subsec: One-Relator Groups}

Let $G$ be a one-relator group. Let $G = \langle (q_i)_{i \in I} \mid r \rangle$ 
be a presentation with one relation. 
There is up to conjugacy one maximal finite subgroup $C$ which is cyclic.
Let $p \colon \ast_{i \in I} \bbZ \to G$ be the epimorphism from the free group
generated by the set $I$ to $G$, which sends the generator $i \in I$ to
$q_i$. Let $Y \to \bigvee_{i \in I} S^1$ be the $G$-covering associated to the epimorphism
$p$. There is a $1$-dimensional  unitary $C$-representation $V$
and a $C$-map $f \colon SV \to \res_G^C Y$ such that the following is true.
The induced action on the unit sphere $SV$ is free. If we equip $SV$ and $DV$ with the obvious
$C$-$CW$-complex structures, the $C$-map $f$ can be chosen to be cellular 
and we obtain a $G$-$CW$-model for $\underline{E}G$
by the $G$-pushout
$$\comsquare{G \times_C SV}{\overline{f}}{Y}{}{}{G \times_C DV}{}{\underline{E}G}$$
where $\overline{f}$ sends $(g,x)$ to $gf(x)$. 
Thus we get a $2$-dimensional $G$-$CW$-model for $\underline{E}G$ such that
$\underline{E}G$ is obtained from $G/C$ for a maximal finite cyclic
subgroup $C \subseteq G$ by attaching free cells of dimensions $\le 2$
and the $CW$-complex structure on the quotient $G\backslash \underline{E}G$ 
has precisely one $0$-cell, precisely one $2$-cell and as many $1$-cells as there are
elements in $I$. All these claims follow from
\cite[Exercise 2 (c) II. 5 on page 44]{Brown(1982)}. 

If $G$ is torsionfree, the 
$2$-dimensional complex associated to a  presentation with one relation
is a model for $BG$ (see also \cite[Chapter III \S\S 9 -11]{Lyndon-Schupp(1977)}).


\subsection{Special Linear Groups of (2,2)-Matrices}
\label{subsec: Special Linear Groups of (2,2)-Matrices}

In order to illustrate some of the general statements above we consider
the special example $SL_2(\bbR)$ and $SL_2(\bbZ)$.

Let $\bbH^2$ be the $2$-dimensional hyperbolic space. We will use
either the upper half-plane model or the Poincar\'e disk model. The
group $SL_2(\bbR)$ acts by isometric diffeomorphisms on the upper half-plane by
Moebius transformations, i.e.\ a matrix $\squarematrix{a}{b}{c}{d}$ 
acts by sending a complex number $z$ with positive imaginary part to
$\frac{az +b}{cz +d}$. This action  is proper and transitive. The isotropy group
of $z = i$ is $SO(2)$. Since $\bbH^2$ is a
simply-connected Riemannian manifold, whose sectional curvature is
constant $-1$, the $SL_2(\bbR)$-space $\bbH^2$ is a model for
$\underline{E}SL_2(\bbR)$ by Theorem 
\ref{the: actions on simply connected non-posively curved manifolds}. 

One easily checks that $SL_2(\bbR)$ is a connected Lie group 
and $SO(2) \subseteq SL_2(\bbR)$ is a maximal compact subgroup. 
Hence $SL_2(\bbR)/SO(2)$ is a model for $\underline{E}SL_2(\bbR)$
by Theorem \ref{the: almost connected groups}.
Since the $SL_2(\bbR)$-action on $\bbH^2$ is transitive and
$SO(2)$ is the isotropy group at $i \in \bbH^2$, we see that
the $SL_2(\bbR)$-manifolds $SL_2(\bbR)/SO(2)$ and $\bbH^2$ are
$SL_2(\bbR)$-diffeomorphic.

Since $SL_2(\bbZ)$ is a discrete subgroup of $SL_2(\bbR)$, 
the space $\bbH^2$ with the obvious $SL_2(\bbZ)$-action is a model for 
$\underline{E}SL_2(\bbZ)$ 
(see Theorem \ref{the: Discrete subgroups of almost connected Lie groups}).

The group $SL_2(\bbZ)$ is isomorphic to the amalgamated product 
$\bbZ/4 \ast_{\bbZ/2} \bbZ/6$. From Example \ref{exa: amalgamated products} 
we conclude that a model for $\underline{E}SL_2(\bbZ)$ is given 
by the following $SL_2(\bbZ)$-pushout
$$\comsquare{SL_2(\bbZ)/(\bbZ/2) \times \{-1,1\}}{F_{-1} \coprod F_1}
{SL_2(\bbZ)/(\bbZ/4) \coprod SL_2(\bbZ)/(\bbZ/6)}{}{}
{SL_2(\bbZ)/(\bbZ/2) \times [-1,1]}{}{\underline{E}SL_2(\bbZ)}
$$
where $F_{-1}$ and $F_1$ are the obvious projections.
This model for $\underline{E}SL_2(\bbZ)$ is a tree, which has
alternately two and three edges emanating from each vertex. The other model 
$\bbH^2$ is a manifold. These two models must be $SL_2(\bbZ)$-homotopy equivalent.
They can explicitly be related by the following construction.

Divide the Poincar\'e disk into fundamental domains for the $SL_2(\bbZ)$-action.
Each fundamental domain is a geodesic triangle with one vertex at infinity, i.e.
a vertex on the boundary sphere, and two vertices in the interior. Then the union of
the edges, whose end points lie in the interior of the Poincar\'e disk, is a tree
$T$ with $SL_2(\bbZ)$-action. This is the tree model above. The tree is a
$SL_2(\bbZ)$-equivariant deformation retraction of the Poincar\'e disk.
A retraction is given by moving a point $p$ in the Poincar\'e disk along a geodesic
starting at the vertex at infinity, which belongs to the triangle containing $p$,
through $p$ to the first intersection point of this geodesic with $T$.

The tree $T$ above can be identified with the Bruhat-Tits building
of $SL_2(\bbQ\widehat{_p})$ and hence is a model for $\underline{E}SL_2(\bbQ\widehat{_p})$
(see \cite[page 134]{Brown(1998)}).
Since $SL_2(\bbZ)$ is a discrete subgroup of $SL_2(\bbQ\widehat{_p})$,
we get another reason why this tree is a model for $SL_2(\bbZ)$.


\subsection{Manifold Models}
\label{subsec: Manifold Models}

It is an interesting question, whether one can find a model
for $\underline{E}G$ which is a smooth $G$-manifold. One may also ask whether such a manifold model
realizes the minimal dimension for $\underline{E}G$ or whether the action is cocompact.
Theorem \ref{the: actions on simply connected non-posively curved manifolds} gives some information about these questions
for simply connected non-positively curved Riemannian manifolds and
Theorem \ref{the: dim(L/K)} for discrete subgroups of Lie groups with finitely many path components.
On the other hand there exists a virtually torsionfree group $G$
such that $G$ acts properly and cocompactly on a contractible manifold (without
 boundary),
but there is no finite $G$-$CW$-model 
for $\underline{E}G$ \cite[Theorem 1.1] {Davis-Leary(2003)}.


\typeout{--------------------   Section 5: Finiteness Conditions --------------------------}

\section{Finiteness Conditions}
\label{sec: Finiteness Conditions}

In this section we investigate whether there are models for 
$\EGF{G}{\calf}$ which satisfy certain finiteness conditions
such as being finite, being of finite type or being of finite dimension as a $G$-$CW$-complex.


\subsection{Review of Finiteness Conditions on $BG$}
\label{subsec: Review of Finiteness Conditions on $BG$}

As an illustration we review the corresponding question for $EG$ for a discrete group $G$.
This is equivalent to the question whether for a given discrete group $G$ there is a $CW$-complex model for
$BG$ which is finite, of finite type or finite dimensional.  

We introduce the following notation. 
Let $R$ be a commutative associative ring with unit.
The trivial $RG$-module is $R$ viewed as $RG$-module by the trivial $G$-action.
A \emph{projective resolution}%
\index{resolution!resolution}
or \emph{free resolution}%
\index{resolution!free}
respectively for an $RG$-module $M$ is an $RG$-chain complex $P_*$
of projective or free respectively $RG$-modules with $P_i = 0$ for $i \le -1$ such that
$H_i(P_*) = 0$ for $i \ge 1$ and $H_0(P_*)$ is $RG$-isomorphic to $M$.
If additionally each $RG$-module $P_i$ is finitely generated
and $P_*$ is finite dimensional,
we call $P_*$ \emph{finite}.%
\index{resolution!finite}

An $RG$-module $M$ has \emph{cohomological dimension}%
\index{cohomological dimension!for modules}
$\cd(M) \le n$,%
\indexnotation{cd(M)}
if there exists a projective resolution of dimension $\le n$ for $M$.
This is equivalent to the condition that for any $RG$-module $N$
we have $\Ext^i_{RG}(M,N) = 0$ for $i \ge n+1$. 
A group $G$ has \emph{cohomological dimension}%
\index{cohomological dimension! for groups}
$cd(G) \le n$
\indexnotation{cd(G)}
 over $R$ if the trivial $RG$-module $R$ has cohomological dimension $\le n$.
An $RG$-module $M$ is of \emph{type $F\!P_n$},%
\index{type $FP_n$!for modules}
\indexnotation{FP_n}
if it admits a  projective $RG$-resolution $P_*$ such that $P_i$ is finitely generated for
$i \le n$ and of type $F\!P_{\infty}$%
\index{type $FP_{\infty}$!for modules} 
\indexnotation{FP_infty}
if it admits a  projective $RG$-resolution $P_*$  such that $P_i$ is finitely generated for
all $i$. A group $G$ is of type $F\!P_n$%
\index{type $FP_n$!for groups}
or  $F\!P_{\infty}$%
\index{type $FP_{\infty}$!for groups}
respectively if the trivial $\bbZ G$-module $\bbZ$ is of type $F\!P_n$ or $F\!P_{\infty}$ respectively.

Here is a summary of well-known statements about finiteness conditions on $BG$.
A key ingredient in the proof of the next result is the fact
that the cellular $RG$-chain complex $C_*(EG)$ is a free and in particular a projective $RG$-resolution
of the trivial $RG$-module $R$ since $EG$ is a free $G$-$CW$-complex and contractible,
and that $C_*(EG)$ is $n$-dimensional or of type $F\!P_n$ respectively if $BG$ is $n$-dimensional
or has finite $n$-skeleton respectively.

\begin{theorem}[Finiteness conditions for $BG$]
\label{the: BG} 
\index{Theorem!Finiteness conditions for $BG$}
Let $G$ be a discrete group.

\begin{enumerate}

\item  \label{the: BG: finite dimensional model}
If there exists a finite dimensional model for $BG$, then $G$ is torsionfree;

\item \label{the: BG: finite type}
\begin{enumerate}

\item \label{the: BG: finite type: 1}
There exists a $CW$-model for $BG$ with finite $1$-skeleton
if and only if $G$ is finitely generated;

\item \label{the: BG: finite type: 2}
There exists a $CW$-model for $BG$ with finite $2$-skeleton
if and only if $G$ is finitely presented;

\item \label{the: BG: finite type: FP_n}
For $n \ge 3$ there exists a $CW$-model for $BG$ with finite $n$-skeleton
if and only if $G$ is finitely presented and of type $F\!P_n$;

\item \label{the: BG: finite type:  FP_infty}
There exists a $CW$-model for $BG$ of finite type, i.e.\ all skeleta are finite, 
if and only if $G$ is finitely presented
and of type $F\!P_{\infty}$; 

\item \label{the: BG: finite type: Bestvinas counterexample}
There exists groups $G$ which are  of type $F\!P_2$ and which are not finitely presented;

\end{enumerate}

\item \label{the: BG: finite model}
There is a finite $CW$-model for $BG$ if and only if $G$ is finitely
presented and there is a finite  free $\bbZ G$-resolution $F_*$
for the trivial $\bbZ G$-module $\bbZ$;

\item  \label{the: BG: Stallings} 
The following assertions are equivalent:

\begin{enumerate} 

\item \label{the: BG: Stallings: cohodim le 1} 
The cohomological dimension of $G$ is $\le 1$;

\item \label{the: BG: Stallings: 1-dim model for BG} 
There is a model for $BG$ of dimension $\le 1$;

\item \label{the: BG: Stallings: G is free} 
$G$ is free.

\end{enumerate}

\item  \label{the: BG: cohomological characterization}

The following assertions are equivalent for $d \ge 3$:

\begin{enumerate}

\item There exists a $CW$-model for $BG$ of dimension $\le d$;

\item $G$ has cohomological dimension $\le d$ over $\bbZ$;

\end{enumerate}

\item  \label{the: Thompson} For Thompson's group $F$ there is a $CW$-model of finite type for $BG$ 
but no finite dimensional model for $BG$.

\end{enumerate}

\end{theorem}
\begin{proof}
\ref{the: BG: finite dimensional model}
Suppose we can choose a finite dimensional model for $BG$.
Let $C \subseteq G$ be a finite cyclic subgroup.
Then $C\backslash\widetilde{BG} = C\backslash EG$ is a finite dimensional model for $BC$.
Hence there is an integer $d$ such that we have $H_i(BC) = 0$ for $i \ge d$.
This implies that $C$ is trivial \cite[(2.1) in II.3 on page 35]{Brown(1982)}.
Hence $G$ is torsionfree.
\\[1mm]
\ref{the: BG: finite type}
See \cite{Bestvina-Brady(1997)}  and
\cite[Theorem 7.1 in VIII.7 on page 205]{Brown(1982)}. 
\\[1mm]
\ref{the: BG: finite model} 
See \cite[Theorem 7.1 in VIII.7 on page 205]{Brown(1982)}.
\\[1mm]
\ref{the: BG: Stallings} 
See \cite{Stallings(1968)} and \cite{Swan(1969)}.
\\[1mm]
\ref{the: BG: cohomological characterization}
See \cite[Theorem 7.1 in VIII.7 on page 205]{Brown(1982)}. \\[1mm]
\ref{the: Thompson} See \cite{Brown-Geoghegan(1984)}. 
\end{proof}


\subsection{Modules over the Orbit Category}
\label{subsec: Modules over the Orbit Category}

Let $G$ be a discrete group and let $\calf$ be a family of subgroups.
The \emph{orbit category}%
\index{orbit category}
$\Or(G)$%
\indexnotation{Or(G)}
of $G$ is the small category, whose objects are
homogeneous $G$-spaces $G/H$ and whose
morphisms are $G$-maps. Let $\OrGF{G}{\calf}$%
\indexnotation{Or(G,calf)}
be the full subcategory of $\Or(G)$ consisting of those objects $G/H$ for which
$H$ belongs to $\calf$. A
\emph{$\bbZ\OrGF{G}{\calf}$-module}%
\index{ZOr(G,F)-module@$\bbZ\OrGF{G}{\calf}$-module}
is a contravariant functor from
$\OrGF{G}{\calf}$ to the category of $\bbZ$-modules.
A morphism of such modules is a natural
transformation. The category of $\bbZ\OrGF{G}{\calf}$-modules
inherits the structure of an abelian
category from the standard structure of an
abelian category on the category of
$\bbZ$-modules. In particular the notion
of a projective $\bbZ\OrGF{G}{\calf}$-module is defined.
The {\it free} $\bbZ\OrGF{G}{\calf}$-module
$\bbZ\map(G/?,G/K)$ \emph{based at the object}%
\index{ZOr(G,F)-module@$\bbZ\OrGF{G}{\calf}$-module!free}
 $G/K$ is the $\bbZ\OrGF{G}{\calf}$-module that
assigns to an object $G/H$ the
free $\bbZ$-module $\bbZ\map_G(G/H,G/K)$ generated
by the set $\map_G(G/H,G/K)$.
The key property of it is that for any
$\bbZ\OrGF{G}{\calf}$-module $N$ there is a natural
bijection of $\bbZ$-modules
$$
\hom_{\bbZ\OrGF{G}{\calf}}(\bbZ\map_G(G/?,G/K), N) \xrightarrow{\cong} N(G/K),
\quad \phi \mapsto \phi(G/K)(\id_{G/K}).
$$
This is a direct consequence of the Yoneda Lemma.
A $\bbZ\OrGF{G}{\calf}$-module is {\it free}
if it is isomorphic to a direct sum
$\bigoplus_{i \in I} \bbZ\map(G/?,G/K_i)$ for
appropriate choice of objects $G/K_i$
and index set $I$. A $\bbZ\OrGF{G}{\calf}$-module
is called {\it finitely generated} if
it is a quotient of a $\bbZ\OrGF{G}{\calf}$-module
of the shape  $\bigoplus_{i \in I} \bbZ\map(G/?,G/K_i)$
with a finite index set $I$.
Notice that a lot of standard facts for $\bbZ$-modules carry over to
 $\bbZ\OrGF{G}{\calf}$-modules. For instance,
 a $\bbZ\OrGF{G}{\calf}$-module is projective or
finitely generated projective respectively
if and only if it is a direct summand in a free
$\bbZ\OrGF{G}{\calf}$-module or a finitely generated
free $\bbZ\OrGF{G}{\calf}$-module respectively.
The notion of a {\it projective resolution}
$P_*$ of a $\bbZ\OrGF{G}{\calf}$-module is obvious
and notions like of cohomological dimension $\le n$ or
of type $F\!P_{\infty}$ carry directly over. 
Each $\bbZ\OrGF{G}{\calf}$-module has a projective resolution.
The trivial $\bbZ \OrGF{G}{\calf}$-module
${\underline \bbZ}$%
\indexnotation{underline{Z}}
is the constant functor 
from $\OrGF{G}{\calf}$ to the category of $\bbZ$-modules, which sends
any morphism to $\id \colon \bbZ \to \bbZ$. More information about
modules over a category can be found for instance in
\cite[Section 9]{Lueck(1989)}.

The next result is 
proved in \cite[Theorem 0.1]{Lueck-Meintrup(2000)}.
A key ingredient in the proof of the next result is the fact
that the cellular $R\OrGF{G}{\calf}$-chain complex $C_*(\EGF{G}{\calf})$ 
is a free and in particular a projective $R\OrGF{G}{\calf}$-resolution
of the trivial $R\OrGF{G}{\calf}$-module $R$.

\begin{theorem}[Algebraic and geometric finiteness conditions]
\index{Theorem!Algebraic and geometric finiteness conditions}
\label{the: algebraic criterion}
Let $G$ be a discrete group and let $d \geq 3$. Then we have:
\begin{enumerate}
\item  \label{the: algebraic criterion: finite dimensional model}
There is $G$-$CW$-model of dimension $\le d$ for $\EGF{G}{\calf}$ if
and only if the trivial $\bbZ \OrGF{G}{\calf}$-module
${\underline \bbZ}$  has cohomological dimension $\le d$;

\item  \label{the: algebraic criterion: of finite type}
There is a $G$-$CW$-model for $\EGF{G}{\calf}$ of finite type if
and only if $\EGF{G}{\calf}$ has a $G$-$CW$-model with finite 2-skeleton and
the trivial $\bbZ \OrGF{G}{\calf}$-module
${\underline \bbZ}$  is of type $F\!P_{\infty}$;

\item \label{the: algebraic criterion: finite model}
There is a finite $G$-$CW$-model for $\EGF{G}{\calf}$ if
and only if $\EGF{G}{\calf}$ has a $G$-$CW$-model with finite 2-skeleton and
the trivial  $\bbZ \OrGF{G}{\calf}$-module ${\underline \bbZ}$ has a finite free resolution over
$\OrGF{G}{\calf}$;

\item \label{the: algebraic criterion: finite 2-skeleton}
There is a $G$-$CW$-model with finite $2$-skeleton for $\underline{E}G = \EGF{G}{\calfin}$
if and only if
there are only finitely many conjugacy classes of finite subgroups $H \subset G$
and for any finite subgroup $H \subset G$ its Weyl group $W_GH := N_GH/H$ is finitely presented.
\end{enumerate}
\end{theorem}

In the case, where we take $\calf$ to be the trivial family,
Theorem \ref{the: algebraic criterion} 
\ref{the: algebraic criterion: finite dimensional model} reduces to 
Theorem \ref{the: BG} \ref{the: BG: cohomological characterization},
Theorem \ref{the: algebraic criterion} 
\ref{the: algebraic criterion: of finite type} to
Theorem \ref{the: BG} \ref{the: BG: finite type:  FP_infty} and
Theorem \ref{the: algebraic criterion} 
\ref{the: algebraic criterion: finite model} to
Theorem \ref{the: BG}
\ref{the: BG: finite model},
and one should compare 
Theorem  \ref{the: algebraic criterion} 
\ref{the: algebraic criterion: finite 2-skeleton} to
Theorem  \ref{the: BG} \ref{the: BG: finite type: 2}.

\begin{remark} \label{rem: Nucinkis} \em 
Nucinkis \cite{Nucinkis(2000)} investigates the notion of
$\calfin$-cohomological dimension and relates it to the question
whether there are finite dimensional modules for $\underline{E}G$.
It gives another lower bound for the dimension of a model for
$\underline{E}G$ but is not sharp in general \cite{Brady-Leary-Nucinkis(2001)}. \em
\end{remark}


\subsection{Reduction from Topological Groups to Discrete Groups}
\label{subsec: Reduction from Topological Groups to Discrete Groups}

The \emph{discretization} $G_d$%
\indexnotation{G_d}
of a topological group
$G$ is  the same group but now with the discrete topology. 
Given a family $\calf$ of (closed) subgroups of $G$, denote by $\calf_d$
the same set of subgroups, but now in connection with $G_d$. Notice that $\calf_d$ is again a
family. We will need the following condition
\\[4mm]
(S) \hspace*{7mm} \begin{minipage}{104mm}
For any closed subgroup $H \subset G$ the projection
$p\colon G \to G/H$ has a local cross section, i.e.\ there is a neighborhood
$U$ of $eH$ together with a map $s\colon U \to G$ satisfying $p \circ s = \id_U$.
\end{minipage}
\\[4mm]
Condition (S) is automatically satisfied
if $G$ is discrete, if $G$ is a Lie group, or
more generally, if $G$ is locally compact and second countable and has finite
covering dimension \cite{Mostert(1953)}.
The metric needed in
\cite{Mostert(1953)} follows under our assumptions, since
a locally compact Hausdorff space is regular and regularity in a second
countable space implies metrizability.

The following two results are proved in 
\cite[Theorem 0.2 and Theorem 0.3]{Lueck-Meintrup(2000)}.

\begin{theorem}[Passage from totally disconnected groups to discrete groups]
\index{Theorem!Passage from totally disconnected groups to discrete groups}
 \label{the: totally disconnected if and only if discrete}
Let $G$ be a locally compact totally disconnected Hausdorff group
and let $\calf$ be a family of subgroups of
$G$. Then there is a  $G$-$CW$-model
for $\EGF{G}{\calf}$ that is $d$-dimensional or finite or  of finite type respectively
if and only if there is a $G_d$-$CW$-model for $\EGF{G_d}{\calf_d}$
 that is $d$-dimensional or finite or of finite type respectively.
\end{theorem}

\begin{theorem}[Passage from topological groups to totally disconnected groups]
\index{Theorem!Passage from topological groups to totally disconnected groups}
\label{the: from general to totally disconnected}
Let $G$ be a locally compact Hausdorff group satisfying condition (S).  
Put $\overline{G} := G/G^0$. Then there is a $G$-$CW$-model for
$\underline{E}G$ that is $d$-dimensional or finite or of finite type respectively
if and only if $\underline{E}\overline{G}$ has a
$\overline{G}$-$CW$-model
that is $d$-dimensional or finite or of finite type respectively.
\end{theorem}

If we combine Theorem \ref{the: algebraic criterion},
Theorem \ref{the: totally disconnected if and only if discrete}
and Theorem \ref{the: from general to totally disconnected}
we get

\begin{theorem}[Passage from topological groups to discrete groups]
\index{Theorem!Passage from topological groups to discrete groups} 
\label{the: reduction from topological to discrete}
Let $G$ be a locally compact group satisfying $(S)$. Denote by
$\overline{\calcom}$ the family of compact subgroups of its component group $\overline{G}$
and let $d \geq 3$. Then

\begin{enumerate}

\item There is a $d$-dimensional $G$-$CW$-model for $\underline{E}G$ if
and only if the trivial $\bbZ \OrGF{\overline{G}_d}{\overline{\calcom}_d}$-module
${\underline \bbZ}$ has cohomological dimension $\le d$;

\item  There is a $G$-$CW$-model for $\underline{G}$ of finite type if
and only if $\EGF{\overline{G}_d}{\overline{\calcom}_d}$ has a $\overline{G}_d$-$CW$-model with
finite 2-skeleton and the trivial $\bbZ \OrGF{\overline{G}_d}{\overline{\calcom}_d}$-module
${\underline \bbZ}$ is of type $F\!P_{\infty}$;

\item There is a finite $G$-$CW$-model for $\underline{E}G$ if
and only if $\EGF{\overline{G}_d}{\overline{\calcom}_d}$ has a $\overline{G}_d$-$CW$-model with
finite 2-skeleton and the trivial $\bbZ \OrGF{\overline{G}_d}{\overline{\calcom}_d}$-module
${\underline \bbZ}$ has a finite free resolution.
\end{enumerate}
\end{theorem}

In particular we see from
Theorem \ref{the: from general to totally disconnected} that,
for a Lie group $G$, type questions about $\underline{E}G$
are equivalent to the corresponding type questions of
$\underline{E}\pi_0(G)$, since $\pi_0(G)=\overline{G}$ is discrete.
In this case the family $\overline{\calcom}_d$ appearing in Theorem 
\ref{the: reduction from topological to discrete}.
is just the family $\calfin$ of finite subgroups of $\pi_0(G)$.


\subsection{Poset of Finite Subgroups}
\label{subsec: Poset of Finite Subgroups}

Throughout this Subsection \ref{subsec: Poset of Finite Subgroups} let $G$ be a discrete group.
Define the $G$-poset \begin{eqnarray}
\calp(G) & := & \{K \mid K \subset G \mbox{ finite}, K \not= 1\}.
\label{definition of Gamma-poset}
\end{eqnarray}
An element $g \in G$ sends $K$ to $g K g^{-1}$
and the poset-structure comes from inclusion of subgroups.
Denote by $|\calp(G)|$ the geometric realization of the category given
by the poset $\calp(G)$. This is a $G$-$CW$-complex
but in general not proper, i.e.\ it can have points with infinite isotropy
groups.\par

Let $N_GH$%
\indexnotation{N_GH}
be the \emph{normalizer} and let
$W_GH := N_GH/H$ be the \emph{Weyl group}%
\indexnotation{W_GH}
of $H \subset G$. Notice for a $G$-space $X$
that $X^H$ inherits a $W_GH$-action.
Denote by $CX$%
\indexnotation{CX}
the \emph{cone} over $X$. Notice that
$C\emptyset$ is the one-point-space.

If $H$ and $K$ are subgroups of $G$ and $H$ is finite, then
$G/K^H$ is a finite union of
$W_GH$-orbits of the shape $W_GH/L$ for finite $L \subset W_GH$.
Now one easily checks

\begin{lemma} \label{E(Gamma,Fin)^H = E(WH,FIN)}
The $W_GH$-space $\underline{E}G^H$ is a $W_GH$-$CW$-model
for $\underline{E}W_GH$. In particular, if $\underline{E}G$
has a $G$-$CW$-model which is finite, of finite type
or $d$-dimensional respectively, then there is a $W_GH$-model for
$\underline{E}W_GH$ which is finite, of finite
type or $d$-dimensional respectively. 
\end{lemma}

\begin{notation}[The condition $b(d)$ and $B(d)$]
 \label{notation of b(d) and B(d)}
Let $d \ge 0$ be an integer.
A group $G$ satisfies the condition
$b(d)$ or $b(<\!\infty)$ respectively
if any $\bbZ G$-module $M$ with
the property that $M$ restricted to $\bbZ K$ is projective for all
finite subgroups $K \subset G$ has a  projective
$\bbZ G$-resolution of dimension $d$ or of finite dimension respectively.
A group $G$ satisfies the condition $B(d)$
if  $W_GH$ satisfies the condition $b(d)$
for any finite subgroup $H \subset G$.

The \emph{length}%
\index{length of a subgroups}
$l(H)%
\indexnotation{l(H)}
\in \{0,1, \ldots\}$
of a finite group $H$ is the supremum
over all $p$ for which there is a nested sequence
$H_0 \subset H_1 \subset \ldots \subset H_p$ of subgroups
$H_i$ of $H$ with $H_i \not= H_{i+1}$.
\end{notation}

\begin{lemma} \label{necessity of B(d)}
Suppose that there is a $d$-dimensional $G$-$CW$-complex $X$ with
finite isotropy groups such that $H_p(X;\bbZ) = H_p(*,\bbZ)$ for all $p \ge 0$
holds. This assumption is for instance satisfied if there is a
$d$-dimensional $G$-$CW$-model for $\underline{E}G$.
Then $G$ satisfies condition $B(d)$.
\end{lemma}
\begin{proof} Let $H \subset G$ be finite.
Then $X/H$ satisfies $H_p(X/H;\bbZ) = H_p(*,\bbZ)$ for all $p \ge 0$
\cite[III.5.4 on page 131]{Bredon(1972)}.  Let $C_*$
be the cellular $\bbZ W_GH$-chain complex of $X/H$.
This is a $d$-dimensional
resolution of the trivial $\bbZ W_GH$-module $\bbZ$
and each chain module is a sum of $\bbZ W_GH$-modules of the
shape $\bbZ[W_GH/K]$ for some finite subgroup
$K \subset W_GH$. Let $N$ be a $\bbZ W_GH$-module
such that $N$ is projective over $\bbZ K$ for
any finite subgroup $K \subset W_GH$.
Then $C_*\otimes_{\bbZ} N$ with the diagonal
$W_GH$-operation is a $d$-dimensional projective $\bbZ W_GH$-resolution of
$N$. 
\end{proof}

\begin{theorem}[An algebraic criterion for finite dimensionality]
\label{Theorem!An algebraic criterion for finite dimensionality}
 \label{the: a criterion for dim(e(Gamma,Fin)^H) = d(H)}
Let $G$ be a discrete group. Suppose that we have for any finite subgroup
$H \subset G$ an integer $d(H) \ge 3$
such that $d(H) \ge d(K)$ for $H \subset K$ and $d(H) = d(K)$ if
$H$ and $K$ are conjugate in $G$.
Consider the following statements:

\begin{enumerate}

\item \label{the: a criterion for dim(e(Gamma,Fin)^H) = d(H): existence of a model}
There is a $G$-$CW$-model $\underline{E}G$ such that
for any finite subgroup $H \subset G$
$$\dim(\underline{E}G^H) ~ = ~ d(H);$$

\item \label{the: a criterion for dim(e(Gamma,Fin)^H) = d(H): cohomological condition}
We have for any finite subgroup $H \subset G$
and for any $\bbZ W_GH$-module M
$$H^{d(H) + 1}_{\bbZ W_GH}
(EW_GH \times (C|\calp(W_GH)|,|\calp(W_GH)|);M)
~ = ~ 0;$$

\item \label{the: a criterion for dim(e(Gamma,Fin)^H) = d(H): cohomological condition and B}
We have for any finite subgroup $H \subset G$
that its Weyl group $W_GH$ satisfies $b(<\! \infty)$ and
that there is a subgroup
$\Delta(H) \subset W_GH$ of finite index such that for any
$\bbZ\Delta(H)$-module $M$
$$H^{d(H) + 1}_{\bbZ \Delta(H)}
(E\Delta(H) \times (C|\calp(W_GH)|,|\calp(W_GH)|);M)
~ = ~ 0.$$

\end{enumerate}

Then \ref{the: a criterion for dim(e(Gamma,Fin)^H) = d(H): existence of a model}
implies both 
\ref{the: a criterion for dim(e(Gamma,Fin)^H) = d(H): cohomological condition} and
\ref{the: a criterion for dim(e(Gamma,Fin)^H) = d(H): cohomological condition and B}.
If there is an upper bound  on the length
$l(H)$ of the finite subgroups
$H$ of $G$, then these statements
\ref{the: a criterion for dim(e(Gamma,Fin)^H) = d(H): existence of a model},
\ref{the: a criterion for dim(e(Gamma,Fin)^H) = d(H): cohomological condition}
and
\ref{the: a criterion for dim(e(Gamma,Fin)^H) = d(H): cohomological condition and B} are equivalent.
\end{theorem}

The proof of Theorem \ref{the: a criterion for dim(e(Gamma,Fin)^H) = d(H)} 
can be found in \cite[Theorem 1.6]{Lueck(2000a)}. 
In the case that $G$ has finite virtual cohomological dimension
a similar result is proved in
 \cite[Theorem III]{Connolly-Kozniewski(1986)}.

\begin{example} \label{exa: torsionfree Gamma} \em
Suppose that $G$ is torsionfree. Then
Theorem \ref{the: a criterion for dim(e(Gamma,Fin)^H) = d(H)}
reduces to the well-known result
\cite[Theorem VIII.3.1 on page 190,Theorem VIII.7.1 on page 205]{Brown(1982)}
that the following assertions are equivalent
for an integer $d \ge 3$:
\begin{enumerate}

\item There is a $d$-dimensional $CW$-model for $BG$;

\item $G$ has cohomological dimension $\le d$;

\item $G$ has virtual cohomological dimension $\le d$.

\end{enumerate}
\em
\end{example}

\begin{remark} \em \label{rem: poset if WH contains a finite normal subgroup}
If $W_GH$ contains a non-trivial normal finite subgroup
$L$, then $|\calp(W_GH)|$ is contractible and
\begin{eqnarray*}
H^{d(H) + 1}_{\bbZ W_GH}(EW_GH \times (C|\calp(W_GH)|,|\calp(W_GH)|);M)
& = & 0;
\\
H^{d(H) + 1}_{\bbZ \Delta(H)}
(E\Delta(H) \times (C|\calp(W_GH)|,|\calp(W_GH)|);M)
& = & 0.
\end{eqnarray*}
The proof of this fact is given in \cite[Example 1.8]{Lueck(2000a)}. 
\em
\end{remark}

The next result is taken from \cite[Theorem 1.10]{Lueck(2000a)}. 
A weaker version of it for certain classes of groups and
in $l$ exponential dimension estimate can be found in 
\cite[Theorem B]{Kropholler-Mislin(1998)}
(see \cite[Remark 1.12]{Lueck(2000a)}). 
\begin{theorem}[An upper bound on the dimension]
\label{Theorem!An upper bound on the dimension}
\label{the: sufficient condition B(d)}
Let $G$ be a group and let $l \ge 0$ and $d \ge 0$ be integers
such that the length $l(H)$ of any finite subgroup $H \subset G$
is bounded by $l$ and $G$ satisfies $B(d)$. Then there is a
$G$-$CW$-model for $\underline{E}G$ such that for any finite subgroup
$H \subset G$
$$\dim(\underline{E}G^H) ~ \le   ~ \max\{3,d\} + (l - l(H))(d+1)$$
holds. In particular $\underline{E}G$ has dimension at most
$\max\{3,d\} + l(d+1)$.
\end{theorem}


\subsection{Extensions of Groups}
\label{subsec: Extensions of Groups}

In this subsection we consider an exact sequence of discrete groups
$1 \to \Delta \to G \to \pi \to 1$. We want to investigate whether
finiteness conditions about the type of a classifying space
for $\calfin$ for $\Delta$ and $\pi$ carry over to the one of $G$.
The proof of the next Theorem  
\ref{the: dimension of underline{E} and exact sequences, original version}
is  taken from  \cite[Theorem 3.1]{Lueck(2000a)}), the proof of Theorem 
\ref{the: dimension of underline{E} and exact sequences}
is an easy variation.

\begin{theorem}[Dimension bounds and extensions]
\label{Theorem!Dimension bounds and extensions}
 \label{the: dimension of underline{E} and exact sequences, original version}
Suppose that there exists a positive integer $d$
which is an upper bound on the orders of finite subgroups
of $\pi$. Suppose that $\underline{E}\Delta$ has a $k$-dimensional
$\Delta$-$CW$-model and $\underline{E}\pi$ has a $m$-dimensional
$\pi$-$CW$-model.  Then $\underline{E}G$ has a
$(dk + m)$-dimensional $G$-$CW$-model.
\end{theorem}

\begin{theorem} \label{the: dimension of underline{E} and exact sequences}
Suppose that $\Delta$ has the property that for any
group $\Gamma$ which contains $\Delta$ as subgroup of finite index,
there is a $k$-dimensional $\Gamma$-$CW$-model for $\underline{E}\Gamma$.
Suppose that $\underline{E}\pi$ has a $m$-dimensional
$\pi$-$CW$-model.  Then $\underline{E}G$ has a
$(k + m)$-dimensional $G$-$CW$-model.
\end{theorem}

We will see in Example \ref{exa: virtually poly-cyclic} 
that the condition about $\Delta$ in 
Theorem \ref{the: dimension of underline{E} and exact sequences}
is automatically satisfied if $\Delta$ is virtually poly-cyclic.

The next two results are taken from 
\cite[Theorem 3.2 and Theorem 3.3]{Lueck(2000a)}).

\begin{theorem}
\label{the: finiteness or finite type of underline{E} and exact sequences}
Suppose for any finite subgroup $\pi' \subset \pi$
and any extension
$1 \to \Delta \to \Delta' \to \pi' \to 1$ that $\underline{E}\Delta'$
has a finite $\Delta'$-$CW$-model or
a $\Delta'$-$CW$-model of finite type respectively and suppose
that $\underline{E}\pi$ has a finite $\pi$-$CW$-model or
a $\pi$-$CW$-model of finite type respectively.  Then
$\underline{E}G$ has a finite $G$-$CW$-model or
a $G$-$CW$-model of finite type respectively.
\end{theorem}

\begin{theorem} 
\label{the: finiteness or finite type of underline{E} and exact sequences for word-hyperbolic Delta}
Suppose that $\Delta$ is word-hyperbolic
or virtually poly-cyclic. Suppose
that $\underline{E}\pi$ has a finite $\pi$-$CW$-model or
a $\pi$-$CW$-model of finite type respectively.  Then
$\underline{E}G$ has a finite $G$-$CW$-model or
a $G$-$CW$-model of finite type respectively. 
\end{theorem}


\subsection{One-Dimensional Models for \underline{E}G}
\label{subsec: One-Dimensional Models for underline E G}

The following result follows from Dunwoody \cite[Theorem 1.1]{Dunwoody(1979)}.

\begin{theorem}[A criterion for $1$-dimensional models]
\index{Theorem!A criterion for $1$-dimensional models}
 \label{the: Dunwoody's characterization}
Let $G$ be a discrete group. Then there exists a $1$-dimensional model
for $\underline{E}G$ if and only the cohomological dimension of $G$ over 
the rationals $\bbQ$ is less or equal to one.
\end{theorem}

If $G$ is finitely generated, then there is a $1$-dimensional model
for $\underline{E}G$ if and only if $G$ contains a finitely generated free subgroup
of finite index \cite[Theorem 1]{Karrass-Pietrowski-Solitar(1973)}. If $G$ is torsionfree,
we rediscover the results due to Swan and Stallings stated in 
Theorem \ref{the: BG} \ref{the: BG: Stallings} from Theorem 
\ref{the: Dunwoody's characterization}.


\subsection{Groups of Finite Virtual Dimension}
\label{subsec: Groups of Finite Virtual Dimension}

In this section we investigate the condition $b(d)$ and $B(d)$ of
Notation \ref{notation of b(d) and B(d)} for a discrete group $G$
and explain how our results specialize
in the case of a group of finite virtual cohomological dimension.

\begin{remark} \label{rem: underlineEG finite does not imply virtually torsionfree} \em
There exists groups $G$ with a finite dimensional  model for $\underline{E}G$,
which do not admit a torsionfree subgroup of finite index.
For instance, let $G$ be a countable locally finite group which is not finite.
Then its cohomological dimension over the rationals is $\le 1$ and
hence it  possesses a $1$-dimensional model for $\underline{E}G$ 
by Theorem \ref{the: Dunwoody's characterization}. Obviously it contains no 
torsionfree subgroup of finite index. An example of a group $G$ 
with a finite $2$-dimensional  model for $\underline{E}G$,
which does not admit a torsionfree subgroup of finite index, is described in 
\cite[page 493]{Brady-Leary-Nucinkis(2001)}. \em
\end{remark}

A discrete group $G$ has
\emph{virtual cohomological dimension $\le d$}%
\index{virtual cohomological dimension}
if and only if it contains a torsionfree subgroup $\Delta$ of finite index
such that $\Delta$ has cohomological dimension $\le d$. This is independent
of the choice of $\Delta \subseteq G$ because for two torsionfree subgroups
$\Delta, \Delta' \subseteq G$ we have that $\Delta$ has cohomological dimension
$\le d$ if and only if $\Delta'$ has cohomological dimension $\le d$.
The next two results are taken from 
\cite[Lemma 6.1, Theorem 6.3, Theorem 6.4]{Lueck(2000a)}. 

\begin{lemma} \label{lem: B(d) and subgroups}
If $G$ satisfies $b(d)$ or $B(d)$ respectively, then any subgroup
$\Delta$ of $G$ satisfies $b(d)$ or $B(d)$ respectively.
\end{lemma}

\begin{theorem}[Virtual cohomological dimension and the condition $B(d)$]
\index{Theorem!Virtual cohomological dimension and the condition $B(d)$}
\label{B(d) and vcd}
If $G$ contains a torsionfree subgroup $\Delta$ of finite index,
then the following assertions are equivalent:
\begin{enumerate}

\item $G$ satisfies $B(d)$;

\item $G$ satisfies $b(d)$;

\item $G$ has virtual cohomological dimension $\le d$.

\end{enumerate}
\end{theorem}

Next we improve Theorem \ref{the: sufficient condition B(d)} in the case of groups
with finite virtual cohomological dimension. Notice that for such a group
there is an upper bound on the length $l(H)$ of finite subgroups
$H \subset G$.

\begin{theorem}[Virtual cohomological dimension and $\dim(\underline{E}G$]
\index{Theorem!Virtual cohomological dimension and $\dim(\underline{E}G$}
 \label{vcd le d implies dim le l+d}
Let $G$ be a discrete group which contains a torsionfree subgroup of finite index
and has \emph{virtual cohomological dimension}%
\index{virtual cohomological dimension}
 $\vcd(G) \le d$.%
\indexnotation{vcd(G)}
Let $l \ge 0$ be an  integer
such that the length $l(H)$ of any finite subgroup $H \subset G$
is bounded by $l$.

Then we have $\vcd(G) \le dim(\underline{E}G)$ for any model for $\underline{E}G$ and
there is a $G$-$CW$-model for $\underline{E}G$ such that for any finite subgroup
$H \subset G$
$$\dim(\underline{E}G^H) ~ = ~ \max\{3,d\} + l - l(H)$$
holds.
In particular there exists a model for $\underline{E}G$ of dimension $\max\{3,d\} + l$.
\end{theorem}

\begin{theorem}[Discrete subgroups of Lie groups]
\index{Theorem!Discrete subgroups of Lie groups}
 \label{the: dim(L/K)} 
Let $L$ be a Lie group with finitely many path components.
Then $L$ contains a maximal compact subgroup $K$ which is unique up to conjugation.
Let $G \subseteq L$ be a discrete subgroup of $L$. Then $L/K$ with the
left $G$-action is a model for $\underline{E}G$. 

Suppose additionally that $G$ contains a torsionfree subgroup $\Delta \subseteq G$ of
finite index. Then we have 
$$\vcd(G) \le \dim(L/K)$$
and equality holds if and only if $G\backslash L$ is compact.
\end{theorem}
\begin{proof} We have already mentioned in 
Theorem \ref{the: Discrete subgroups of almost connected Lie groups}
that $L/K$ is a model for $\underline{E}G$. The restriction of
$\underline{E}G$ to $\Delta$ is a $\Delta$-$CW$-model for
$\underline{E}\Delta$ and hence $\Delta\backslash\underline{E}G$ is a
$CW$-model for $B\Delta$. This implies $\vcd(G) := \cd(\Delta) \le
\dim(L/K)$. Obviously $\Delta\backslash L/K$ is a manifold without
boundary. Suppose that $\Delta\backslash L/K$ is compact.
Then $\Delta\backslash L/K$ is a closed manifold and hence
its homology with $\bbZ/2$-coefficients in the top dimension is
non-trivial. This implies $\cd(\Delta) \ge \dim(\Delta\backslash L/K)$
and hence $\vcd(G) = \dim(L/K)$. If $\Delta\backslash L/K$ is not
compact, it contains a $CW$-complex $X \subseteq \Delta\backslash L/K$
of dimension smaller than $\Delta\backslash L/K$ such that
the inclusion of $X$ into $\Delta\backslash L/K$ is a homotopy
equivalence. Hence $X$ is another model for $B\Delta$. This implies
$\cd(\Delta) < \dim(L/K)$ and hence $\vcd(G) < \dim(L/K)$. 
\end{proof}

\begin{remark} \label{rem: strategy to find smaller models} \em 
An often useful strategy to find smaller models for $\EGF{G}{\calf}$
is to look for a $G$-$CW$-subcomplex $X \subseteq \EGF{G}{\calf}$ such
that there exists a $G$-retraction $r \colon \EGF{G}{\calf} \to X$,
i.e.\ a $G$-map $r$ with $r|_X = \id_X$. Then $X$ is automatically
another model for $\EGF{G}{\calf}$. We have seen this already in the
case $SL_2(\bbZ)$, where we found a tree inside 
$\bbH^2 = SL_2(\bbR)/SO(2)$ as explained in 
Subsection \ref{subsec: Special Linear Groups of (2,2)-Matrices}.
This method can be used to construct a model for
$\underline{E}SL_n(\bbZ)$ of dimension $\frac{n(n-1)}{2}$ and to show
that the virtual cohomological dimension of $SL_n(\bbZ)$ is 
$\frac{n(n-1)}{2}$.  Notice that $SL_n(\bbR)/SO(n)$ is also a model 
for $\underline{E}SL_n(\bbZ)$ by 
Theorem \ref{the: Discrete subgroups of almost connected Lie groups}
but has dimension $\frac{n(n+1)}{2} -1$. \em
\end{remark}

\begin{example}[Virtually poly-cyclic groups] \em \label{exa: virtually poly-cyclic} 
Let the group $\Delta$ be  \emph{virtually poly-cyclic},%
\index{group!virtually poly-cyclic}
i.e.\ $\Delta$ contains a subgroup $\Delta'$ of finite index for which there is
a finite sequence $\{1\} =  \Delta_0' \subseteq \Delta_1' \subseteq \ldots \subseteq \Delta_n' =
\Delta'$ of subgroups such that $\Delta_{i-1}'$ is normal in $\Delta_{i}'$ with cyclic
quotient $\Delta_{i}'/\Delta_{i-1}'$ for $i = 1,2, \ldots , n$. 
Denote by $r$ the number of elements $i \in \{1,2, \ldots, n\}$ with 
$\Delta_i'/\Delta_{i-1}' \cong \bbZ$. The number $r$ is called the \emph{Hirsch rank}.%
\index{Hirsch rank}
 The group $\Delta$ contains a torsionfree
subgroup of finite index. We call $\Delta'$ \emph{poly-$\bbZ$}%
\index{group!poly-$\bbZ$}
if $r = n$, i.e.\ all quotients $\Delta_i'/\Delta_{i-1}'$ are infinite cyclic. We want to show:
\begin{enumerate} 
\item $r = \vcd(\Delta)$;

\item $r = \max\{i \mid H_i(\Delta';\bbZ/2) \not= 0\}$ for one (and hence all)
poly-$\bbZ$ subgroup $\Delta' \subset \Delta$ of finite index;

\item There exists a finite $r$-dimensional model for $\underline{E}\Delta$
and for any model $\underline{E}\Delta$ we have $\dim(\underline{E}\Delta) \ge r$.

\end{enumerate}

We use induction over the number $r$. If $r = 0$, then $\Delta$ is finite
and all the claims are obviously true. Next we explain the induction step from
$(r-1)$ to $r \ge 1$. We can choose an extension
$1 \to \Delta_0 \to \Delta \xrightarrow{p} V \to 1$ for some virtually poly-cyclic group $\Delta_0$
with $r(\Delta_0) = r(\Delta) -1 $ and some group $V$ which contains $\bbZ$ as subgroup of
finite index. The induction
hypothesis applies to any group $\Gamma$ which contains $\Delta_0$ as subgroup of finite
index.  Since $V$ maps surjectively to
$\bbZ$ or the infinite dihedral group $D_{\infty}$ with finite kernel and both
$\bbZ$ and $D_{\infty}$ have $1$-dimensional models for their classifying space
for proper group actions, there is a $1$-dimensional model for $\underline{E}V$.
We conclude from Theorem 
\ref{the: dimension of underline{E} and exact sequences}
that there is a $r$-dimensional model for $\underline{E}\Delta$.

The existence of a $r$-dimensional model for $\underline{E}\Delta$
implies $\vcd(\Delta) \le r$. 

For any torsionfree
subgroup $\Delta' \subset \Delta$ of finite index we have
$\max \{i \mid H_i(\Delta';\bbZ/2) \not= 0\} \le \vcd(\Delta)$, 

It is not hard to check by induction over $r$
that we can find a sequence of torsionfree subgroups
$\{1\} \subseteq \Delta_0 \subseteq \Delta_1 \subseteq \ldots \subseteq \Delta_r \subseteq
\Delta$ such that $\Delta_{i-1}$ is normal in $\Delta_i$ with $\Delta_i/\Delta_{i-1} \cong
\bbZ$ for $i \in \{1,2, \ldots , r\}$ and $\Delta_r$ has finite index
in $\Delta$.
We show by induction over $i$ that $H_i(\Delta_i;\bbZ/2) = \bbZ/2$ for  
$i = 0,1, \ldots, r$. The induction beginning $i = 0$ is trivial. The induction step
from $(i-1)$ to $i$ follows from the part of the long exact Wang sequence
\begin{multline*}
H_i(\Delta_{i-1};\bbZ/2) ~ = ~ 0 \to H_i(\Delta_i;\bbZ/2) \to 
H_{i-1}(\Delta_{i-1};\bbZ/2) = \bbZ/2
\\
\xrightarrow{\id - H_{i-1}(f;\bbZ/2) ~ = ~ 0} 
H_{i-1}(\Delta_{i-1};\bbZ/2)
\end{multline*}
which comes from the Hochschild-Serre spectral sequence associated to the 
extension $ 1 \to \Delta_{i-1} \to \Delta_i \to \bbZ \to 1$ for 
$f \colon \Delta_{i-1} \to \Delta_{i-1}$ the automorphism induced by conjugation 
with some preimage in $\Delta_i$ of the generator of $\bbZ$.
This implies 
$$r = \max \{i \mid H_i(\Delta_r;\bbZ/2) \not= 0\} = \cd(\Delta_r) = \vcd(\Delta).$$
Now the claim follows. 

The existence of a $r$-dimensional model for 
$\underline{E}G$  is proved for finitely
generated nilpotent groups with $\vcd(G) \le r$ for $r \not= 2$ in 
\cite{Nucinkis(2003g)}, where also not necessarily finitely generated
nilpotent groups are studied.

The work of Dekimpe-Igodt \cite{Dekimpe-Igodt(1997)} or 
Wilking \cite[Theorem 3]{Wilking(2000b)} implies
that there is a model for $\EGF{\Delta}{\calfin}$ whose underlying space is
$\bbR^r$. 
\em
\end{example}


\subsection{Counterexamples}
\label{subsec: Counterexampples}

The following problem is  stated by Brown \cite[page 32]{Brown(1979)}.
It created a lot of activities and many of the results stated above were motivated
by it.

\begin{problem} \label{prob: Brown's problem}
For which discrete groups $G$, which contain a torsionfree subgroup
of finite index and has virtual cohomological dimension $\le d$, does there exist
a $d$-dimensional $G$-$CW$-model for $\underline{E}G$?
\end{problem}

The following four problems for discrete groups $G$ are stated in the problem lists appearing in 
\cite{Lueck(2000a)} and \cite{Wall(1979b)}.

\begin{problem}
\label{pro: subgroups of finite index and finite type}
Let $H \subseteq G$ be a subgroup of finite index. Suppose
that $\underline{E}H$ has a $H$-$CW$-model of finite type or
a finite $H$-$CW$-model respectively. Does then
 $\underline{E}G$ have a $G$-$CW$-model of finite type
or a finite $G$-$CW$-model respectively?
\end{problem}

\begin{problem} \label{pro: finitely many conjugacy classes I}
If the group $G$ contains a subgroup of finite index
$H$ which has a $H$-$CW$-model of finite type
for $\underline{E}H$, does then $G$ contain only finitely
many conjugacy classes of finite subgroups?
\end{problem}

\begin{problem} \label{pro: NH for finite groups H and FP_infty}
Let $G$ be a group such that $BG$ has a model of finite type.
Is then $BW_GH$ of finite type
for any finite subgroup $H \subset G$?
\end{problem}

\begin{problem} \label{pro: conclusion from Delta and Gamma to pi}
Let $1 \to \Delta \xrightarrow{i} G \xrightarrow{p} \pi \to 1$ be an
exact sequence of groups. Suppose that there is
a $\Delta$-$CW$-model of finite type for $\underline{E}\Delta$
and a $G$-$CW$-model of finite type for $\underline{E}G$.
Is then there a $\pi$-$CW$-model of finite type for $\underline{E}\pi$?
\end{problem}

Leary and Nucinkis \cite{Leary-Nucinkis(2003)}  have constructed many
very interesting examples of discrete groups some of which are listed below. 
Their main technical input is an equivariant version of the constructions
due to Bestvina and Brady \cite{Bestvina-Brady(1997)}.
These examples show
that the answer to the Problems \ref{prob: Brown's problem},
\ref{pro: subgroups of finite index and finite type}, 
\ref{pro: finitely many conjugacy classes I}, 
\ref{pro: NH for finite groups H and FP_infty} 
and \ref{pro: conclusion from Delta and Gamma to pi} above
is \emph{not} positive in general.
A group $G$ is \emph{of type $V\!F$}%
\index{type ${V\hspace{-1mm}F}$}
\indexnotation{VF}
if it contains a subgroup $H \subseteq G$ of finite index
for which there is a finite model for $BH$. 

\begin{enumerate} 

\item For any positive integer $d$ there exist a group $G$ of type $V\!F$  which has 
virtually cohomological dimension $\le 3d$, but for which any model for
$\underline{E}G$ has dimension $\ge 4d$;

\item There exists a group $G$ with a finite cyclic subgroup $H \subseteq G$ 
such that $G$ is of type $V\!F$ but the centralizer $C_GH$ of $H$ in $G$ is
not of type $F\!P_{\infty}$;

\item There exists a group $G$ of type $V\!F$ which contains infinitely many conjugacy
classes of finite subgroups;

\item There exists an extension $1 \to \Delta \to G \to \pi \to 1$
such that $\underline{E}\Delta$ and $\underline{E}G$ have finite
$G$-$CW$-models
but there is no $G$-$CW$-model for $\underline{E}\pi$ of finite type.

\end{enumerate}


\typeout{--------------------   Section 6: On the space G backslash underline E G --------------------------}

\section{The Orbit Space of \underline{E}G}
\label{sec: On the space G backslash underline E G}

We will see that in many computations of the group (co-)homology,
of the algebraic $K$- and $L$-theory of the group ring or the topological $K$-theory 
of the reduced $C^*$-algebra of a discrete group $G$ a key problem is to determine the homotopy type
of the quotient space $G\backslash \underline{E}G$ of $\underline{E}G$.
 The following result shows that this is a difficult problem in general
and can only be solved in special cases. 
It was proved by Leary and Nucinkis \cite{Leary-Nucinkis(2001a)} based on ideas due to
Baumslag-Dyer-Heller \cite{Baumslag-Dyer-Heller(1980)} and Kan and Thurston \cite{Kan-Thurston(1976)}.

\begin{theorem}[The homotopy type of $G\backslash \underline{E}G$]
\index{Theorem!The homotopy type of $G\backslash \underline{E}G$}
\label{the: Gbackslash underline E G can be anything}
Let $X$ be a connected $CW$-complex. Then there exists a group $G$ such that
$G\backslash \underline{E}G$ is homotopy equivalent to $X$.
\end{theorem}

There are some cases, where the quotient $G\backslash \underline{E}G$ has been
determined explicitly using geometric input. We mention a few examples.

\begin{enumerate}

\item Let $G$ be a planar group (sometimes also called NEC) group, i.e.
a discontinuous group of isometries of the two-sphere $S^2$, the Euclidean plane $\bbR^2$,
or the hyperbolic plane $\bbH^2$. Examples are Fuchsian groups and two-dimensional
crystallographic groups. If $G$ acts on $\bbR^2$ or $\bbH^2$ and the action is
cocompact, then $\bbR^2$ or $\bbH^2$ is a model for $\underline{E}G$ and the
quotient space $G\backslash\underline{E}G$ is a compact $2$-dimensional surface. 
The number of boundary components, its genus and the answer to the question, whether
$G\backslash\underline{E}G$ is orientable, can be read off from an explicit presentation of $G$. 
A summary of these details can be found in \cite[Section 4]{Lueck-Stamm(2000)},
where further references to papers containing proofs of the stated facts are given;

\item Let $G = \langle (q_i)_{i \in I} \mid r \rangle$  be a one-relator group.
Let $F$ be the free group on the letters
$\{q_i \mid i \in I\}$. Then $r$ is an element in
$F$. There exists an element $s \in F$ and an integer $m \ge 1$
such that $r = s^m$, the cyclic  subgroup $C$
generated by the class $\overline{s} \in G$ represented 
by $s$ has order $m$, any finite subgroup of $G$ is subconjugated to $C$ and 
for any $g \in G$ the implication
$g^{-1}Cg \cap C \not= \{1\} \Rightarrow g \in C$ holds
(see \cite[Propositions 5.17, 5.18 and 5.19 in II.5 on pages 107 and 108]{Lyndon-Schupp(1977)}).

In the sequel we use the two-dimensional model for $\underline{E}G$ described in
Subsection \ref{subsec: One-Relator Groups}.
Let us compute the integral homology of $BG$ and $G\backslash \underline{E}G$.
Since $G\backslash \underline{E}G$  has precisely one
$2$-cell and is two-dimensional, $H_2(G\backslash \underline{E}G)$ is either trivial or
infinite cyclic and $H_k(G\backslash\underline{E}G) = 0$ for $k \ge 3$.
We obtain the short exact sequence
\begin{multline*}
0 \to H_2(BG) \xrightarrow{H_2(q)}{H_2(G\backslash \underline{E}G)}
\xrightarrow{\partial_2} H_1(BC) \xrightarrow{H_1(Bi)} H_1(BG) 
\\ \xrightarrow{H_1(q)} H_1(G\backslash\underline{E}G) \to 0
\end{multline*}
and for $k \ge 3$ isomorphisms
$$H_k(Bi) \colon H_k(BC) \xrightarrow{\cong} H_k(BG)$$
from the pushout coming from \eqref{pushout for underline E G}
$$\comsquare{BC}{i}{BG}{}{}{\pt}{}{G\backslash\underline{E}G}$$
Hence $H_2(G\backslash \underline{E}G) = 0$ and the sequence
$$0 \to H_1(BC) \xrightarrow{H_1(Bi)} H_1(BG) 
\xrightarrow{H_1(q)} H_1(G\backslash\underline{E}G) \to 0$$ is exact,
provided that $H_2(BG) = 0$. Suppose that $H_2(BG) \not = 0$. 
Hopf's Theorem says that
$H_2(BG) \cong R \cap [F,F]/[F,R]$ if $R$ is the subgroup of $G$ normally generated by 
$r \in F$ (see \cite[Theorem 5.3 in II.5 on page 42]{Brown(1982)}). 
For every element in $R \cap [F,F]/[F,R]$ there exists $n \in \bbZ$ such that
$r^n$ belongs to $[F,F]$ and the element is represented by $r^n$.
Hence there is $n \ge 1$ such that $r^n$ does belong to $[F,F]$. 
Since $F/[F,F]$ is torsionfree, also $s$ and $r$ belong to $[F,F]$.
We conclude that both $H_2(BG)$ and
$H_2(G\backslash\underline{E}G)$ are infinite cyclic groups, $H_1(BC) \to H_1(BG)$ is trivial
and ${H_1(q) \colon H_1(BG) 
\xrightarrow{\cong}} H_1(G\backslash\underline{E}G)$ is bijective. 
We also see that $H_2(BG) = 0$ if and only if $r$ does not belong to
$[F,F]$.
 
\item Let $\Hei$ be the three-dimensional discrete
Heisenberg group which is the subgroup of $GL_3(\bbZ)$ consisting of
upper triangular matrices with $1$ on the diagonals. Consider the $\bbZ/4$-action 
given by
$$\left(\begin{array}{ccc} 1& x & y \\ 0 & 1 & z \\ 0 & 0 & 1\end{array}\right)
~ \mapsto ~ 
\left(\begin{array}{ccc} 1& -z & y-xz \\ 0 & 1 & x \\ 0 & 0 & 1\end{array}\right).$$
Then a key result in \cite{Lueck(2004g)} is that
$G\backslash\underline{E}G$ is homeomorphic to $S^3$ for $G = \Hei \rtimes \bbZ/4$;

\item A key result in \cite[Corollary on page 8]{Soule(1978)}
implies that for $G = SL_3(\bbZ)$ the quotient space $G\backslash\underline{E}G$ is contractible.

\end{enumerate}


\typeout{--------   Section 7: Relevance and Applications of Classifying Spaces for Families--------}

\section{Relevance and Applications of Classifying Spa\-ces for Families}
\label{sec: Relevance and Applications of Classifying Spaces for Families}

In this section we discuss some theoretical aspects which involve and rely on the notion
of a classifying space for a family of subgroups.


\subsection{Baum-Connes Conjecture}
\label{subsec: Baum-Connes Conjecture}

Let $G$ be a locally compact second countable Hausdorff group. 
Using the equivariant $KK$-theory due to Kasparov one can assign to a 
$\calcom$-numerable $G$-space $X$ its equivariant $K$-theory $K^G_n(X)$. Let
$C_r^*(G)$%
\indexnotation{C_r^*(G)}
be the reduced group $C^*$-algebra associated to $G$.
The goal of the Baum-Connes Conjecture is to compute the topological
$K$-theory $K_p(C_r^*(G))$. The following formulation is taken from
\cite[Conjecture 3.15]{Baum-Connes-Higson(1994)}.

\begin{conjecture}[Baum-Connes Conjecture] 
\label{con: Baum-Connes Conjecture} \index{Conjecture!Baum-Connes
  Conjecture} The assembly map defined by taking the equivariant index
$$\asmb \colon K_n^G(\underline{J}G) \xrightarrow{\cong}
K_n(C^*_r(G))$$
is bijective for all $n \in \bbZ$.
\end{conjecture}

More information about this conjecture and its relation and application to other
conjectures and problems can be found for instance in 
\cite{Baum-Connes-Higson(1994)}, \cite{Higson(1998a)}, \cite{Lueck-Reich(2004g)},
\cite{Mislin-Valette(2003)}, \cite{Valette(2002)}.


\subsection{Farrell-Jones Conjecture}
\label{subsec: Farrell-Jones Conjecture}

Let $G$ be a discrete group. Let $R$ be a associative ring with unit.
One can construct a $G$-homology theory $\calh^G_*(X;\bfK)$ graded over the integers
and defined for $G$-$CW$-complexes $X$ such that for any subgroup $H \subseteq G$ the
abelian group $\calh^G_n(G/H;\bfK)$ is isomorphic to the algebraic $K$-groups $K_n(RH)$ for $n \in
\bbZ$. If $R$ comes with an involution of rings, one can also 
construct a $G$-homology theory $\calh^G_*(X;\bfL^{\langle -\infty\rangle})$ graded over the integers
and defined for $G$-$CW$-complexes $X$ such that for any subgroup $H \subseteq G$ the
abelian group $\calh^G_n(G/H;\bfL^{\langle -\infty\rangle})$ is isomorphic to the algebraic $L$-groups
$L_n^{- \infty}(RH)$ for $n \in \bbZ$. Let $\calvcyc$ be the family of
virtually cyclic subgroups of $G$.
The goal of the Farrell-Jones Conjecture is to compute the 
algebraic $K$-groups $K_n(RH)$ and the algebraic  $L$-groups. The following formulation is
equivalent to the original one appearing in \cite[1.6 on page 257]{Farrell-Jones(1993a)}.

\begin{conjecture}[Farrell-Jones Conjecture] 
\label{con: Farrell-Jones Conjecture} 
\index{Conjecture!Farrell-Jones Conjecture} 
The assembly maps induced by the projection $\EGF{G}{\calvcyc} \to G/G$
\begin{eqnarray}
\asmb \colon \calh^G_n(\EGF{G}{\calvcyc},\bfK) & \to &
\calh^G_n(G/G,\bfK) = K_n(RG);
\\
\asmb \colon \calh^G_n(\EGF{G}{\calvcyc},\bfL^{-\infty}) & \to &
\calh^G_n(G/G,\bfL^{-\infty}) = L^{-\infty}_n(RG),
\end{eqnarray}
are bijective for all $n \in \bbZ$.
\end{conjecture}

More information about this conjecture and its relation and application to other
conjectures and problems can be found for instance in 
\cite{Farrell-Jones(1993a)} and \cite{Lueck-Reich(2004g)}.

We mention that for a discrete group $G$ one can formulate the Baum-Connes Conjecture in a
similar fashion. Namely, 
one can also  construct a $G$-homology theory $\calh^G_*(X;\bfK^{\topo})$ graded over the integers
and defined for $G$-$CW$-complexes $X$ such that for any subgroup $H \subseteq G$ the
abelian group $\calh^G_n(G/H;\bfK^{\topo})$ is isomorphic to the topological $K$-groups
$K_n(C_r^*(H))$ for $n \in \bbZ$ and the assembly map appearing in the Baum-Connes
Conjecture can be identified with the map induced by the projection
$\underline{J}G = \underline{E}G \to G/G$ (see \cite{Davis-Lueck(1998)},
\cite{Hambleton-Pedersen(2003)}). If the ring $R$ is regular and contains $\bbQ$ as
subring, then one can replace in the Farrell-Jones Conjecture
\ref{con: Farrell-Jones Conjecture} 
$\EGF{G}{\calvcyc}$ by $\underline{E}G$ but this is not possible for arbitrary rings
such as $R = \bbZ$. This comes from the appearance of Nil-terms in the Bass-Heller-Swan
decomposition which do not occur in the context of the topological $K$-theory of reduced
$C^*$-algebras.

Both the Baum-Connes Conjecture \ref{con: Baum-Connes Conjecture} and the Farrell-Jones 
Conjecture \ref{con: Farrell-Jones Conjecture} allow
to reduce the computation of certain $K$-and $L$-groups of the group ring or the
reduced $C^*$-algebra of a group $G$ to the computation of
certain $G$-homology theories applied to $\underline{J}G$, $\underline{E}G$ or 
$\EGF{G}{\calvcyc}$. Hence it is important to find good models for these spaces
or to make predictions about their dimension or whether they are finite or of finite type.


\subsection{Completion Theorem}
\label{subsec: Completion Theorem}

Let $G$ be a discrete group. For a proper finite $G$-$CW$-complex let $K_G^*(X)$ be its
equivariant $K$-theory defined in terms of equivariant finite dimensional complex vector 
bundles over $X$ (see \cite[Theorem 3.2]{Lueck-Oliver(2001a)}).
It is a $G$-cohomology theory with a multiplicative structure.
Assume that $\underline{E}G$ has a finite $G$-$CW$-model. 
Let $I \subseteq K_G^0(\underline{E}G)$ be the augmentation ideal, i.e.\ the kernel of the
map $K^0(\underline{E}G)\to \bbZ$ sending the class of an equivariant complex vector bundle
to its complex dimension. Let $K_G^*(\underline{E}G)\widehat{_I}$ be the $I$-adic
completion of $K_G^*(\underline{E}G)$ and let $K^*(BG)$ be the topological $K$-theory of $BG$.

\begin{theorem}[Completion Theorem for discrete groups]
\label{the: Completion Theorem for discrete groups}
\index{Theorem!Completion Theorem for discrete groups}
Let $G$ be a discrete group such that there exists a finite model for $\underline{E}G$.
Then there is a canonical isomorphism
$$K^*(BG) \xrightarrow{\cong} K_G^*(\underline{E}G)\widehat{_I}.$$
\end{theorem}

This result is proved in \cite[Theorem 4.4]{Lueck-Oliver(2001a)}, where a more general
statement is given provided that there is a finite dimensional model for $\underline{E}G$
and an upper bound on the orders of finite subgroups of $G$. 
In the case where $G$ is finite, Theorem \ref{the: Completion Theorem for discrete groups}
reduces to the Completion Theorem due to Atiyah and Segal \cite{Atiyah(1961)},
\cite{Atiyah-Segal(1969)}. A Cocompletion Theorem for equivariant 
$K$-homology will appear in \cite{Joachim-Lueck(2004)}.


\subsection{Classifying Spaces for Equivariant Bundles}
\label{subsec: Classifying Spaces for Equivariant Bundles}

In \cite{Lueck-Oliver(2001b)} the equivariant $K$-theory for 
finite proper $G$-$CW$-complexes appearing in Subsection
\ref{subsec: Completion Theorem} above is extended 
to arbitrary proper $G$-$CW$-complexes (including the multiplicative structure)
using $\Gamma$-spaces in the sense of Segal and involving classifying spaces
for equivariant vector bundles. These classifying spaces for equivariant vector bundles
are again classifying spaces of certain Lie groups and certain families
(see \cite[Section 8 and 9 in Chapter I]{Dieck(1987)},
\cite[Lemma 2.4]{Lueck-Oliver(2001a)}).


\subsection{Equivariant Homology and Cohomology}
\label{subsec: Equivariant Homology and Cohomology}

Classifying spaces for families play a role in computations of equivariant homology
and cohomology for compact Lie groups such as equivariant bordism 
as explained in \cite[Chapter 7]{Dieck(1979)}, \cite[Chapter
III]{Dieck(1987)}. Rational computations of equivariant (co-)-homology groups
are possible in general using Chern characters
for discrete groups and proper $G$-$CW$-complexes (see
\cite{Lueck(2002b)}, \cite {Lueck(2002d)}, \cite{Lueck(2004k)}).


\typeout{---------   Section 8: Computations using on classifying spaces of families----------}

\section{Computations using Classifying Spaces for Families}
\label{sec: Computations using Classifying Spaces for Families}

In this section we discuss some computations  which involve and rely on the notion
of a classifying space for a family of subgroups. These computations are possible
since one understands in the cases of interest the geometry of $\underline{E}G$ and
$G\backslash \underline{E}G$. We focus on the case described in Subsection 
\ref{subsec: Groups with Appropriate Maximal Finite Subgroups}, namely of a discrete group 
$G$ satisfying the conditions (M) and (NM). Let $s \colon EG \to \underline{E}G$ be the up
to $G$-homotopy unique $G$-map. Denote by $j_i \colon M_i \to G$ the inclusion.


\subsection{Group Homology}
\label{subsec: Group Homology}

We begin with the group homology
$H_n(BG)$ (with integer coefficients). 
Let $\widetilde{H}_p(X)$ be the reduced homology, i.e.\ the kernel of the map
$H_n(X) \to H_n(\pt)$ induced by the projection $X \to \pt$.
The Mayer-Vietoris sequence applied to the pushout, which is obtained from the $G$-pushout
\eqref{pushout for underline E G} by dividing out the $G$-action, yields
the long Mayer-Vietoris sequence
\begin{multline} \ldots \to H_{p+1}(G\backslash \underline{E}G)) \xrightarrow{\partial_{p+1}}
\bigoplus_{i \in I} \widetilde{H}_p(BM_i) \xrightarrow{\bigoplus_{i \in I} H_p(Bj_i)}
H_p(BG) 
\\
\xrightarrow{H_p(G\backslash s)} H_p(G\backslash\underline{E}G) \xrightarrow{\partial_p} \ldots 
\label{Mayer-Vietoris for H_p(BG)}
\end{multline}
In particular we obtain an isomorphism for $p \ge \dim(\underline{E}G) + 2$
\begin{eqnarray}
\bigoplus_{i \in I} H_p(Bj_i) \colon \bigoplus_{i \in I} \widetilde{H}_p(BM_i) 
& \xrightarrow{\cong} & H_p(BG).
\end{eqnarray}
This example and the forthcoming ones show 
why it is important to get upper bounds on  the dimension of $\underline{E}G$ and
to understand the quotient space $G\backslash \underline{E}G$. For Fuchsian groups and 
for one-relator groups we have $\dim(G\backslash \underline{E}G) \le 2$ and it is easy to compute
the homology of $G\backslash \underline{E}G$ in this case as explained in Section 
\ref{sec: On the space G backslash underline E G}. 


\subsection{Topological $K$-Theory of Group $C^*$-Algebras}
\label{subsec: Topological K-theory of Group C^*-Algebras}

Analogously one can compute the source of the assembly map appearing in the
Baum-Connes Conjecture \ref{con: Baum-Connes Conjecture}. Namely, the Mayer-Vietoris
sequence associated to the $G$-pushout
\eqref{pushout for underline E G} and the one associated to its quotient under the
$G$-action look like
\begin{multline} \ldots \to K_{p+1}^G(\underline{E}G) \to
\bigoplus_{i \in I} K_p^G(G\times_{M_i} EM_i) 
\\
\to 
\left(\bigoplus_{i \in I} K_p^G(G/M_i)\right) \bigoplus K_p^G(EG) 
\to
 K_p^G(\underline{E}G) \to \ldots
\label{Mayer-Vietoris for K_p^G of G-pushout}
\end{multline}
and
\begin{multline} \ldots \to K_{p+1}(G\backslash \underline{E}G) \to
\bigoplus_{i \in I} K_p(BM_i) 
\\
\to 
\left(\bigoplus_{i \in I} K_p(\pt)\right) \bigoplus K_p(BG) 
\to
\bigoplus_{i \in I} K_p(G\backslash \underline{E}G) \to \ldots
\label{Mayer-Vietoris for K_p of the quotient pushout}
\end{multline}
Notice that for a free $G$-$CW$-complex $X$ there is a canonical
isomorphisms $K_p^G(X) \cong K_p(G\backslash X)$.
We can splice these sequences together and obtain the long exact sequence
\begin{multline}
\ldots \to K_{p+1}(G\backslash \underline{E}G) \to \bigoplus_{i \in I} K_p^G(G/M_i)
\to \bigoplus_{i \in I} K_p(\pt) \bigoplus K_p^G(\underline{E}G) 
\\ \to K_p(G\backslash\underline{E}G) \to \ldots
\label{spliced Mayer-Vietoris sequence}
\end{multline}
There are identification of $K_0^G(G/M_i)$ with the complex representation ring 
$R_{\bbC}(M_i)$ of the finite group $M_i$ and of $K_0(\pt)$ with $\bbZ$. Under 
these identification the map $K_0^G(G/M_i) \to K_0(\pt)$ becomes the split surjective
map $\epsilon \colon R_{\bbC}(M_i) \to \bbZ$ which sends the class of a complex
$M_i$-representation $V$ to the complex dimension of $\bbC \otimes_{\bbC[M_i]} V$. 
The kernel of this map
is denoted by $\widetilde{R}_{\bbC}(M_i)$. The groups
$K_1^G(G/M_i)$ and  $K_1(\pt)$ vanish.
The abelian group $R_{\bbC}(M_i)$ and hence
also $\widetilde{R}_{\bbC}(M_i)$ are finitely generated free abelian groups.
If $\bbZ \subseteq \Lambda \subseteq \bbQ$ is ring such that the order
of any finite subgroup of $G$ is invertible in $\Lambda$, then the map
$$\Lambda \otimes_{\bbZ} K_p^G(s) \colon \Lambda \otimes_{\bbZ} K_p^G(EG) 
\to \Lambda \otimes_{\bbZ} K_p(G\backslash \underline{E}G)$$ 
is  an isomorphism for all $p \in \bbZ$ \cite[Lemma 2.8 (a)]{Lueck-Stamm(2000)}. 
Hence the long exact 
sequence \eqref{spliced Mayer-Vietoris sequence} splits after applying 
$ \Lambda \otimes_{\bbZ} -$.
We conclude from the long exact 
sequence \eqref{spliced Mayer-Vietoris sequence} since the representation ring of a finite
group is torsionfree

\begin{theorem} \label{the: computation of K_p(C_r^*(G)) for special G}
Let $G$ be a discrete group which satisfies the conditions (M) and (NM) appearing in 
Subsection \ref{subsec: Groups with Appropriate Maximal Finite Subgroups}.
Suppose that the Baum-Connes Conjecture \ref{con: Baum-Connes Conjecture} is true for
$G$. Let $\{(M_i) \mid i \in I\}$ be the set of conjugacy classes of maximal finite
subgroups of $G$. Then there is an isomorphism
$$K_1(C^*_r(G)) \xrightarrow{\cong} K_1(G\backslash\underline{E}G)$$
and a short exact sequence
$$
0 \to \bigoplus_{i \in I} \widetilde{R}_{\bbC}(M_i) \to 
K_0(C^*_r(G)) \to K_0(G\backslash\underline{E}G) \to 0,
$$
which splits if we invert the orders of all finite subgroups of $G$.
\end{theorem}


\subsection{Algebraic $K$-and $L$-Theory of Group Rings}
\label{subsec: Topological K-Theory of Group Rings}

Suppose that $G$ satisfies the Farrell-Jones 
Conjecture \ref{con: Farrell-Jones Conjecture}.  Then the computation
of the relevant groups $K_n(RG)$ or $L_n^{\langle -  \infty\rangle}(RG)$ 
respectively is equivalent to the computation of $\calh^G_n(\EGF{G}{\calvcyc},\bfK)$
or  $\calh^G_n(\EGF{G}{\calvcyc},\bfL^{-\infty})$ respectively. The following result
is due to Bartels \cite{Bartels(2003)}. Recall that $\underline{E}G$ is the same as
$\EGF{G}{\calfin}$.

\begin{theorem} \label{the: splitting of calfin in K- and L-theory}
\begin{enumerate}
\item \label{the: splitting of calfin in K- and L-theory: K}
For every group $G$, every ring $R$ and every $n \in \bbZ$ the up to $G$-homotopy unique
$G$-map
$f \colon \EGF{G}{\calfin} \to \EGF{G}{\calvcyc}$ induces a split injection
$$H_n^G(f;\bfK_R ) \colon  H_n^G(\EGF{G}{\calfin};\bfK_R) \to
H_n^G(\EGF{G}{\calvcyc};\bfK_R); 
$$

\item \label{the: splitting of calfin in K- and L-theory: L}
Suppose $R$ is such that 
$K_{-i} ( RV )= 0$ for all virtually cyclic subgroups $V$ of $G$ and for sufficiently large $i$
(for example $R= \bbZ$ will do). Then we get a split injection
$$
H_n^G (f;\bfL_R^{\langle - \infty \rangle}) \colon
H_n^G(\EGF{G}{\calfin};\bfL_R^{\langle - \infty \rangle}) 
\to
H_n^G(\EGF{G}{\calvcyc};\bfL_R^{\langle - \infty \rangle}).
$$
\end{enumerate}
\end{theorem}
It remains to compute 
$H_n^G(\EGF{G}{\calfin};\bfK)$
and 
$H_n^G(\EGF{G}{\calvcyc},\EGF{G}{\calfin};\bfK)$,
if we arrange $f$ to be a $G$-cofibration and think of
$\EGF{G}{\calfin}$ as a $G$-$CW$-subcomplex of $\EGF{G}{\calvcyc}$.
Namely, we get from the Farrell-Jones 
Conjecture~\ref{con: Farrell-Jones Conjecture} and 
Theorem~\ref{the: splitting of calfin in K- and L-theory} an isomorphism
$$H_n^G(\EGF{G}{\calfin};\bfK) \bigoplus 
H_n^G(\EGF{G}{\calvcyc},\EGF{G}{\calfin};\bfK) \xrightarrow{\cong}
K_n(RG).$$ 
The analogous statement holds for $\bfL_R^{\langle - \infty \rangle})$,
provided $R$ satisfies the conditions appearing in 
Theorem \ref{the: splitting of calfin in K- and L-theory}~
\ref{the: splitting of calfin in K- and L-theory: L}.

Analogously to Theorem \ref{the: computation of K_p(C_r^*(G)) for special G} 
one obtains
\begin{theorem} \label{the: computation of K_p(RG) and L_p(RG)for special G}
Let $G$ be a discrete group which satisfies the conditions (M) and (NM) appearing in 
Subsection \ref{subsec: Groups with Appropriate Maximal Finite Subgroups}.
Let $\{(M_i) \mid i \in I\}$ be the set of conjugacy classes of maximal finite
subgroups of $G$. Then
\begin{enumerate}

\item 
There is a long exact sequence
\begin{multline*}
\ldots \to H_{p+1}(G\backslash\EGF{G}{\calfin};\bfK(R)) \to \bigoplus_{i \in I} K_p(R[M_i])
\\\to \bigoplus_{i \in I} K_p(R) \bigoplus H_p^G(\EGF{G}{\calfin};\bfK_{\bfR}) 
\to 
\bigoplus H_p(G\backslash\EGF{G}{\calfin};\bfK(R))  \to \ldots
\end{multline*}
and analogously for $\bfL^{\langle -\infty \rangle}_{R}$.

\item
For $R = \bbZ$ there are isomorphisms
$$\bigoplus_{i \in I} \Wh_n(M_i) \bigoplus 
H_n^G(\EGF{G}{\calfin},\EGF{G}{\calvcyc};\bfK_{\bbZ}) \xrightarrow{\cong} \Wh_n(G).$$

\end{enumerate}

\end{theorem}

\begin{remark} \label{rem: Heisenberg} 
\em These results about groups satisfying 
conditions (M) and (NM) are extended in \cite{Lueck(2004g)} to groups
which map surjectively to groups satisfying conditions (M) and (NM) 
with special focus on the semi-direct product of the discrete three-dimensional
Heisenberg group with $\bbZ/4$.  \em
\end{remark}

\begin{remark} \label{rem: Soule's computation} \em  
In \cite{Soule(1978)} a special model
for $\underline{E}SL_3(\bbZ)$ is presented which allows to compute the
integral group homology. Information about the algebraic $K$-theory of
$SL_3(\bbZ)$ can be found in \cite[Chapter 7]{Stamm(1999)}, \cite{Upadhyay(1996)}.
\em \end{remark}

The analysis of the other term $H_n^G(\EGF{G}{\calvcyc},\EGF{G}{\calfin};\bfK)$
simplifies considerably under certain assumptions on $G$.

\begin{theorem}[On the structure of $\EGF{G}{\calvcyc}$]
\index{Theorem!On the structure of $\EGF{G}{\calvcyc}$}
\label{the: hyperbolic groups and fin to vcyc}
Suppose that $G$ satisfies the following conditions:

\begin{itemize} 

\item Every infinite cyclic subgroup $C \subseteq G$ has finite index in its centralizer $C_GC$;

\item There is an upper bound on the orders of finite subgroups.

\end{itemize}

(Each word-hyperbolic group satisfies these two conditions.) Then

\begin{enumerate}

\item \label{the: hyperbolic groups and fin to vcyc: maximal}
For an infinite  virtually cyclic subgroup $V \subseteq G$ define
$$V_{\max}~=~ \bigcup~\{N_GC \mid C \subset V \text{ infinite cyclic normal}\}.$$
Then 

\begin{enumerate}

\item \label{the: hyperbolic groups and fin to vcyc: maximal: V_max virt. cyc.}
$V_{\max}$ is an infinite virtually cyclic subgroup of $G$ and contains $V$;

\item \label{the: hyperbolic groups and fin to vcyc: maximal: inclusions}
If $V \subseteq W \subseteq G$ are infinite virtually cyclic subgroups of $G$, then
$V_{\max} = W_{\max}$;

\item \label{the: hyperbolic groups and fin to vcyc: maximal: V_max = N_GV_max}
Each infinite virtually cyclic subgroup $V$ is contained in a unique maximal infinite virtually
cyclic subgroup, namely $V_{\max}$, and $N_GV_{\max} = V_{\max}$;

\end{enumerate}

\item \label{the: hyperbolic groups and fin to vcyc: maximal: model for calvcyc}
Let $\{V_i \mid i \in I \}$ be a complete system of representatives of
conjugacy classes of maximal infinite virtually cyclic subgroups. Then there exists a 
$G$-pushout
$$
\comsquare{\coprod_{i \in I} G \times_{V_i}\EGF{V_i}{\calfin}}{}{\EGF{G}{\calfin}}
{\pr}{}
{\coprod_{i \in I} G/V_i}{}{\EGF{G}{\calvcyc}}
$$
whose upper horizontal arrow is an inclusion of $G$-$CW$-complexes.

\item \label{the: hyperbolic groups and fin to vcyc: maximal: H_n^g(vcyc,fin)}
There are natural isomorphisms
\begin{eqnarray*}
\bigoplus_{i \in I} H_n^{V_i}(\EGF{V_i}{\calvcyc},\EGF{V_i}{\calfin};\bfK_R) 
 & \xrightarrow{\cong} &H_n^G(\EGF{G}{\calvcyc},\EGF{G}{\calfin};\bfK_R)
\\
\bigoplus_{i \in I} H_n^{V_i}(\EGF{V_i}{\calvcyc},\EGF{V_i}{\calfin};\bfL^{\langle -\infty \rangle}_R) 
& \xrightarrow{\cong} &
H_n^G(\EGF{G}{\calvcyc},\EGF{G}{\calfin};\bfL^{\langle -\infty \rangle}_R).
\end{eqnarray*}
\end{enumerate}
\end{theorem}
\begin{proof} Each word-hyperbolic group $G$ satisfies these two conditions by
\cite[Theorem 3.2 in III.$\Gamma$.3 on page 459 and Corollary 3.10 in 
III.$\Gamma$.3 on page 462]{Bridson-Haefliger(1999)}.
\\[2mm]
\ref{the: hyperbolic groups and fin to vcyc: maximal}
Let $V$ be an infinite virtually cyclic subgroup $V \subseteq G$.
Fix a normal infinite cyclic subgroup $C \subseteq V$. Let $b$ be a common multiple
of the orders of finite subgroups of $G$. Put $d :=b\cdot b!$. 
Let $e$ be the index of the infinite cyclic
group $dC = \{d \cdot x \mid x \in C\}$ in its centralizer $C_GdC$. Let $D \subset dC$ be any non-trivial
subgroup. Obviously $dC \subseteq C_GD$. We want to show
\begin{eqnarray}
[C_GD : dC] & \le & b \cdot e^2. \label{estimate on C_GD : dC}
\end{eqnarray}
Since $D$ is central in $C_GD$ and $C_GD$ is virtually cyclic and hence $|C_GD/D| <
\infty$, the spectral sequence associated to the extension $1 \to D \to C_GD \to C_GD/D
\to 1$ implies that the map $D = H_1(D) \to H_1(C_GD)$ is injective and has finite
cokernel. In particular the quotient of $H_1(C_GD)$ by its torsion subgroup 
$H_1(C_GD)/\tors$  is an infinite cyclic group. Let $p_{C_GD} \colon C_GD \to H_1(C_GD)/\tors$
be the canonical epimorphism. Its kernel is a finite normal subgroup.
The following diagram commutes and has exact rows 
$$
\begin{CD}
1 @>>> \ker(p_C) @>>> C_GC @>p_C >>H_1(C_GC)/\tors @>>> 1
\\
& &@VVV @VVV @VVV
\\
1 @>>> \ker(p_D) @>>> C_GD @> p_D >>H_1(C_GD)/\tors @>>> 1
\end{CD}
$$
where the vertical maps are induced by the inclusions
$C_GC \subseteq C_GD$. All vertical maps are injections with finite cokernel.
Fix elements $z_C \in C_GC$ and $z_D \in C_GD$ such that $p_C(z_C)$ 
and $p_D(z_D)$ are generators. Choose $l \in \bbZ$ such that $p_C(z_C)$ is send to 
$l \cdot p_D(z_D)$. Then there is $k \in \ker(p_D)$ with
$z_C = k \cdot z_D^l$.  The order of   $\ker(p_D)$ divides $b$ by assumption.
If $\phi \colon \ker(p_D) \to
\ker(p_D)$ is any automorphism, then $\phi^{b!} = \id$. This implies for any element
$k \in \ker(p_D)$ that
$$\prod_{i=0}^{d-1} \phi^i(k) = \left(\prod_{i=0}^{b!-1}\phi^i(k)\right)^b = 1.$$
Hence we get in $C_GD$ if $\phi$ is conjugation with $z_D^l$
$$z_C^d = (k \cdot z_D^l)^d = \prod_{i=0}^{d-1} \phi^i(k) \cdot z_D^{dl} = z_D^{dl}.$$
Obviously $z_D\in C_GdC$ since $z_C^d = z_D^{dl}$ generates $dC$. Hence $z_D^e$ lies in $dC$ and we get
$z_D^e = z_C^{df}$ for some integer $f$. This implies 
$z_D^{e} = z_D^{ldf}$ and hence that $l$ divides $e$. We conclude that
the cokernel of the map $H_1(C_GC)/\tors \to H_1(C_GD)/\tors$  is bounded by $e$. 
Hence the index $[C_GD : C_GC]$ is bounded by $b \cdot e$ since the order of
$ker(p_D)$ divides $b$. Since $dC \subseteq C_GC \subseteq C_GdC \subseteq C_GD$ holds,
equation~\ref{estimate on C_GD : dC} follows. 

Next we show that there is a normal infinite cyclic subgroup $C_0 \subseteq V$ such that
$V_{\max} = N_GC_0$ holds. If $C'$ and $C''$ are infinite cyclic normal subgroups of $V$, then
both $C_GC'$ and $C_GC''$ are contained in $C_G(C' \cap C'')$ and $C' \cap C''$ is
again an infinite cyclic normal subgroup. 
Hence there is a sequence of normal infinite
cyclic subgroups of $V$
$$dC \supseteq C_1 \supseteq C_2 \supseteq C_3 \supseteq \ldots$$
which yields a  sequence
$C_GdC \subseteq C_GC_1 \subseteq C_GC_2 \subset \ldots$ satisfying
$$\bigcup ~ \{C_GC_n \mid n \ge 1\} ~ = ~  
\bigcup ~ \{C_GC \mid C \subset V \text{ infinite cyclic normal}\}.$$
Because of \ref{estimate on C_GD : dC}  there is an upper bound on
$[C_GC_n : C_GdC]$ which is independent of $n$. Hence there is an index $n_0$ with
$$C_GC_{n_0} ~ = ~ \bigcup~\{C_GC \mid C \subset V \text{ infinite cyclic  normal}\}.$$
For any infinite cyclic subgroup
$C \subseteq G$ the index of $C_GC$ in $N_GC$ is $1$ or $2$. Hence there
is an index $n_1$ with 
$$N_GC_{n_1} ~ =  ~ \bigcup ~ \{N_GC \mid C \subset V \text{ infinite cyclic  normal}\}.$$
Thus we have shown the existence of a normal infinite cyclic subgroup $C \subseteq V$
with $V_{\max} = C$. Now
assertion~\ref{the: hyperbolic groups and fin to vcyc: maximal: V_max virt. cyc.} follows.

We conclude assertion \ref{the: hyperbolic groups and fin to vcyc: maximal: inclusions}
from the fact that for an inclusion of infinite virtually cyclic group $V \subseteq W$
there exists a normal infinite cyclic subgroup $C \subseteq W$ such that $C \subseteq V$
holds. Assertion \ref{the: hyperbolic groups and fin to vcyc: maximal: V_max = N_GV_max}
is now obviously true.
This finishes the proof of assertion 
\ref{the: hyperbolic groups and fin to vcyc: maximal}.
\\[2mm]
\ref{the: hyperbolic groups and fin to vcyc: maximal: model for calvcyc}
Construct a $G$-pushout
$$
\comsquare{\coprod_{i \in I} G \times_{V_i}\EGF{V_i}{\calfin}}{j}{\EGF{G}{\calfin}}
{\pr}{}
{\coprod_{i \in I} G/V_i}{}{X}
$$
with $j$ an inclusion of $G$-$CW$-complexes. Obviously $X$ is a $G$-$CW$-complex
whose isotropy groups are virtually cyclic. It remains to prove for
virtually cyclic $H \subseteq G$ that $X^H$ is contractible. 

Given a $V_i$-space $Y$ and a subgroup $H\subseteq G$, there is after a choice of a map of sets
$s \colon G/V_i \to G$, whose composition with the projection 
$G \to G/V_i$ is the identity, a  $G$-homeomorphism
\begin{eqnarray}
\coprod_{\substack{w \in G/V_i\\s(w)^{-1}Hs(w) \subseteq V_i}} 
Y^{s(w)^{-1}Hs(w)} 
& \xrightarrow{\cong} &
\left(G \times_{V_i} Y\right)^H,
\label{computation of (G times_V_i Y)^H}
\end{eqnarray}
which sends $y \in Y^{s(w)^{-1}Hs(w)}$ to $(s(w),y)$.

If $H$ is infinite, the $H$-fixed point set of the upper right and upper left corner is empty and 
of the lower left corner is the one-point space because of assertion
\ref{the: hyperbolic groups and fin to vcyc: maximal: V_max = N_GV_max} and equation
\ref{computation of (G times_V_i Y)^H}.
Hence $X^H$ is a point for an infinite virtually cyclic subgroup $H \subseteq G$.

If $H$ is finite, one checks using equation
\ref{computation of (G times_V_i Y)^H} that the left vertical map induces a homotopy 
equivalence on the $H$-fixed point set. Since the upper horizontal arrow induces a cofibration on the
$H$-fixed point set, the right vertical arrow induces a homotopy equivalence on the
$H$-fixed point sets. Hence $X^H$ is contractible for finite $H \subseteq G$.
This shows that $X$ is a model for $\EGF{G}{\calvcyc}$.
\\[2mm]
\ref{the: hyperbolic groups and fin to vcyc: maximal: H_n^g(vcyc,fin)}
follows from excision and the induction structure.
This finishes the proof of Theorem \ref{the: hyperbolic groups and fin to vcyc}. 
\end{proof}

Theorem~\ref{the: hyperbolic groups and fin to vcyc} has also been proved by Daniel Juan-Pineda
and Ian Leary \cite{Juan-Pineda-Leary(2003)} under the stronger condition
that every infinite  subgroup of $G$, which is not virtually cyclic, contains a non-abelian free subgroup.
The case, where $G$ is the fundamental group of a closed Riemannian
manifold with negative sectional curvature is treated in
\cite{Bartels-Reich(2003a)}. 

\begin{remark} \label{rem: Pindea-Leary} \em
In Theorem~\ref{the: hyperbolic groups and fin to vcyc} the terms 
$H_n^{V_i}(\EGF{V_i}{\calvcyc},\EGF{V_i}{\calfin};\bfK_R)$ and
$H_n^{V_i}(\EGF{V_i}{\calvcyc},\EGF{V_i}{\calfin};\bfL_R^{\langle -\infty \rangle})$
occur. They also appear in the direct sum decomposition
\begin{eqnarray*}
K_n(RV_i) & \cong & H_n^{V_i}(\EGF{V_i}{\calfin};\bfK_R) \bigoplus
H_n^{V_i}(\EGF{V_i}{\calvcyc},\EGF{V_i}{\calfin};\bfK_R);
\\
L_n(RV_i) & \cong & H_n^{V_i}(\EGF{V_i}{\calfin};\bfL_R^{\langle -\infty \rangle}) \bigoplus
H_n^{V_i}(\EGF{V_i}{\calvcyc},\EGF{V_i}{\calfin};\bfL_R^{\langle -\infty \rangle}).
\end{eqnarray*}
They can be analysed further and
contain information about and are build from the Nil and UNIL-terms in algebraic $K$-theory and $L$-theory
of the infinite virtually cyclic group $V_i$. They vanish for
$L$-theory after inverting $2$ by results of \cite{Cappell(1974b)}. 
For $R = \bbZ$ they vanishes rationally 
for algebraic $K$-theory by results of \cite {Kuku-Tang(2003)}.\em
\end{remark} 

\typeout{-------------------- References -------------------------------}

\addcontentsline{toc}{section}{References}
\bibliographystyle{abbrv}
\bibliography{dbdef,dbpub,dbpre,dbclassiextra}

\typeout{-------------------- Notation -------------------------------}
\twocolumn
\section*{Notation}
\addcontentsline{toc}{section}{Notation}


\typeout{--------------------   Notation --------------------}
\nopagebreak
\noindent
\entry{$\cd(G)$}{cd(G)}
\\
\entry{$\cd(M)$}{cd(M)}
\\
\entry{$CX$}{CX}
\\
\entry{$C_r^*(G)$}{C_r^*(G)}
\\
\entry{$\EGF{G}{\calf}$}{EGF(G)(calf)}
\\
\entry{$\underline{E}G$}{underline{E}G}
\\
\entry{$F_n$}{free group in n letters}
\\
\entry{$F\!P_n$}{FP_n}
\\
\entry{$F\!P_{\infty}$}{FP_infty}
\\
\entry{$G^0$}{G^0}
\\
\entry{$G_d$}{G_d}
\\
\entry{$\overline{G}$}{overline(G)}
\\
\entry{$\JGF{G}{\calf}$}{JGF(G)(calf)}
\\
\entry{$\underline{J}G$}{underline{J}G}
\\
\entry{l(H)}{l(H)}
\\
\entry{$N_GH$}{N_GH}
\\
\entry{$\Or(G)$}{Or(G)}
\\
\entry{$\OrGF{G}{\calf}$}{Or(G,calf)}
\\
\entry{$\Out(F_n)$}{outer group of free group}
\\
\entry{$PC_0(G)$}{PC_0(G)} 
\\
\entry{$P_d(G,S)$}{P_d(G,S)}
\\
\entry{$T = T(\calg,X,X_0)$}{T(calg,X,X_0)}
\\
\entry{$\vcd(G)$}{vcd(G)}
\\
\entry{$W_GH$}{W_GH}
\\
\entry{$[X,Y]^G$}{[X,Y]^G}
\\
\entry{$V\!F$}{VF}
\\
\entry{$\underline{\bbZ}$}{underline{Z}}
\\
\entry{$\pi = \pi(\calg,X,X_0)$}{pi(calg,X,X_0)}
\\
\entry{$\Gamma^s_{g,r}$}{Gamma^s_g,r}
\\
\entry{$\calall$}{calall}
\\
\entry{$\calcom$}{calcom}
\\
\entry{$\calcomop$}{calcomop}
\\
\entry{$\calfin$}{calfin} 
\\
\entry{$\calt^s_{g,r}$}{calt^s_g,r}
\\
\entry{$\caltr$}{caltr} 
\\
\entry{$\calvcyc$}{calvcyc}


\onecolumn

\typeout{-------------------- Index ---------------------------------}

\flushbottom
\addcontentsline{toc}{section}{Index}
\printindex                                  

\end{document}